\documentclass[12pt]{book}
\usepackage[hmarginratio=1:1]{geometry}
\usepackage{amsmath}
\usepackage{amsthm}
\usepackage{amssymb}
\usepackage{graphicx} 
\usepackage{multicol} 
\usepackage{enumerate} 
\usepackage[british]{babel}

\theoremstyle{plain}  
\newtheorem{theorem}{Theorem}[section]
\newtheorem{proposition}[theorem]{Proposition}
\newtheorem{corollary}[theorem]{Corollary}
\newtheorem{lemma}[theorem]{Lemma}

\theoremstyle{definition} 
\newtheorem{definition}[theorem]{Definition}
\newtheorem{example}[theorem]{Example}
\newtheorem{exercise}[theorem]{Exercise}

\theoremstyle{remark} 
\newtheorem{remark}[theorem]{Remark}

\newcommand{\R}{\mathbb{R}}  
\newcommand{\Q}{\mathbb{Q}}  
\newcommand{\Z}{\mathbb{Z}}  
\newcommand{\N}{\mathbb{N}}  
\newcommand{\C}{\mathbb{C}}  
\newcommand{\CP}{\mathbb{CP}}  
\newcommand{\ord}{\text{ord}\,} 
\newcommand{\codim}{\text{codim}} 
\renewcommand{\Im}{\text{Im}\,}
\renewcommand{\S}{S}
\newcommand{\pro}{\text{pro}}
\newcommand{\barbeta}{\bar{\beta}}
\newcommand{\cont}{\mathcal{O}}

\usepackage{txfonts} \usepackage[T1]{fontenc}

\usepackage{transparent}
\usepackage{eso-pic}
\AddToShipoutPicture*{
    \put(0,0){
        \parbox[b][\paperheight]{\paperwidth}{%
            \vfill
            \centering
            {\transparent{0.5}\includegraphics[width=0.8\textwidth]{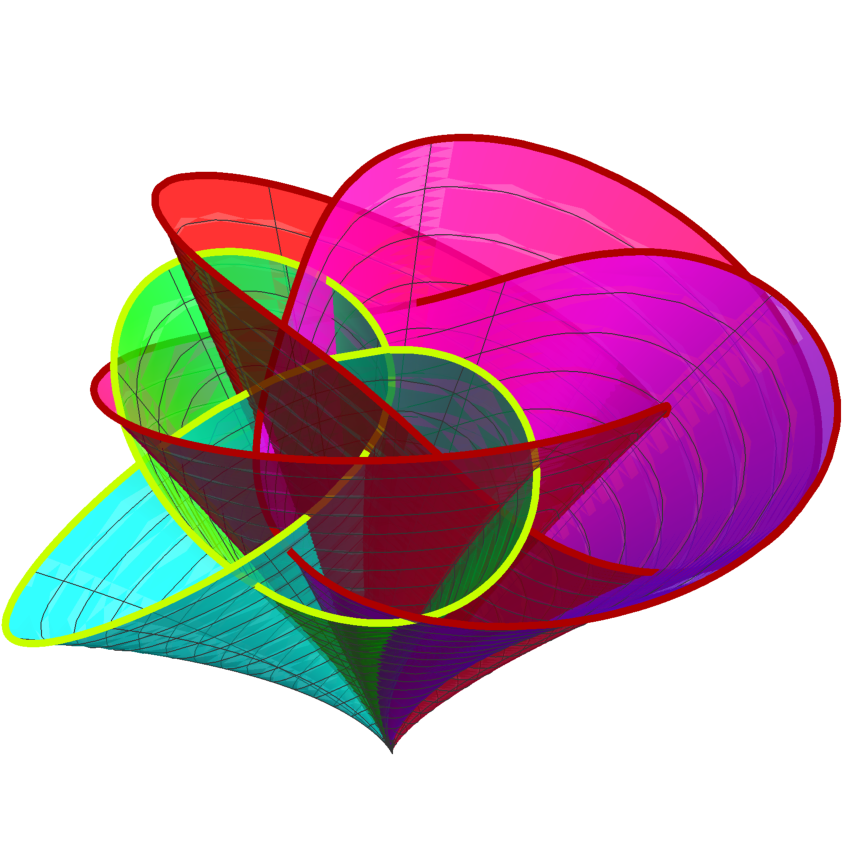}}%
            \vfill
        }
    }
}

\usepackage[ backend=biber, style=alphabetic, sorting=nyvt]{biblatex}
\bibstyle{alpha} 
\usepackage{csquotes}

\begin{document}
\pagenumbering{roman}
\begin{titlepage}
    \vspace{0.5in}
    \begin{center}
    \begin{figure}[ht]
    \begin{center}
    \includegraphics[width=\textwidth]{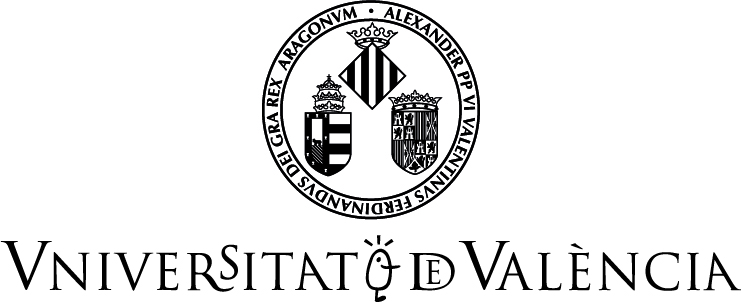}
    \end{center}
    \end{figure}
    \begin{huge}
    Degree Final Project - 2021/2022

\vspace{0.4in}    
    \textbf{On the topological classification of \\ complex plane curve singularities} \\
    \end{huge}
    
\vspace{0.8in}   
    \begin{huge}
    Author: {\bf
    Alberto Fern\'andez Hern\'andez} \\
    \end{huge}
    
    \vspace*{0.2in}
    
    \begin{huge}
    Advisor: {\bf \sc
    Juan Jos\'e Nu\~{n}o Ballesteros } \\
    \end{huge}
    \vspace{0.5in}
    \rule{110mm}{0.3mm}
    \begin{figure}[ht]
        \centering
        \includegraphics[width=\textwidth]{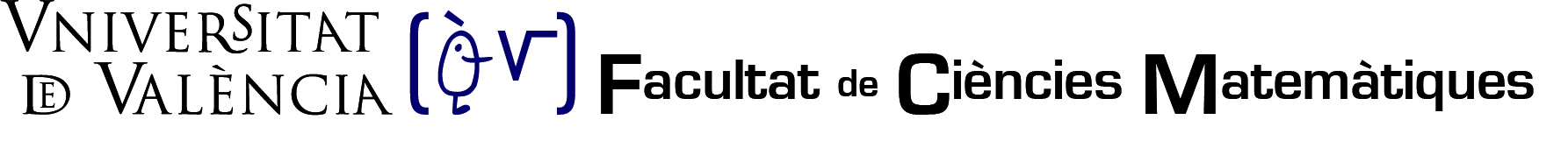}
    \end{figure}

    \end{center}
    \end{titlepage}
    
\newpage
\pagestyle{empty}
\textcolor{white}{blank page}
\newpage

\begin{flushright}
    “My work always tried to unite the truth with the beautiful, but when I had\\ 
    to choose one or the other, I usually chose the beautiful” - Hermann Weyl
\end{flushright} 
\vfill
\begin{center} \begin{Large}\textbf{Abstract} \end{Large}\end{center}
This project is devoted to study the topological classification of complex plane curves. These are subsets of $\C^2$ that can be described by an equation $f(x,y)=0$. Loosely speaking, curves are said to be equivalent in a topological sense whenever they are ambient homeomorphic, \textit{i.e.}, there exists an orientation-preserving homeomorphism of the ambient space carrying one curve to the other. The project's aim is to develop operative and clear conditions to determine whether two curves are equivalent. Curves will be shown to be decomposable into branches: sets that can be explicitly parametrised in the form $x=\phi(t), y=\psi(t)$. These parametric expressions will be analysed to extract a complete numerical invariant for the classification of branches: the Puiseux characteristic. The intermediate key result to lend this notion with topological weight is that a branch can be completely described through its associated knot, arising from the intersection of the branch with a small enough 3-sphere. The combination of these above-mentioned facts will then culminate in the project's most powerful result, which assures that two curves are equivalent if and only if their branches share the same Puiseux characteristics and intersection numbers.
\vspace{-0.3cm}
\begin{center}
    \includegraphics[scale=0.6]{curve2.eps}

    \vspace{-1cm}
    \parbox{0.85\textwidth}{{\bf A mathematical flower:} a representation of the curve $x^5-x^3y^3-x^2y^4+y^7=0$, which is given by the branches $x^2-y^3=0$ (green-blue) and $x^3-y^4=0$ (red-purple). Its link consists of two connected components: the $(3,2)$ and $(4,3)$ torus knots, painted yellow and red, respectively.}
\end{center}
\newpage
\thispagestyle{empty}
\textcolor{white}{blank page}
\newpage
\thispagestyle{empty}
\tableofcontents
\newpage
\pagestyle{empty}
\textcolor{white}{blank page}
\newpage
\pagenumbering{arabic}
\thispagestyle{plain}
\setcounter{page}{1}
\chapter*{Introduction}
\markboth{Introduction}{Introduction}
\addcontentsline{toc}{chapter}{Introduction}
The main purpose of this project is to study one of the most beautiful objects that geometry has ever provided to the mathematicians: complex plane curves. A plane curve is the zero set of a nonzero holomorphic function $f:U\subset \C^2\rightarrow \C$, where $U$ is an open neighbourhood of $(0,0)$ and $f(0,0)=0$. In other words, a curve is a subset of $\C^2$ that is defined by a nontrivial condition $f(x,y)=0$. Since our purpose is to analyse them locally, a curve would be considered as a set-germ $(C,0)\subset (\C^2,0)$. It will be shown that a curve can be decomposed in a finite number of branches, each of them corresponding to an irreducible factor of the defining function $f$. These branches can, in turn, be fully described by parametric equations of the form $x=t^m, y=\phi(t)$, where $\phi$ is a holomorphic function of a single variable, or just by $x=0,y=t$. In the first chapter, the formal treatment of these preliminary statements is analysed, as well as the theoretical groundwork required to develop the topological classification of plane curves. \\

A tuple of natural numbers can be extracted from each branch, namely the Puiseux characteristic, and it turns out to encapsulate a great amount of information from it. If two branches share the same Puiseux characteristic, they will be referred to be equisingular. In general, two curves will be said to be equisingular whenever the Puiseux characteristic and the intersection numbers of their branches coincide. The second chapter develops the above-mentioned equisingularity concept and provides some arithmetic tools that will play a fundamental role in the topological classification of branches and curves: the exponent of contact and the semigroup of a branch. Equisingularity offers an operative classification method for branches and curves, due to the fact that the Puiseux characteristic of a branch can be obtained through a routine calculation. However, its arithmetic nature does not offer a topological intuition of what do equisingular branches share in its topological structure.  \\

It turns out that the most fruitful way to distinguish the topology of curves is by looking at the way they are embedded in $\C^2$. More precisely, if $D_\epsilon\subset \C^2$ is the closed disk of radius $\epsilon$, the curves $C$ and $\tilde{C}$ will be said to be topologically equivalent whenever there exist orientation-preserving homeomorphisms $h:D_\epsilon\rightarrow D_\epsilon$ such that $h(C\cap D_\epsilon)=\tilde{C}\cap D_\epsilon$, for every $\epsilon>0$ small enough. In the third chapter, the main result of the project is presented and developed:
\begin{center}\textbf{ Two curves are equisingular if and only if they are topologically equivalent.}\end{center} 

\vspace{\baselineskip}
The conical structure of the curves will be shown to force the topological class of a curve to be fully described by the understanding of its link, which is the intersection of the curve with a small enough 3-dimensional sphere $\S^3$, the boundary of $D_\epsilon$. While, in general, this intersection is a link (a disjoint union of knots), for a branch this intersection defines a single knot. Moreover, parametric equations for a given branch induce parametric equations for the associated knot. In fact, the study of these equations yields that the isotopy type of the knot depends only on the Puiseux characteristic of the branch. Hence, equisingular branches define the same knot up to an isotopy. An analysis of the kind of knots that appear associated to the plane branches, namely the cable knots, enlightens that one can recover the Puiseux characteristic of a branch from the Alexander polynomial of the knot that defines in $\S^3$. This forces that the isotopy type of the knot also determines the equisingularity class of the branch. Therefore, a classification via equisingularity induces exactly the same classes as the topological equivalence for branches. With the tools introduced in the second chapter, this result is extended for general curves, and the main result of the project follows. \\

The exposition is strongly inspired in Wall's book \cite{Wall}, chapters 1, 2, 4 and 5. It has been intended to reduce the information provided \cite{Wall} in order to restrict to our purposes. Following the same reasoning, some concepts have been noticeably expanded. For this latter purpose, reference \cite{AnalyticGeometry} is being followed for the development of the fundamental theorems in Analytic Geometry, \cite{Juanjo} for some Singularity Theory results and concepts, together with \cite{burde} and \cite{cromwell} for a reference in Knot Theory. \\


Along the whole project, several solved exercises and examples are provided to illustrate the involved defins and theoretical results.
\chapter{Preliminaries}
The aim of this first chapter is to introduce the main concepts of the project, as well as to present some basic properties that will be useful later. During this work, our attention will be restricted to the understanding of curves defined through an equation of the form $f(x,y)=0$, or through parametric equations $x=\phi(t), y=\psi(t)$. Newton’s algorithm, studied in the second section, provides an effective method for computing parametric equations of the form $x=t^n, y=\phi(t)$ from a given power series $f(x,y)$. However, it does not deal with the convergence of $\phi$ provided $f$ is convergent. This will be solved in the third section, where an application of Weierstrass Preparation Theorem provides a non-constructive proof for the fact that $\C\{x,y\}$ is a unique factorisation domain, and that every curve that does not contain the germ of vertical line $x=0$ admits, at least, a parametric equation of the form $x=t^n, y=\phi(t)$, with $\phi(t)$ being a holomorphic function. With this, in the fourth section, branches are presented as the curves defined by an irreducible function. It is shown that a curve is a finite union of branches that just meet at the origin, and that branches are exactly the same as sets given by parametric equations. Finally, the concept of tangent lines and multiplicity of a curve are presented, as well as the intersection number between two branches.
\section{Concept of curve}
This work will be focusing primarily on the complex singularities in $\C^2$ near the origin. Consequently, the germ concept, presented below, will play a fundamental role.
\begin{definition} Let $x\in \C^n$. 
    \begin{enumerate}[(i)]
        \item Let $f:U\rightarrow \C^m$, $g:V\rightarrow \C^m$ be maps, where $U$ and $V$ are open neighbourhoods of $x$ in $\C^n$. We say that $f$ and $g$ \textit{define the same germ at} $x$ if there is a  neighbourhood $W\subset U \cap V$ of $x$ in $\C^n$ such that $f$ and $g$ coincide on $W$. This clearly defines an equivalence relation, and a \textit{germ of mapping at} $x$ is an equivalence class under this relation. 
        \item Two subsets $C$ and $D$ of $\C^n$ \textit{have the same germ at} $x$ if there is a neighbourhood $U$ of $x$ in $\C^n$ such that $C\cap U=D\cap U$. A \textit{set-germ at} $x$ is an equivalence class under this relation.
        \item We denote a germ at $x$ of a mapping $f:U\rightarrow \C^m$ by $f:(\C^n, x)\rightarrow \C^m$, or $f:(\C^n, x)\rightarrow (\C^m, y)$ so as to specify that $f(x)=y$. Analogously, the germ at $x$ of a subset $C$ is denoted by $(C,x)$. 
        \item A map-germ $f:(\C^n, x)\rightarrow \C^m$ is said to be \textit{continuous} whenever there exists a representative map of $f$ that is continuous at $x$. 
    \end{enumerate}
\end{definition}
\noindent This concept has the following basic properties:
\begin{proposition}\label{germs} Let $f:(\C^n,x)\rightarrow (\C^m,y)$ be a continuous germ of mapping. Then,
    \begin{enumerate}[(i)]
        \item For any set-germ $(D,y)$ in $\C^m$ the inverse image $f^{-1}(D,y)$ is well-defined as $(f^{-1}(D), y)$, that is, this set-germ does not depend on the representative of the map germ.
        \item If $f$ is open onto its image, then the image of any germ $(C,x)$ in $\C^n$ is well-defined as $(f(C),y)$.
        \item If $(f^{-1}(y),x)=(\{x\}, x)$, then $f$ is open onto its image.
    \end{enumerate}
\end{proposition}
\noindent For a proof, see Lemma 2.1 and Proposition 2.1 of \cite{Juanjo}. 
In what follows, the definition of a plane curve is presented.
\begin{definition}
    A \textit{curve} in $(\C^2,0)$ is a set-germ $(C,0)$ given by the inverse image of $0$ by a nonzero holomorphic function germ $f:(\C^2,0)\rightarrow (\C,0)$, \textit{i.e.}
    \begin{equation*}
        (C,0)=(f^{-1}(0), 0).
    \end{equation*}
\end{definition}
\noindent Notice that $0$ means $(0,0)$ when it appears in $(\C^2,0)$ and $(C,0)$. It will be common to shorten the previous concept just by calling it \textit{a curve }$(C,0)$ \textit{defined by} $f$, or \textit{defined by} $f(x,y)=0$. 
\begin{example}
    Let $f(x,y)=x^2+y^2$. Then, $f$ defines the curve $(C,0)$ given by the germ at $(0,0)$ of the set
    \begin{equation*}
        \{(x,y)\in \C^2 : x^2+y^2=0\},
    \end{equation*}
    which is the union of the line germs at $(0,0)$ of $y=ix$ and $y=-ix$. \\
\end{example}
As we stated in the previous definition, we are concerned about the curves defined by holomorphic functions in a neighbourhood of $(0,0)$. However, a precise definition of what ``holomorphic'' means applied to a function of more than one variable has not been given yet. We already know that a function $f:A\subset\C \rightarrow \C$ defined in an open neighbourhood of $0$ is holomorphic at the origin if and only if it is analytic at that point. That is, there exists a radius $R>0$ such that, if $|x|<R$, then $x\in A$ and $f(x)$ can be written as a convergent power series $f(x)=\sum_{n=0}^{+\infty}a_n x^n$. Therefore, a function $f:\Omega \subset \C^2 \rightarrow \C$ defined in an open neighbourhood of $(0,0)$ is \textit{holomorphic} at $(0,0)$ if there exists some radii $R,S>0$ such that, if $|x|<R$ and $|y|<S$, then $(x,y)\in \Omega$ and $f(x,y)$ can be written as
\begin{equation*}
    f(x,y)=\sum_{n,m=0}^{+\infty} a_{n,m}x^n y^m.
\end{equation*}
In general, a map $f:\Omega\subset \C^n\rightarrow \C^m$ defined in a neighbourhood of the origin is holomorphic at $0\in \C^n$ if there exists $R_1,\dots,R_n>0$ such that, for $|x_1|<R_1, \dots, |x_n|<R_n$, its coordinates $f=(f_1,\ldots, f_m)$ can be expressed as convergent power series in the variables $x_1, \ldots, x_n$ of $\C^n$. \\
\begin{remark} As the reader might realise, alternative notions for a map $f:\Omega\subset \C^n\rightarrow \C^m$ to be holomorphic appear along the literature. For instance, one could define this concept through the complex differential, or the generalised Cauchy-Riemann equations. As it happens in the single variable case, all of them provide equivalent definitions.
\end{remark}
\noindent In terms of germs, one says that the map-germ $f:(\C^n,0)\rightarrow \C^m$ is \textit{holomorphic} whenever there exists a representative map of $f$ that is holomorphic at $0$. \\

One could also think of a possible definition of the curve concept by means of parametric equations. 
\begin{definition}
    A \textit{parametric equation} for a curve $(C,0)$ is a germ of an injective holomorphic mapping $\alpha:(\C,0)\rightarrow (\C^2,0)$ whose image is contained in $(C,0)$. 
\end{definition}
\noindent Recall that a germ of mapping is said to be injective if there exists a representative mapping that is injective. That is, injectivity follows for a small enough neighbourhood of the basepoint. \\

Consequently, if $(C,0)$ is defined by the function $f:(\C^2, 0)\rightarrow (\C,0)$, then $\alpha(t)=(\phi(t), \psi(t))$ is a parametric equation for $(C,0)$ when $f(\phi(t),\psi(t))=0$ in a neighbourhood of $t=0$.\\
\begin{remark}
    Since $\alpha$ is injective, it follows by Proposition \ref{germs} that $\alpha$ is open onto its image, and therefore the germ of the image $\alpha(\C,0)\subset (\C^2,0)$ is well-defined.  
\end{remark}
\noindent The sets that can be fully described through parametric equations will be shown to be ``branches'', and it will be studied that any curve can be described as a finite union of branches.
\section{Newton's algorithm}
In this section, our objective is to find parametric equations for a curve given by $f(x,y)=0$. In order to do so, it will be fruitful to give an effective method to express $y$ in terms of $x$. This first approach does not deal with the convergence near the origin of the series that appear.
\begin{theorem}[Newton's algorithm]\label{thm:newton} Let $f(x,y)$ be a nonzero formal power series with $f(0,0)=0$. Then, the equation $f(x,y)=0$ admits at least one solution of the form
\begin{equation*}
    x=t^n, \hspace*{0.3cm} y=\sum_{r=1}^{+\infty} a_rt^r,
\end{equation*}
where $n\in \N$, or of the form $x=0, y=t$. 
\end{theorem}
\begin{remark}
    This algorithm lets us write the equation $f(x,y)=0$ expressing $y$ in terms of a ``fractional power series'' of $x$, namely $y=\sum_{r=1}^{+\infty} a_rx^{r/n}$. However, $n$-th roots in $\C$ have some peculiarities, as in general there is not a unique $n$-th root for a given complex number. In this first chapter the fractional exponents in a variable will be considered to be formal, that is, in the ring of fractional formal power series in variables $x^{1/n}$ and $y$. In the second chapter, a more natural approach to define properly $x^{1/n}$ by means of $n$-th roots will be given.
\end{remark}
\noindent The proof of Theorem \ref{thm:newton} is included in this project so as to explain the general procedure of Newton's algorithm. Hence, since the result is not going to be applied later, the proof will be merely sketched. 
\begin{proof}[Proof of Theorem \ref{thm:newton}] Let $f(x,y)=\sum_{r,s=0}^{+\infty} a_{r,s}x^ry^s$. The construction of the expansion of $y$ in terms of fractional powers of $x$ is done term by term, and the main idea is to start with the part of $f$ that ``dominates'' the rest. In order to clarify this notion, assume that $$y=c_0x^\alpha + \text{terms of higher order,}$$ where $\alpha\in \Q$. Then, each term $a_{r,s}x^ry^s$ of $f$ becomes $a_{r,s}c_0^sx^{r+s\alpha}$ plus terms of higher order, when performing the substitution of $y$. Our objective will be to make a choice of $c_0$ and $\alpha$ in such a way that the terms of order at most $r+s\alpha$ disappear after the substitution.
    
In order to do so, consider the real plane with coordinates $(r,s)$, and draw a point for the values $(r,s)\in \N^2$ such that $a_{r,s}\neq 0$. The set of these points will be called the \textit{Newton diagram}, and denoted as $\Delta (f)$. Now, let us focus on the polygonal line whose vertices are points of $\Delta(f)$ and whose edges leave the origin and the points of $\Delta(f)$ in different regions. More precisely, consider the convex hull of $\Delta(f)+(\R^+)^2$. Its topological border consists of two half lines parallel to the axes and a polygonal line joining them. The \textit{Newton polygon} is defined as this polygonal line (see Figure 1).
\begin{figure}
    \centering
    \includegraphics[width=0.38\textwidth]{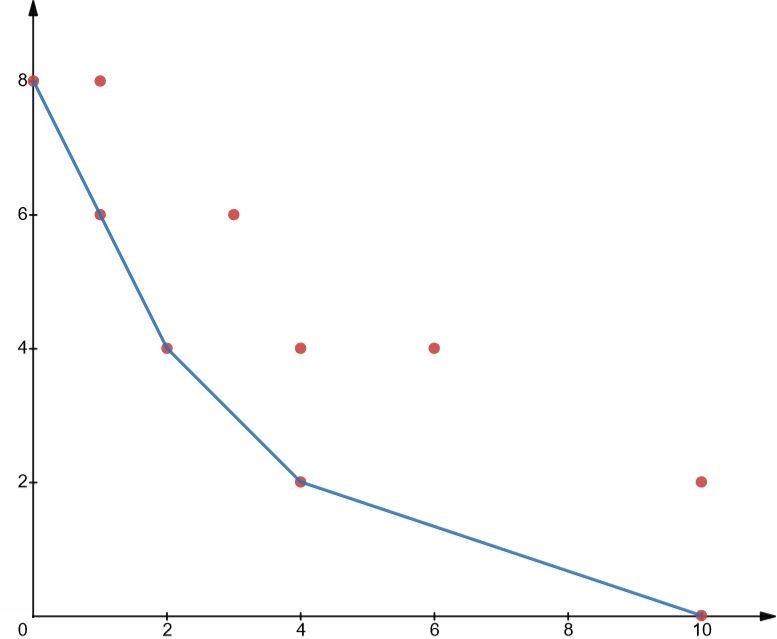}
    \caption{A Newton Polygon}
\end{figure}
Choose a line $r+s\alpha=C$ lying along an edge of this polygon, and write the vertices of the segment chosen as $(r_0,s_0)$ and $(r_0+ka, s_0-kb)$ where $a,b,k$ are positive integers with $a$ and $b$ being coprime. Accordingly, the points of the segment with integer coordinates have the precise form $(r_0+la,s_0-lb)$, where $l\in \{0,\dots,k\}$. Thus, after the substitution of $y=c_0x^{\alpha}+\dots$ in $f(x,y)$, no powers of $x$ with order lower than $C=r+s\alpha$ appear. This follows from the fact that there are no points $(r,s)$ ``under'' the Newton polygon. Moreover, the coefficient of $x^C$ would be $\sum_{l=0}^k a_{r_0+la,s_0-lb}c_0^{s_0-lb}$. Let us denote the involved coefficients as $u_l=a_{r_0+la,s_0-lb}$, and $\phi(T)=\sum_{l=0}^k u_l T^{k-l}$. Then, the coefficient of $x^C$ is $c_0^{s_0-kb}\phi(c_0^b)$. In addition, the endpoints of the segment are chosen so that $u_0$ and $u_k$ are nonzero. \\

With all the above-mentioned considerations, we are now able to start with the algorithm for the construction of the fractional power series of $y$ in terms of $x$. Let us check first whether $f$ is divisible by $x$: in this case, $x=0,y=t$ are parametric equations for $f(x,y)=0$. In general, write $f(x,y)=x^m\tilde{f}(x,y)$ where $m\in \N$, such that $x$ does not divide $\tilde{f}$, and work with $\tilde{f}$. Notice that the Newton polygon of $f$ will intersect the $s$-axis in the point $(0,m_0)$, where $m_0$ is the order of $f(0,y)$ (which is nonzero since $x$ does not divide $f$). Choose an edge of the Newton polygon, a solution $T_0$ of $\phi(T)=0$ and let $c_0$ be a $b$-th root of $T_0$. Notice that $T_0$ and $c_0$ are nonzero since $u_k\neq 0$ in the above discussion. Take $y=c_0x^{a/b}$ as a first approximation for $f(x,y)=0$. The choice of $c_0$ and $a/b$ forces $f(x, c_0x^{a/b})$ to have an order greater than $r_0+s_0a/b$. If the previous substitution yields the zero function as a result, then $y=c_0x^{a/b}$ is an exact solution for $f(x,y)$. Otherwise, it can be improved by adding extra terms. In order to improve the approximation, set $x=x_1^b$ and $y=x_1^a(c_0+y_1)$, and substitute in $f(x,y)$. It follows that the expression $ f(x_1^b, x_1^a(c_0+y_1))$ is divisible by $(x_1^b)^{r_0}(x_1^a)^{s_0}=x_1^{br_0+as_0}$. Then, divide by this term and denote the result as $f_1(x_1,y_1)$. Repeat the above procedure but replacing $f,x,y$ with $f_1,x_1,y_1$. Write $m_1$ as the order of $f_1(0,y_1)$ and, since $f_1(0,y_1)$ is obtained from the terms on the edge of the previous Newton polygon by dividing by $x_1^{br_0+as_0}$, then it is given by $(c_0+y_1)^{s_0-kb}\phi((c_0+y_1)^b)$. Since the degree of $\phi$ is $k$, one has that $m_1\leq kb\leq m_0$. \\

Repeat this procedure. After $r$ iterations, the resulting expression has the form 
\begin{equation*}
    Y_r=x^{\frac{a_0}{b_0}}\left( c_0+x^{\frac{a_1}{b_0b_1}} \left( c_1+x^{\frac{a_2}{b_0b_1b_2}}\left( c_2+ \ldots (c_r+x^{\frac{a_r}{b_0b_1\ldots b_r}})\ldots\right) \right) \right),
\end{equation*}
which can be manipulated to be an expansion of increasing fractional powers of $x$. Notice that $Y_r$ differs from $Y_{r-1}$ only in coefficients of $x$ with order of, at least, $a_r/(b_0b_1\ldots b_r)$. If the lowest nonzero coefficient of the substitution $f(x,Y_{r-1})$ is of order $M_{r-1}$, then, since the obtention of $Y_r$ is constructed to force these monomials of order $M_{r-1}$ to disappear, it follows that $M_r>M_{r-1}$. \\

We claim that, if $b_r>1$, then $m_r>m_{r+1}$. If this holds, since $m_r$ is a monotone decreasing sequence of positive integers, it must be constant at some point onwards. Then, the claim would yield that $b_r=1$ for large enough $r$. Thus, the $Y_r$ are all power series in $x_d=x^{1/d}$ for some fixed positive integer $d$. In addition, the series expansion of $Y_{r}$ and $Y_{r-1}$ agree up to the power $(x_d)^{dM_{r-1}}$, and the numbers $dM_{r}$ tend to infinity as $r$ does. Hence, the sequence of power series $Y_r$ for $r\in \N$ are constructed preserving their first $dM_{r-1}$ coefficients. One may thus define a power series $Y_\infty$ where each coefficient is defined as the ``limit'' of the coefficients of the $Y_r$. This is truly a trivial limit, since the sequence is eventually constant. This notion is known as a limit in the $\mathfrak{m}$-adic sense. Consequently, one has that $f(x,Y_\infty)=0$, and therefore $x=(x_d)^d, y=Y_\infty$ defines a parametric equation in power series for $f(x,y)=0$, as desired. \\

It only remains to check that $b_r>1$ implies that $m_r>m_{r+1}$. Since we are repeating the same procedure at each step, it suffices to consider the case $r=1$. Since $m_0\geq m_1$ in general, let us assume that $m_0=m_1$ for showing that $b:=b_1=1$. We may have that the least power of $y_1$ in 
\begin{equation*}
    (c_0+y_1)^{m_0-kb}\phi((c_0+y_1)^b)
\end{equation*}
has degree $m_0$. Thus, the previous polynomial is just a constant multiple of $y_1^{m_0}$. Since $c_0\neq 0$, $m_0=kb$ and then the Newton polygon has just one edge. Substitute $z=y_1+c_0$. This forces $(z-c_0)^{bk}$ to be equal to a polynomial in $z^b$. Thus, $b=1$, since the coefficient of $z$ in $(z-c_0)^{bk}$ is nonzero. This proves the claim.
\end{proof}
\begin{remark}
    Notice that this proof does not hold if, instead of working over $\C$, we do it over a field with finite characteristic. However, since this scenario is beyond the scope of this work, we will skip the details.\\
\end{remark}
It can be shown that, if $f$ is holomorphic, the Newton algorithm produces holomorphic equations. However, we are not going to prove this convergence claim. In practice, we will apply the algorithm to find an approximation of a parametric equation, and in case the algorithm stops, holomorphy will follow. Nonetheless, and for the sake of completeness, in the next chapter we will give a non-constructive proof showing that every curve defined by a holomorphic function $f$ has holomorphic parametric equations. \\

Now, an example is presented to illustrate the Newton's algorithm in a particular case.
\begin{example} Let $f(x,y)=x^4+3x^2y+x^2-y^3$. Then, the Newton polygon consists of a single edge joining the points $(0,3)$ and $(2,0)$. Thus, $(r_0,s_0)=(0,3)$, $k=1$ and $(a,b)=(2,3)$. We will then have $\alpha=2/3$ and $y=c_0x^{2/3}$ to be our first approximation. To obtain the value of $c_0$, notice that substitution in $f(x,y)$ yields a polynomial in $x$ of the form $(1-c_0^3)x^2$ plus terms of higher order. Then, any third root of the unity would be a valid option for $c_0$. Thus, choose $c_0=1$. Our first approximation to $f(x,y)=0$ would therefore be $y=x^{2/3}$. Now, for a second approximation, write $x=x_1^3$ and $y=x_1^2(1+y_1)$. Then, a substitution in $f(x,y)$ brings
\begin{equation*}
    x_1^{12}+3x_1^8(1+y_1)-x_1^6(3y_1+3y_1^2+y_1^3),
\end{equation*}
which is divisible by $x_1^6$ (as the theory anticipated, since $6=3\cdot 0 + 2\cdot 3=br_0+as_0$). Then, 
\begin{equation*}
    f_1(x_1,y_1)=x_1^6+3x_1^2+3x_1^2y_1-3y_1-3y_1^2-y_1^3.
\end{equation*}
We have that the Newton polygon of $f_1$ consists again of a single edge joining $(0,1)$ and $(2,0)$. Thus, $(r_1,s_1)=(0,1)$, $k=1$ and $(a,b)=(2,1)$. Therefore, $y_1=c_1x_1^2$ and a substitution in $f_1(x_1,y_1)$ yields the term of least order to be $3(1-c_1)x_1^2$. Hence, $c_1=1$ and the second approximation for $f(x,y)=0$ is $y_1=x_1^2$, which yields $x=x_1^3$ and $y=x_1^2(1+x_1^2)=x_1^2+x_1^4$. For a third approximation, we should write $x_1=x_2$ and $y_1=x_2^2(1+y_2)$. However, notice that, if $y_1=x_1^2$ one has that $f_1(x_1,y_1)=0$. Therefore, the addition of $y_2$ would immediately yield $y_2=0$. This is a particular case where the algorithm ends in a finite number of steps, giving the exact solution $x=t^3, y=t^2+t^4$. 
\end{example}
\begin{remark}
    Notice that, in general, the algorithm would not stop in a finite number of steps. Thus, in practice, one has to settle for an approximation with a certain order. Nonetheless, approximations will be shown to be enough  to extract the essential topological information that is encoded in a curve: the Puiseux characteristic. This will be detailed in the following chapter.
\end{remark}
\section{Convergence}
In this section, it is intended to show that, if the function $f$ that defines the equation $f(x,y)=0$ is holomorphic at $(0,0)$ (that is, if the power series $f(x,y)$ is convergent in a neighbourhood of the origin), then one can write $y$ as a convergent power series of $t$, where $x=t^m$ for some positive integer $m$. Some of the results that are about to appear along this chapter can be generalised for functions and power series in $n$ variables. Let us start with some technical results that will be applied later. Firstly, recall that the order of a power series $f(t)=\sum_{n=0}^{+\infty}a_nt^n$ is the least $m\in\N$ such that $a_m\neq 0$. 
\begin{lemma}\label{lemma:order}
    If $f:(\C,0)\rightarrow (\C,0)$ is a holomorphic function of order $n$, then there exists a function $g:(\C,0)\rightarrow (\C,0)$ of order 1 such that $g^n=f$. 
\end{lemma}
\begin{proof} Since $f$ has order $n\geq 1$, we can write $f=ax^n(1+E)$ for $a\neq 0$ and $E$ a holomorphic function that vanishes at $0$. Let $b\in \C$ such that $b^n=a$, and expand $E=\sum_{r=1}^{+\infty} E_rx^r$. Notice that we have an expansion for $(1+E)^{1/n}$ given by 
    \begin{equation*}
        (1+E)^{n^{-1}}=1+n^{-1}E+\dfrac{n^{-1}(n^{-1}-1)}{2}E^2+ \dfrac{n^{-1}(n^{-1}-1)(n^{-1}-2)}{3!}E^3+\dots  
    \end{equation*}

    Notice that for a fixed $r\in \N$, there are only a finite number of terms contributing in the previous expansion to the coefficient of $x^r$. Thus, the previous defines a power series $G(x)$. If we check that $G$ is in fact a convergent power series, then the function $g=bxG(x)$ has clearly order $1$ since $G(0)=1$ and $b\neq 0$, and it satisfies that $g^n=ax^n(1+E)=f$. \\

    It only remains to check the convergence of $G$. Notice that, since $E$ is holomorphic and $E(0)=0$, we have that there exists a radius $R>0$ for which $\sum_{r=1}^{+\infty}|E_r|R^r<1$. Now, the binomial expansion has radius 1, and hence for $|x|<R$ the substitution of $E(x)$ in the binomial expansion yields a convergent power series, as it is the composition of holomorphic functions. Thus, $G(x)$ is holomorphic for $|x|<R$. 
\end{proof}
\begin{definition}
    A holomorphic function germ $f:(\C^2,0)\rightarrow (\C,0)$ is said to be \textit{regular in} $y$ \textit{of order} $m$ if we can write 
    \begin{equation*}
        f(0,y)=y^mA(y),
    \end{equation*}
    with $A$ being a holomorphic function with $A(0)\neq 0$ and $m$ being a positive integer. 
\end{definition}
\begin{remark}
    Notice that, for any function $f$, unless $f(0,y)=0$ in a neighbourhood of $y=0$, we may write $f(0,y)$ as in the previous definition for some $m\in \N$. \\
\end{remark}
The following result reduces the arguments concerning holomorphic functions to polynomials, and plays a fundamental role in the theoretical groundwork of plane curves. 
\begin{theorem}[Weierstrass Preparation Theorem]
    Let $G:(\C^2,0)\rightarrow (\C,0)$ be a holomorphic function regular in $y$ of order $s$. Then, there exists a holomorphic function $U:(\C^2,0)\rightarrow \C$ with $U(0,0)\neq 0$ and $A_r:(\C,0)\rightarrow (\C,0)$ holomorphic functions for $r\in \{0,\dots, s-1\}$ such that
    \begin{equation*}
        G(x,y)=U(x,y)\left(y^s+\sum_{r=0}^{s-1}A_r(x)y^r\right).
    \end{equation*}
    In addition, the functions $U$ and $A_r$ are unique with this property.
\end{theorem}
\noindent The proof of the theorem requires first to develop a particular division algorithm, and for the sake of brevity, the proofs of these results are omitted. We refer the reader to \cite{AnalyticGeometry}, Theorems 3.2.3 and 3.2.4 for a proof. The kind of polynomials that arise in Weierstrass Preparation Theorem are called \textit{Weierstrass polynomials}. On the other hand, the previous results let us perform the division in the usual way for holomorphic functions (see \cite{Wall}, Corollary 2.2.4):
\begin{corollary}\label{cor:weierstrass} Let $F,G:(\C^2, 0)\rightarrow \C$ be holomorphic functions. If $G(x,y)$ is regular in $y$ of order $s$, then 
    \begin{align*}
        F(x,y)=G(x,y)D(x,y)+\sum_{i=1}^s C_i(x)y^{s-i}.\\
    \end{align*}
\end{corollary}
Notice that the set of holomorphic functions $f:(\C^2,0)\rightarrow \C$ has the structure of a ring (this ring can be also understood as the set of convergent complex power series in $x$ and $y$). As it is usual, this ring will be denoted as $\C\{x,y\}$. Since a unit in the ring is an element with a multiplicative inverse, it follows that the units are precisely the functions $U(x,y)$ with $U(0,0)\neq 0$. With Weierstrass Preparation Theorem, we are now able to show that the ring of holomorphic functions is a unique factorisation domain, that is, an integral domain where any nonzero element can be expressed as a product of irreducible factors (recall that an irreducible element is a noninvertible element such that, whenever it is expressed as a product of two factors, one of them is necessarily a unit), and where irreducible elements are prime (\textit{i.e.} if an irreducible element divides a product, then it divides at least one of the factors).
\begin{theorem}
    $\C\{x,y\}$ is a unique factorisation domain.
\end{theorem} 
\begin{proof} 
First of all, we should check that $\C\{x,y\}$ is an integral domain. In order to do so, notice that an elementary result of algebra states that $\C[x,y]$ (the ring of polynomials in $x,y$) is an integral domain since so is $\C$. With this, we can show that the ring of power series in $x,y$, denoted as $\C[[x,y]]$, is also an integral domain. Since $\C\{x\}$ and $\C\{x,y\}$ are subrings of $\C[[x,y]]$, they will also be integral domains. Let us show then that $\C[[x,y]]$ is an integral domain. Let $A(x,y)=\sum_{n,m=0}^{+\infty} a_{n,m}x^ny^m$ and $B(x)=\sum_{n,m=0}^{+\infty} b_{n,m}x^ny^m$. Notice that a series is identically zero if and only if all the polynomials obtained by cutting the series up to a certain degree are zero. Therefore, if $A$ and $B$ were nonzero series, there should exist $k_1,k_2\in \N$ for which $\sum_{n,m=0}^{k_1}a_{n,m}x^ny^m$ and $\sum_{n,m=0}^{k_2}b_{n,m}x^ny^m$ are nonzero. Thus, if $k=k_1+k_2$, one has that the product series $A(x,y)B(x,y)=\sum_{n,m=0}^{+\infty}c_{n,m}x^ny^m$ satisfies that $$\sum_{n,m=0}^{k}c_{n,m}x^ny^m=\left(\sum_{n,m=0}^{k_1}a_{n,m}x^ny^m\right)\left(\sum_{n,m=0}^{k_2}b_{n,m}x^ny^m\right).$$
Since both factors are nonzero and $\C[x,y]$ is an integral domain, it follows that the previous polynomial is nonzero. Thus, the product series $A(x,y)B(x,y)$ is nonzero. This proves that $\C[[x,y]]$ is an integral domain, and then so are $\C\{x\}$ and $\C\{x,y\}$. \\

Now, let us show that $\C\{x,y\}$ is a unique factorisation domain (UFD). Notice first that $\C\{x\}$ is clearly a UFD, since it is, in particular, a principal ideal domain. In addition, Gauss' Lemma claims that, if $A$ is a UFD, then so is $A[t]$. Thus, the ring $\C\{x\}[y]$ is a UFD. Let $f\in \C\{x,y\}$ be a nonzero element. Then, one can write $f=x^k\tilde{f}$ for some natural number $k$ such that $x$ does not divide $\tilde{f}$ (this is straightforward in the ring of power series, and, consequently, also in the convergent ones). Then, $\tilde{f}$ is regular in $y$ of some order $s$, and Weierstrass Preparation Theorem lets us write $\tilde{f}=Ug$ where $U$ is a unit in $\C\{x,y\}$ and $g\in \C\{x\}[y]$ is a monic polynomial in $y$ of degree $s$. Hence, $g$ can be expressed as a product of irreducible factors $g_i$ in $\C\{x\}[y]$. We claim that these factors are also irreducible in $\C\{x,y\}$. Indeed, let $g_i=hk$ for $h,k\in \C\{x,y\}$. Then, $x$ does not divide $h$ or $k$, since otherwise $x$ would divide $g_i$ and hence $x$ would also divide $\tilde{f}$. Thus, $h$ and $k$ are regular in $y$ and one can apply Weierstrass Preparation Theorem to infer that $h=U_h \tilde{h}$ and $k=U_k \tilde{k}$ for $U_h,U_k$ units in $\C\{x,y\}$ and $\tilde{h},\tilde{k}\in \C\{x\}[y]$. Then, $g_i=U_hU_k \tilde{h}\tilde{k}$ and, by the uniqueness clause in Weierstrass Preparation Theorem, one has that $g_i=\tilde{h}\tilde{k}$. However, $g_i$ is irreducible in $\C\{x\}[y]$, and hence $\tilde{h}$ or $\tilde{k}$ is a unit. Thus, $g_i$ is also irreducible in $\C\{x,y\}$. This shows that $f$ is a product of irreducible functions of $\C\{x,y\}$.  \\

Finally, suppose that $g$ is an irreducible factor dividing a product $h_1h_2$ in $\C\{x,y\}$. If $x$ divides $g$, then $g(x,y)=x$, and in this case it is clear that $x$ has to divide $h_1$ or $h_2$. Indeed, if $g(0,y)=0$, then $h_1(0,y)h_2(0,y)=0$. Then, since $\C\{y\}$ is a domain, at least one of the factors is zero, so $x$ divides it. Hence, $x$ divides $h_1$ or $h_2$. By Corollary \ref{cor:weierstrass}, there exists $R_1,R_2\in \C\{x\}[y]$ such that $h_1=gQ_1+R_1$ and $h_2=gQ_2+R_2$. Since $g$ divides $h_1h_2$, it also divides $R_1R_2$, but irreducibles are prime in $\C\{x\}[y]$, yielding that $g$ divides at least one of the $R_i$. In this case, it follows that $g$ would divide $h_i$, so irreducible elements are prime in $\C\{x,y\}$.
\end{proof}
\noindent As we claimed before, the previous theorems can be extended to $\C\{x_1,\ldots, x_n\}$ without further complications. In fact, the reference \cite{AnalyticGeometry} proves Weierstrass Preparation Theorem in its general version. However, we are mainly interested in the case of $\C\{x,y\}$ along this exposition. \\

The following concepts, first stated in the XIX century, will be needed in the proof of this chapter's main theorem. Let $U$ be a unique factorisation domain with characteristic 0 (which will be, in our case $U=\C$ or $U=\C\{x\}$). For $f(y)=a_0+a_1y+\ldots + a_my^m$ and $g(y)=b_0+b_1y+\ldots +b_ny^n$, define the \textit{resultant} of $f$ and $g$ as the determinant 
\begin{equation*}
    \begin{vmatrix}
        a_0 & 0 & \cdots & 0 & b_0 & 0 & \cdots & 0 \\
        a_1 & a_0 & \ddots & 0 & b_1 & b_0 & \cdots & 0 \\ 
        a_2 & a_1 & \ddots & 0 & b_2 & b_1 & \ddots & 0 \\
        \vdots & \vdots & \ddots & a_0 & \vdots & \vdots & \ddots & b_0 \\ 
        a_m & a_{m-1} & \cdots & \vdots & b_n & b_{n-1} & \cdots & \vdots \\
        0 & a_m & \ddots & \vdots & 0 & b_n & \ddots & \vdots \\
        \vdots & \vdots & \ddots & a_{m-1} & \vdots & \vdots & \ddots & b_{n-1} \\ 
        0 & 0 & \cdots & a_m & 0 & 0 & \cdots & b_n \\
    \end{vmatrix}
\end{equation*}
of dimension $m+n$, and denote it as $R(f,g)$. It can be shown that $R(f,g)=0$ if and only if $f$ and $g$ have a nontrivial common factor in $U[y]$. \\

For a polynomial $f\in U[y]$ given by $f(y)=a_0+a_1y+\ldots+ a_my^m$, define its formal derivative as the polynomial $f'(y)=a_1+2a_2y+\ldots + ma_my^{m-1}\in U[y]$. With this, the \textit{discriminant} of $f$ is the element $D(f)=(-1)^{m(m-1)/2}R(f,f')\in U$. It can be shown that $D(f)\neq 0$ if and only if $f$ has no multiple roots in $U$. Moreover, a root of $f$ is multiple if and only if it is a root of $f$ and $f'$. Hence, if $h$ is the highest common factor between $f$ and $f'$, it turns out that $f/h$ is a polynomial in $U[y]$ with the same roots as $f$, and all of them are simple. For an in-depth exposition of these notions and results, see \cite{Resultant}, Appendix 1. \\

The following result can be considered to be the most fundamental theorem of the chapter, and concerns the existence of parametric equations for a curve. The relevance of this result is that a similar statement does not hold in the general case of $n$ variables. Namely, hypersurfaces in $\C^n$ described by an expression $f(x_1,\dots,x_n)=0$ do not admit, in general, parametric equations of the form $\gamma:(\C^{n-1},0)\rightarrow (\C^n,0)$. Plane curves are thus a beautiful exception, where richer conclusions could be achieved thanks to this phenomenon.
\begin{theorem}\label{th:par_eq} Let $(C,0)$ be a curve. Then, unless the germ of $x=0$ is contained in $(C,0)$, there exists a positive integer $n$ and a holomorphic function $\phi:(\C,0)\rightarrow (\C,0)$ such that $x=t^n, y=\phi(t)$ is a parametric equation for $(C,0)$.
\end{theorem}
\begin{proof} Let $f:(\C^2,0)\rightarrow (\C,0)$ be the holomorphic function that defines $f$. If $f(0,y)$ vanishes identically in a neighbourhood of $0$, then $(C,0)$ contains the vertical line germ of $x=0$. In other case, $f$ will be regular in $y$ of some order $s\geq 1$. Weierstrass Preparation Theorem lets us write $f(x,y)=U(x,y)A(x,y)$, where $A(x,y)=y^s+\sum_{r=0}^{s-1}A_r(x)y^r$ and $A_r:(\C,0)\rightarrow (\C,0)$ are holomorphic functions, and $U$ is a unit of $\C\{x,y\}$. Thus, since $U(0,0)\neq 0$, it follows that $A$ and $f$ define the same germ at $(0,0)$. \\

Now, $A\in \C\{x\}[y]$ is a polynomial in $y$ with coefficients in the UFD $\C\{x\}$, and hence it has a discriminant $D(A)\in \C\{x\}$. If $D(A)$ is identically zero in $\C\{x\}$, then one has that $A$ has a multiple root (\textit{i.e.} it has a factor $(y-a(x))^k$ for some $a(x)\in \C\{x\}$ and $k\geq 2$). In that case, if $h$ is the greatest common divisor between $A$ and $\partial A/\partial y$ in $\C\{x\}[y]$, then $A/h$ is a polynomial in $y$ with the same roots as $A$, and no repeated factors. Thus, relabelling $A$ with $A/h$ if necessary, we can assume that $A$ has no multiple roots as a polynomial in $y$, and hence $D(A)$ is not the zero function. \\

Notice that $D(A)(x)$ is a holomorphic function of $x$ near the origin, and that $D(A)(0)$ can be 0 because $A(0,y)=y^s$ (then, for $s>1$, we will have that $D(A)(0)=0$). However, $D(A)$ is nonzero as an element in $\C\{x\}$, so there exists an $\epsilon>0$ such that for $0<|x|<\epsilon$ one has that $D(A)$ is well-defined and $D(A)(x)\neq 0$. In particular, for a fixed $x_0\in \C$ with $0<|x_0|<\epsilon$, it holds that $D(A)(x_0)\neq 0$, and therefore the polynomial $A(x_0,y)\in \C[y]$ has distinct roots $y_{1,x_0}, \ldots, y_{s,x_0}\in \C$. In fact, for $i\in \{1,\ldots, s\}$ we have that $A(x_0,y_{i,x_0})=0$ and, in addition, $\partial A / \partial y \, (x_0,y_{i,x_0})\neq 0$ (otherwise, $y_{i,x_0}$ would be a multiple root of $A(x_0,y)$ and thus $D(A)(x_0)=0$). Therefore, the Implicit Function Theorem grants us that there exists a neighbourhood $U_{x_0}$ of $x_0$ in the region of points with $0<|x|<\epsilon$, and a holomorphic function $y_i:U_{x_0}\rightarrow \C$ such that $y_i(x_0)=y_{i,x_0}$ and, for every $x\in U_{x_0}$ and $y\in \C$, we have that $A(x,y)=0$ if and only if $y=y_j(x)$ for some $j\in \{1,\ldots, s\}$. \\

Up to this point we have defined the $y_1, \ldots, y_s$ locally, that is, in small neighbourhoods of a point $x$ whenever $0<|x|<\epsilon$. Hence, it would be of our interest to extend them to functions defined in the whole punctured ball. However, it is a matter of homotopy that this extension can be performed or not depending on the fundamental group of the domain. Since in our case the punctured ball $B_{\epsilon}^*=\{x\in \C: 0<|x|<\epsilon\}$ is not simply connected, it would be impossible to perform such an extension. Thus, our approach will be to remove from $B_{\epsilon}^*$ the segment $(0,\epsilon)$ to assure that the domain is a simply connected space. \\

In order to formalise the previous ideas, notice that $C^*=\{(x,y)\in B_{\epsilon}^* \times \C : A(x,y)=0\}$ is a representative set for the set-germ $(C\setminus\{0\},0)$, and the projection to the first coordinate $\pi :C^*\rightarrow B_{\epsilon}^*$ is a covering map. Indeed, for a fixed $x_0\in B_{\epsilon}^*$ we have shown the existence of a neighbourhood $U_{x_0}\subset B_{\epsilon}^*$ such that for $(x,y)\in U_{x_0}\times \C$ one has that $A(x,y)=0$ if and only if $y=y_i(x)$ for some $i\in \{1,\ldots, s\}$. Thus, if we define $\tilde{U_i}=U_{x_0}\times y_i(U_{x_0})\subset C^*$, since the $y_i(x_0)$ are all different, the continuity of the $y_i$'s lets us shrink $U_{x_0}$ if necessary so as to force the $y_1(U_{x_0}),\dots,y_s(U_{x_0})$ to be disjoint. Thus, $\pi^{-1}(U_{x_0})=\tilde{U}_1\cup \ldots\cup \tilde{U}_s$ is a union of disjoint open sets, and in addition, the restrictions $\pi:\tilde{U}_i\rightarrow U_i$ are homeomorphisms since the map $x\mapsto (x,y_i(x))$ is continuous and is the inverse of $\pi|_{\tilde{U_i}}$. Hence, $\pi : C^*\rightarrow B_{\epsilon}^*$ is a covering map.\\

With this, since $\Omega = B_\epsilon \setminus [0,\epsilon)$ is a path-connected and locally path-connected topological space with trivial fundamental group, by applying Proposition 1.33 of \cite{Hatcher}, one has that for the basepoint $x_0=\epsilon/2\,i \in \Omega$, the inclusion map $h:\Omega \rightarrow B_{\epsilon}^* $ has lifts $\hat{h}_j:\Omega \rightarrow C$ with $\hat{h}_j(x_0)=(x_0,y_{j,x_0})$, for $j\in \{1,\ldots, s\}$. Hence, if we conveniently write $\hat{h}_j(x)=(x,y_j(x))$, then the $y_j$ are well-defined functions in $\Omega$ with $A(x,y_j(x))=0$ for every $x\in \Omega$. In addition, the $y_j$ are unique with this property as a consequence of Proposition 1.34 of \cite{Hatcher}. Moreover, the previous extension can be performed \textit{mutatis mutandis} over the simply connected domain $\tilde{\Omega}=B_\epsilon \setminus (-\epsilon,0]$ and the same basepoint $x_0$, yielding $s$ functions $\tilde{y}_1, \ldots, \tilde{y}_s:\tilde{\Omega}\rightarrow \C$ such that $\tilde{y}_j(x_0)=y_{j,x_0}$ and with the property that they are the unique functions that satisfy $A(x,\tilde{y}_j(x))=0$ for $x\in \tilde{\Omega}$. \\

The following problem would be to fill the gap we have induced by cutting the punctured ball. Firstly, notice that $\Omega \cap \tilde{\Omega}=L^+\cup L^-$, where $L^+=\{x\in B_\epsilon : \Im x >0\}$ and $L^-=\{x\in B_\epsilon : \Im x <0\}$. In addition, the chosen basepoint $x_0=\epsilon/2\, i$ lies in $L^+$. Hence, $y_j$ and $\tilde{y}_j$ coincide in $L^+$ for every $j\in \{1,\ldots, s\}$. Furthermore, the restrictions of $y_1,\ldots, y_s$ to $L^-$ are the unique functions in $L^-$ that satisfy the condition $A(x,y(x))=0$ for $x\in L^-$. Since the same holds for the restrictions of the $\tilde{y}_j$ in $L^-$, there has to exist a permutation $\sigma \in \Sigma_s$ such that $\tilde{y}_j = y_{\sigma(j)}$ in $L^-$, for every $j\in \{1,\ldots, s\}$. \\

It is a well-known fact that the elements of the symmetric group $\Sigma_s$ split in disjoint cycles. Hence, relabelling the elements of $\{1,\dots,s\}$ if necessary, we may write $\sigma(1)=2, \sigma^2(1)=3, \ldots, \sigma^{n-1}(1)=n, \sigma^{n}(1)=1$, for some $n\leq s$ being the least nonzero natural number for which the condition $\sigma^n(1)=1$ holds. With this, it turns out that, for every $j\in \{1,\dots,n\}$, we have that $\tilde{y}_j=y_j$ in $L^+$, and that $\tilde{y}_j=y_{j+1}$ in $L^-$ (accepting by convenience that $y_{m+1}=y_1$). \\

Our last task is to define a holomorphic function $\phi$ in the punctured ball $B_{\delta}^*$ by means of the functions $y_1, \ldots, y_n$. In order to do so, divide the ball in $n$ sectors 
\begin{equation*}
    S\!_j=\left\lbrace re^{i\theta}: 0<r<\delta , \theta \in \left(\dfrac{2\pi j}{n}, \dfrac{2\pi (j+1)}{n}\right)\,\right\rbrace
\end{equation*}
for $j\in \{1,\dots,n\}$. Since the map $z\mapsto z^n$ sends $S\!_j$ to $\Omega=B_\epsilon\setminus [0,\epsilon)$ if $\delta=\epsilon^{1/n}$ for every $j\in \{1,\dots, n\}$, we could think about defining $\phi(t)=y_j(t^n)$ for $t\in S\!_j$, as it is shown in Figure 2. However, this does not give a complete description of $\phi$ in $B_{\delta}$, since the open sectors do not cover the whole punctured ball. In order to do so, notice that any sector $\tilde{S}\!_j=\lbrace re^{i\theta}: 0<r<\delta , \theta \in \left(2\pi (j-1/2)/n,2\pi (j+1/2)/n\right)\,\rbrace$ lies exactly in the middle of the consecutive sectors $S\!_j$ and $S\!_{j+1}$, for $j\in \{1,\dots,n\}$ (accepting again that $S\!_{n+1}:= S\!_1$). In addition, $z\mapsto z^m$ sends $\tilde{S}\!_j$ to $\tilde{\Omega}=B_\epsilon \setminus (-\epsilon,0]$ if $\delta=\epsilon^{1/n}$. Therefore, define $\phi:B_{\delta}^*\rightarrow \C$ as $\phi(t)=y_j(t^n)$ for $t\in S\!_j$ and $\phi(t)=\tilde{y}_j(t^n)$ if $t\in \tilde{S}\!_j$. It is clear that $\phi$ is a well-defined holomorphic function, since in the intersections of different sectors, both definitions coincide. Indeed, the map $z\mapsto z^m$ sends $\tilde{S}\!_j\cap S\!_j$ to $L^+$ and $\tilde{S}\!_j\cap S\!_{j+1}$ to $L^-$, and, in addition, $y_j=\tilde{y}_j$ in $L^+$ and $\tilde{y}_j=y_{j+1}$ in $L^-$. \\

\begin{figure}
    \centering
    \includegraphics[scale=0.39]{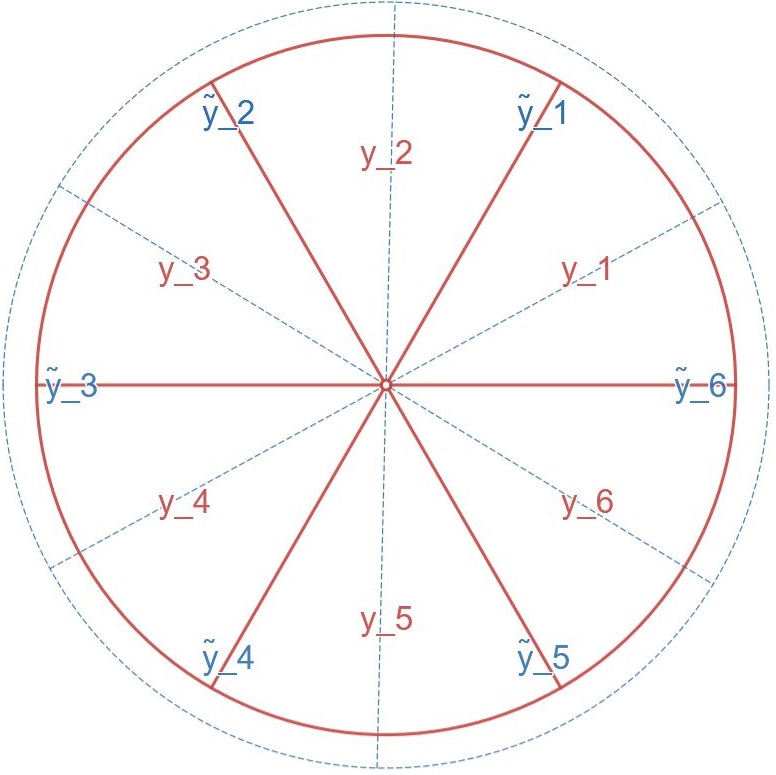}

    \caption{Representation of $\phi$ with respect to the $y_j$ (red) and $\tilde{y}_j$ (blue) for $n=6$. }
    \end{figure}
It then follows that $\phi$ is a holomorphic function in $B_\delta^*$ that satisfies $A(t^n, \phi(t))=0$ as a consequence of having $A(t^n, y_j(t^n))=0$ for $t\in S\!_j$, and the equivalent statement for $\tilde{S}\!_j$. Therefore, it just remains to extend $\phi$ to the origin. However, notice that the only root of $A(0,y)=y^s$ is precisely $y=0$. Hence, since the roots of a Weierstrass polynomial depend continuously on its coefficients (see Lemma 3.4.11 of \cite{AnalyticGeometry}), it follows that $\lim_{t\to 0}\phi(t)=0$ as a limit in $\C$. \\

Thus, $\phi$ is a holomorphic function in a punctured disk, with limit equal to 0 as $t\to 0$. By Riemann's removable singularity Theorem it follows that $\phi$ can be extended to a holomorphic function $\phi:B_{\delta}\rightarrow \C$ such that $\phi(0)=0$ and, for every $t\in B_\delta$, one has that $A(t^n, \phi(t))=0$. Therefore, $x=t^n,y=\phi(t)$ define parametric equations for $A$, and hence for the curve $(C,0)$.
\end{proof}
\begin{remark}
    Notice that the previous theorem grants us that any curve $(C,0)$ admits parametric equations, but in general the image of these equations will not be the whole curve: it is going to be contained in $(C,0)$. For a general curve, it would be required to have different parametric equations to cover the curve. This will be clarified later with the \textit{branch} concept. \\
\end{remark}
Notice that, in the previous proof, we have not checked whether the parametric equations are injective. It turns out that, in general this may not happen, but we can modify the equations slightly so as to force injectivity. 
\begin{theorem}\label{thrm:alwaysinjective}
    Let $\gamma:(\C,0)\rightarrow (\C^2,0)$ given by $\gamma(t)=(t^m, \phi(t))$ for $\phi(t)=\sum_{r=1}^{+\infty}a_rt^r$ a holomorphic function. Then, $\gamma$ is injective if and only if 
    \begin{equation*}
        \gcd \left(\{m\}\cup \{r:a_r\neq 0\}\right)=1.
    \end{equation*}
\end{theorem}
\begin{proof}
    Let $a$ be the previous greatest common divisor. Notice first that, if $a>1$, then $\lambda=e^{2\pi i /a}\neq 1$ is a complex number such that $\lambda^m=1$ and $\lambda^r=1$ for every $r$ satisfying $a_r\neq 0$. Thus, for every $t\in \C$, it follows that
    \begin{equation*}
        \gamma(\lambda t)=\left((\lambda t)^m, \sum_{r=1}^{+\infty} a_r(\lambda t)^r \right)=\left(t^m, \sum_{r=1}^{+\infty} a_r t^r\right)=\gamma (t).
    \end{equation*}
    Therefore, $\gamma$ cannot be injective as a germ. \\

    Conversely, assume that $a=1$. Let $s\neq t$ be complex numbers, and let us show that $\gamma(s)\neq \gamma(t)$ if $s,t$ are small enough. If $s^m\neq t^m$, then the result follows. Thus, assume that $s^m=t^m$. Then, there exists $\lambda\in \C$ such that $\lambda^m=1$ and $s=\lambda t$. Now, in this case, $\gamma(s)=\gamma(t)$ if and only if $\phi(s)=\phi(t)$, which is in turn equivalent to have $\sum_{r=1}^{+\infty}a_r(\lambda^r-1)t^r\neq 0.$ \\
    Assume by contradiction that the latter power series is zero, that is, for every $r$ such that $a_r\neq 0$, one has that $\lambda^r=1$. Then, if we write $q\in \N$ for the order of $\lambda$, it follows that $q$ divides both $m$ and every $r$ for which $a_r\neq 0$. Therefore, $q$ divides the greatest common divisor of the previous numbers. Hence, in particular, $q\leq a$. Since $s\neq t$, it follows that $q>1$. But in this case, $a>1$, which contradicts our hypothesis. Therefore, the power series $\sum_{r=1}^{+\infty}a_r(\lambda^r-1)t^r\in \C\{t\}$ is nonzero. Since the roots of a nonzero holomorphic function cannot have accumulation, it follows that there exists a radius for which the only root of the previous series is possibly $t=0$. However, if $t=0$, then $s=\lambda t=0$, which contradicts that $s\neq t$. It follows that $\gamma(s)\neq \gamma(t)$ for $s,t$ in a small enough neighbourhood of $0$. 
 \end{proof}
In conclusion, if we have a holomorphic function $\gamma(t)=(t^m, \phi(t))$ that is not injective, one can always compute the previous greatest common divisor $a>1$. Then, define $u=t^a$ as a new parameter and express $\gamma$ in terms of $u$ as $\gamma(u)$. It then holds that the gcd in this case is 1, so $\gamma(u)$ is injective. Hence, every $\gamma$ can be easily modified to be injective, without changing the set-germ of image defined by the original equation $\gamma$.
\section{Branches, multiplicities and tangents}
Along this section, some concepts related to curves are presented. The first of them is the notion of branch, which consists of a curve defined by an irreducible equation. It turns out that a decomposition of a general function in irreducible terms yields that a curve can be decomposed in unions of branches. Moreover, it is shown that a curve is a branch if and only if it can be completely described as the image of a single parametric equation, and that this parametric equation is unique up to a change of parameters. Lastly, the concepts of multiplicity and tangent lines for a curve are introduced. \\

Recall that a curve $(C,0)$ is the set-germ of points that satisfy $f(x,y)=0$ for a nonzero holomorphic function germ $f:(\C^2,0)\rightarrow (\C,0)$. Now, since the ring of holomorphic function germs is a unique factorisation domain, it is possible to write 
\begin{equation*}
    f=\prod_{j=1}^m g_j\,^{a_j},
\end{equation*}
with positive integers $a_j$ and distinct irreducible functions $g_j$ (in the sense that, if we could write $g_k=Ug_j$ for a certain unit $U$, then $j=k$ and $U=1$)  in an essentially unique way. The function $f(x,y)$ will be said to be in \textit{reduced form} if $a_j=1$ for any $j\in \{1,\dots,m\}$. Notice that, if $(B_j,0)$ is the curve defined by the irreducible equation $g_j(x,y)=0$, then 
\begin{equation*}
    (C,0)=\bigcup_{j=1}^m \,(B_j,0).
\end{equation*}
This motivates the following definition:
\begin{definition}
    Let $(B,0)$ be a curve defined by $g(x,y)=0$. We will say that $(B,0)$ is a \textit{branch} if the function $g$ is irreducible in the ring of holomorphic function germs. In addition, $(B,0)$ is said to be a \textit{branch of a curve} $(C,0)$ given by $f(x,y)=0$ if $g$ is an irreducible factor of $f$. In this case, $(B,0)\subset (C,0)$ as set-germs.
\end{definition}
\noindent The above definition does not depend on the choice of $f$ in virtue of the following theorem:
\begin{proposition}\label{ruckert} Let $f,g:(\C^2,0)\rightarrow (\C,0)$ be holomorphic functions. Then, the following are equivalent:
\begin{enumerate}[(i)]
    \item They define the same curve $(C,0)$, that is, the germs at $(0,0)$ of $f^{-1}(0)$ and $g^{-1}(0)$ are equal.
    \item $f$ and $g$ have the same irreducible factors possibly occurring with different multiplicities, up to multiplication by units. 
\end{enumerate}
In particular, if $f$ and $g$ are in reduced form, then $f=Ag$ for some unit $A$. \\
\end{proposition}
The proof of this result requires to introduce a few concepts of general analytic geometry. First, given an ideal $I\subset \C\{x,y\}$, we define the set-germ $(V(I),0)$ as the germ of the set of points in $\C^2$ for which all the functions of $I$ vanish. In addition, for a given set-germ $(X,0)$ of $(\C^2,0)$, define $\mathcal{I}(X,0)\subset \C\{x,y\}$ as the functions that vanish at some neighbourhood of $(0,0)$ in $X$. Lastly, for a given ideal $I\subset \C\{x,y\}$ define the \textit{radical ideal} of $I$, and denote it as $\sqrt{I}$, as the set of functions $f$ for which there exists a natural number $n\in\N$ such that $f^n\in I$. Then, the following holds: 
\begin{theorem}[Rückert-Nullstellensatz] Let $I\subset \C\{x,y\}$ be an ideal. Then, $\mathcal{I}(V(I),0)=\sqrt{I}$.
\end{theorem}
\noindent This result is far from being trivial. A proof can be found in Theorem 3.4.4 of \cite{AnalyticGeometry}.\\

\begin{proof}[Proof of Proposition \ref{ruckert}] The functions $f$ and $g$ define the same germ if and only if $\mathcal{I}(V(f))=\mathcal{I}(V(g))$. By Rückert-Nullstellensatz Theorem this holds if and only if $\sqrt{(f)}=\sqrt{(g)}$. We claim that, if $f=f_1^{n_1}\ldots f_m^{n_m}$ is the decomposition in irreducible functions of $f$, then $\sqrt{(f)}=(f_1\cdot\ldots \cdot f_m)$. With this, if $g=g_1^{k_1}\cdots g_r^{k_r}$ is the irreducible decomposition of $g$, then $(f_1\cdot\ldots \cdot f_m)=(g_1\cdot\ldots \cdot g_r)$ and thus $f_1\cdot\ldots \cdot f_m=Ag_1\cdot\ldots \cdot g_r$ for a unit $A$. By the irreducibility of the $f_i$ and the $g_j$, it follows that $m=r$ and $\{f_1,\ldots, f_m\}=\{A_1g_1,\ldots, A_rg_r\}$ for certain units $A_1,\ldots, A_r$. Thus, the result follows. \\

It only remains to check that $\sqrt{(f)}=(f_1\cdot\ldots \cdot f_m)$. For the direct inclusion, take $h\in \sqrt{f}$. Then, $h^n$ is divisible by $f$ for some natural number $n$. Since $f_i$ divides $f$ for $i\in \{1,\ldots, m\}$, then $f_i$ divides $h^n$, and then it divides $h$ too. Therefore, $f_1\cdot \ldots \cdot f_m$ divides $h$. For the other inclusion, it suffices to check that there exists a natural number $n$ such that $f_1^n\cdot \ldots \cdot f_m^n\in (f)$, but taking $n$ as the maximum of the $n_i$, we would have 
\begin{equation*}
    f_1^n\cdot \ldots \cdot f_m^n=f f_1^{n-n_1}\cdot \ldots \cdot f_m^{n-n_m}\in (f),
\end{equation*}
proving the desired result. The case in which $f$ and $g$ are in reduced form follows since $f=f_1\cdot \ldots \cdot f_m$ and $g=g_1\cdot\ldots \cdot g_r$, and as we saw before, $f=Ag$ for the unit $A$. 
\end{proof}
Therefore, for a given curve $(C,0)$, the irreducible factors $g_j$ do not depend on the choice of $f$ up to multiplication by units. Thus, any curve can be decomposed in different branches, and these branches do not depend on the function $f$ that defines the curve $(C,0)$. In addition, the different branches will only intersect at the origin. Moreover, it is possible to show with these ideas that, if $(C,0)$ is a curve contained in $(D,0)$, then the equation of $(C,0)$ has as irreducible factors some of the ones in the equation for $(D,0)$. 
\begin{proposition}\label{ruckert2}
    Let $f,g:(\C^2,0)\rightarrow (\C,0)$ be holomorphic functions. Then, the following are equivalent:
    \begin{enumerate}[(i)]
        \item The curves $(C,0)$ and $(D,0)$ defined by $f$ and $g$ respectively satisfy $(C,0)\subset (D,0)$. 
        \item The set of irreducible factors of $f$ is a subset of the irreducible factors of $g$ up to multiplication by units.
    \end{enumerate}
    In particular, if $f$ and $g$ are in reduced form, then $f$ divides $g$.
\end{proposition}
\begin{proof}
    Notice that with the notation introduced previously, the first statement can be written as $(V(f),0)\subset (V(g),0)$. Now, it is straightforward to verify that $I(V(g),0)\subset (V(f),0)$ (every function vanishing at $(V(g),0)$ will vanish at $(V(f),0)$, so $I(X,0)$ reverses inclusions). Rückert-Nullstellensatz's Theorem yields $\sqrt{(g)}\subset \sqrt{(f)}$. In the proof of the previous theorem, we have shown that, if $f=f_1^{n_1}\ldots f_m^{n_m}$ and 
    $g=g_1^{k_1}\cdots g_r^{k_r}$ are the irreducible factorisations of $f$ and $g$, then $\sqrt{(f)}=(f_1\cdot\ldots \cdot f_m)$ and $\sqrt{(g)}=(g_1\cdot\ldots \cdot g_r)$. Hence, $(g_1\cdot\ldots \cdot g_r)\subset (f_1\cdot\ldots \cdot f_m)$, forcing $f_1\cdot\ldots \cdot f_m$ to divide $g_1\cdot\ldots \cdot g_r$. Thus, $\{f_1, \ldots, f_m\}\subset \{A_1g_1, \ldots, A_rg_r\}$ for some units $A_1,\dots,A_r$, as we wanted to show. Once more, if $f$ and $g$ are in reduced form, then $f=f_1\cdot \ldots \cdot f_m$ divides $g=g_1\cdot \ldots \cdot g_r$ as we have shown before. 
\end{proof}

Now, we are about to prove that the germ of image of a parametric equation is a branch. 
\begin{theorem}\label{pareq_branch}
    Consider parametric equations 
    \begin{equation*}
        x=t^m, \hspace{0.3cm} y=\phi(t),
    \end{equation*}
    with $\phi:(\C,0)\rightarrow (\C,0)$ holomorphic and $m$ a positive integer. Then, the germ of image defined in $(\C^2, 0)$ by the parametric equations is a branch.  
\end{theorem}
\begin{proof} In order to prove the statement, we should construct an irreducible holomorphic function $D:(\C^2,0)\rightarrow (\C,0)$ whose preimage set-germ of $0$ is precisely the image of the parametric equations. Write $y=\phi(t)=\sum_{r=1}^{+\infty} a_rt^r$ as a convergent power series, and reorder the terms by congruence modulo $m$, writing $r=mq+s$ for $0\leq s \leq m-1$ and $q\in \N$ as 
    \begin{equation*}
        y=\sum_{s=0}^{m-1}t^s \sum_{q=0}^{+\infty} a_{mq+s}t^{mq}. 
    \end{equation*}
If $\phi_s(x)=\sum_{q=0}^{+\infty}a_{mq+s}x^q$, one can write $y=\sum_{t=0}^{m-1}t^s\phi_s(x)$. Notice that each $\phi_s$ is given by a convergent power series near $x=0$. Indeed, since the series that defines $y$ is convergent, there exists $R>0$ and $M>0$ such that $|a_rR^r|<M$ for $r\in \N$. Now, choose $\tilde{R}=R^m>0$. This lets us bound $|a_{mq+s}\tilde{R}^q|=|a_{mq+s}R^{mq}|=|a_{mq+s}R^{mq+s}|R^{-s}<MR^{-s}$. Hence, $\phi_s$ is holomorphic at the origin. Now, since $x=t^m$, for $a\in \{0,\dots,m-1\}$ we can write 
    \begin{equation*}
        t^a y=\sum_{s=0}^{m-a-1}t^{a+s} \phi_s(x)+ \sum_{s=m-a}^{m-1}t^{a+s-m}x\phi_s(x),
    \end{equation*}
    which just consists of multiplying by $t^a$ the equation of $y$, and substituting $x=t^m$ in the powers of $t$ with degree greater or equal than $m$. Consider the previous as $m$ equations with unknowns $t^0,t^1,\ldots, t^{m-1}$, and with coefficients in $\C\{x,y\}$. Then, one has a square homogeneous system of equations, and the nontrivial solutions for $t$ given by the parametric equation forces the determinant of the matrix $D(x,y)$ to vanish for $x=t^m$, $y=\phi(t)$.     Notice that $D$ is a monic polynomial in $y$ of degree $m$ (this will be clearer with an example, but it follows since the matrix of the system has $y$'s in the main diagonal, and nowhere else). In addition, since substituting $t$ by $te^{2\pi i k/m}$ yields parametric equations defining the same set-germ of image, it follows that one can factor
    \begin{equation*}
        D(x,y)=\prod_{k=0}^{m-1} \left(y-\phi(e^{2\pi i k/m}x^{1/m})\right)
    \end{equation*}
    in the ring of fractional formal power series $\C\{x^{1/m}, y\}$ (there are no more possible factors since the degree of $D$ in $y$ is $m$, and $D$ is monic). Therefore, it follows that the germ of curve defined by $D(x,y)=0$ is the same as the one defined by $x=t^m, y=\phi(t)$. \\

    It only remains to check that $D$ is irreducible. Assume that $D=fg$ for functions $f,g\in \C\{x,y\}$. Then, $0=f(t^m,\phi(t))g(t^m,\phi(t))$, and thus since $\C\{t\}$ is a domain, at least one of the terms must vanish identically in a neighbourhood of $0$. Suppose that it is the case for $f(t^m, \phi(t))$. Then, $f^{-1}(0)$ contains the germ of image defined by $x=t^m,y=\phi(t)$, which coincides with $D^{-1}(0)$. It follows by Proposition \ref{ruckert2} that $D$ divides $f$. However, $D=fg$ divides $f$, and hence $g$ has to be a unit. This shows that $D$ is irreducible. 
\end{proof}
\begin{remark}
    Notice that, in the proof, we have seen that the number $m$ that appears in the parametric equations in $x=t^m$ is precisely the order of regularity in $y$ of the function $D$ that defines the branch. In addition, an important fact that arises from the previous proof is the following: for a branch defined by the parametric equations $x=t^m, y=\phi(t)$, the irreducible function $D$ has order $\min\{m, \ord \phi\}$. Indeed, $D$ consists of a product of $m$ factors, each of them having order $\min\{1, \ord\phi /m\}$.\\
\end{remark}
The constructive part of the previous proof can be followed for finding an explicit irreducible equation $f(x,y)=0$ for a set-germ given by parametric equations, as the following example shows.
\begin{exercise} (2.6.1 of \cite{Wall}) Let us consider the curve given by parametric equations $x=t^3$ and $y=t^2+t^4$. Then, following the notation introduced in the proof, we can write $y=tx+t^2$, where $\phi_0(x)=\phi_3(x)=0$, $\phi_1(x)=x$ and $\phi_2(x)=1$. Then, multiplying the previous equation by $t^0, t, t^2$ yields 
    \begin{equation*}
        \begin{cases}
            t^0y=tx+t^2 \\ ty=t^2x+x \\ t^2y = x^2+tx 
        \end{cases}
        \Rightarrow 
        \begin{bmatrix}
            y & -x & -1 \\ -x & y & -x \\ -x^2 & -x & y 
        \end{bmatrix}
        \begin{bmatrix}
            t^0 \\ t \\ t^2
        \end{bmatrix}
        = \begin{bmatrix}
            0 \\ 0 \\ 0
        \end{bmatrix}.
    \end{equation*}
    Thus, as the previous lemma shows, the curve is defined by the determinant of the matrix that appears in the homogeneous system, which is $-x^4-3x^2y-x^2+y^3$. \\
\end{exercise}
Moreover, we can check that the concept of ``set-germ of image of a parametric equation'' is exactly the same as the concept of branch. Indeed, we have shown that every parametric equation defines a branch. It just remains to verify that a branch can always be described as the set-germ of image of a single parametric equation.
\begin{proposition}
    Let $(B,0)$ be a branch. Then, $(B,0)$ can be described as the germ of image of a single parametric equation.
\end{proposition}
\begin{proof}
    Let $f$ be an irreducible function defining $(B,0)$. Assume first that the germ $x=0$ is contained in $(B,0)$. We may thus apply Proposition \ref{ruckert2} to infer that $x$ divides $f$. The irreducibility of $f$ yields $f(x,y)=x$, which can be completely parametrised by $x=0, y=t$. Assume now that $x=0$ is not contained in $(B,0)$. In virtue of Theorem \ref{th:par_eq}, there exist parametric equations $x=t^m, y=\phi(t)$ for $(B,0)$. Hence, it is enough to check that the germ of image defined by these equations is exactly $(B,0)$, since \textit{a priori} one can only assure that the germ of image is contained in $(B,0)$. Construct the irreducible function $D$ of Theorem \ref{pareq_branch} associated to $x=t^m, y=\phi(t)$. Then, $D$ defines the same set-germ as the parametric equations does. Since this germ is contained in $(B,0)$, it follows by Proposition \ref{ruckert2} that $D$ divides $f$. Since $f$ is irreducible and $D$ is not a unit (otherwise the germ defined by $D$ would be empty), we must have that $D=Uf$ for some unit $U$. Hence, $f$ and $D$ define the same branch, that is, $(B,0)$ is defined as the germ of image of the parametric equations $x=t^m, y=\phi(t)$. 
\end{proof}
In addition, it can be shown that the parametric equations that define a branch are unique up to a change of parameters.
\begin{lemma}\label{lemma:uniqueparam}
    Let $t$ and $u$ be the parameters of different parametric equations for the same branch $(B,0)$. Then, there exists a germ of biholomorphism $\omega:(\C,0)\rightarrow (\C,0)$ (that is, with $\omega'(0)\neq 0$) such that $u=\omega(t)$.
\end{lemma}
\begin{proof} Let $f$ be an irreducible equation defining $(B,0)$. If $f(0,y)=0$, then $x$ divides $f$ and thus $f(x,y)=x$. This forces $(B,0)$ to be the curve $x=0$, for which the statement is trivial. In other case, $f$ must be regular in $y$ of some order $m$. Thus, there exists a parametric equation of the form $\xi(u)=(u^m, \phi(u))$ in virtue of Theorem \ref{th:par_eq}. Let $\psi(t)=(\psi_1(t), \psi_2(t))$ be any other parametric equation for $(B,0)$, and write $n\geq 1$ as the order of $\psi_1$. Then, by Lemma \ref{lemma:order}, there exists a holomorphic function $\omega(t)$ of order 1 (that is, a germ of biholomorphism) such that $\psi_1(t)=\omega(t)^n$. Now, both $\xi, \psi : (\C,0)\rightarrow (B,0)$ define bijections as germs of mappings. Write $B$ for a representative set-germ of the branch where the parametric equations $\xi, \psi$ are well-defined and bijective in $\xi^{-1}(B)$ and $\psi^{-1}(B)$, respectively. Hence, for every $(x,y)\in B$ such that $(x,y)\neq (0,0)$, there exist $m$ distinct values of $u\in \xi^{-1}(B)$ such that $x=u^n$, and $n$ distinct values of $t\in \psi^{-1}(B)$ for which $x=\psi_1(t)=\omega(t)^n$. Since $\xi^{-1}\circ \psi : \psi^{-1}(B) \rightarrow \xi^{-1}(B)$ is a bijection, it follows that $n=m$ and that $u^m=\omega(t)^m$. Thus, $u=\lambda\omega(t)$ for some $\lambda\in \C$ with $\lambda^n=1$. Then, the result is proved. 
\end{proof}
\begin{remark}
    The previous results show that branches are precisely the germs of images of a single set of parametric equations. In addition, these equations are unique up to a reparametrization. Therefore, a general curve has as many ``different'' parametric equations as irreducible terms has the equation that defines it. \\
\end{remark}
In what follows, the concepts of multiplicity and tangent lines are going to be described. Let $(C,0)$ be a curve defined by a nonzero holomorphic function germ $f:(\C^2,0)\rightarrow (\C,0)$. Write $f$ as a power series of the form $f(x,y)=\sum_{i,j=0}^{+\infty} a_{i,j}x^i y^j$
for some $a_{i,j}\in \C$, and $a_{0,0}=0$. Recall that the order of $f$ is defined as the least sum of indexes $i+j$ for which $a_{i,j}\neq 0$. Notice that, if the equation that defines $(C,0)$ is changed, the orders may not be the same. However, if $f$ and $g$ are in reduced form (that is, the irreducible factors of them appear with exponent $1$), then $f$ and $g$ define the same curve if and only if there exists a unit $U$ such that $g=Uf$. Since the order of a product is the sum of the orders, and the order of a unit is $0$ as $U(0,0)\neq 0$, then the orders of $f$ and $g$ are equal. This motivates the following definition:
\begin{definition}
    Let $(C,0)$ be a curve defined by a function $f(x,y)=0$ in its reduced form. Then, the \textit{multiplicity of }$(C,0)$\textit{ at }$(0,0)$ is defined as the order of $f$, and it is denoted as $m_{0}(C)$. 
\end{definition}
\noindent Since the whole project will be focused on studying the curves in a neighbourhood of $(0,0)$, the last definition can be shortened just by calling it the \textit{multiplicity of }$(C,0)$ and denoting it as $m(C)$. \\

As we stated before, the order of a product is the sum of the orders. Therefore, if a curve $(C,0)$ is decomposed in its branches $(B_1,0),\dots, (B_m,0)$, it follows that $m(C)=\sum_{j=1}^m m(B_j)$.\\
\begin{remark}
    There is an important relationship between the order of a function defining a branch and the order of a given parametrisation. Namely, if $f$ is an irreducible function defining a branch $(B,0)$ and $x=t^m, y=\phi(t)$ are parametric equations for the branch, then by the remark appearing below Theorem \ref{pareq_branch}, one has that $m(B)=\ord f=\min\{m, \ord \phi\}$. Hence, the minimum of the orders of the parametric equations is the multiplicity of the branch. \\
\end{remark}
Now, we aim to introduce the concept of tangent lines of a curve $(C,0)$. Let us write $m=m(C)$ and $f(x,y)=f_m(x,y)+F(x,y)$, where $f_m(x,y)$ consists of all the terms of $f(x,y)$ with order $m$ (that is, the terms $a_{i,j}x^iy^j$ with $i+j=m$) and $F(x,y)$ with the terms of order greater than $m$. Then, $f_m(x,y)$ is a homogeneous polynomial of degree $m$. As a consequence of the Fundamental Theorem of Algebra, one can show that homogeneous polynomials in two variables split in linear terms, as the following lemma states.
\begin{lemma}
    Let $P(x,y)$ be a homogeneous polynomial of degree $m$. If $x$ does not divide $P(x,y)$, then $P(x,y)$ can be uniquely factorised as
    \begin{equation*}
        P(x,y)=A\prod_{i=1}^n (y-\alpha_i x)^{m_i},
    \end{equation*}
    for $A\in \C\backslash\{0\}$, $\alpha_1,\ldots, \alpha_n\in \C$ distinct and $m_1,\dots, m_n$ positive integers such that $m_1+\ldots+m_n=m$.
\end{lemma}
\begin{proof}
    Write $P(x,y)=\sum_{i=0}^m a_ix^i y^{m-i}$ with complex coefficients. Since $x$ does not divide $P(x,y)$, it follows that $a_0\neq 0$. Now, for all $x\neq 0$, the homogeneity of $P$ yields $P(x,y)=x^m P(1,y/x)$. Since $P(1,t)=\sum_{i=0}^m a_it^{m-i}$ is a polynomial with degree $m$ in $\C[t]$, then the Fundamental Theorem of Algebra gives a unique factorisation of $P(1,t)$ as $P(1,t)=A\prod_{i=1}^n (t-\alpha_i)^{m_i}$, with $n\leq m$, $A\in \C\backslash \{0\}$, $\alpha_1,\ldots, \alpha_n\in \C$ the distinct roots of $P(1,t)$ and $m_1,\ldots, m_n$ the positive integers corresponding to their multiplicities satisfying $m=m_1+\ldots+m_n$. Therefore, 
    \begin{equation*}
        P(x,y)=x^mP(1,y/x)=A\prod_{i=1}^n x^{m_i}(y/x-\alpha_i)^{m_i}=A\prod_{i=1}^n (y-\alpha_ix)^{m_i}
    \end{equation*}
    for $x\neq 0$. With this, the result follows extending the previous by continuity for $x=0$. 
\end{proof}
\begin{remark}
    Notice that, if $x$ divides $P(x,y)$, one can write $P(x,y)=x^p Q(x,y)$ for a positive integer $p\leq m$ and a homogeneous polynomial $Q(x,y)$ of degree $m-p$ that is not divisible by $x$. Then, the previous lemma can be applied to $Q$ so as to give a complete factorisation in linear terms of $P$. \\
\end{remark}
With this, if $f_m$ is the homogeneous polynomial obtained by removing from $f(x,y)$ the terms with order greater than $m$, the lines in $\C^2$ given by the factors of $f_m$ are  defined as the \textit{tangent lines} of the curve defined by $f(x,y)=0$. Notice that this definition does not depend on $f$. Indeed, if $f,g$ both define $(C,0)$ and are in reduced form, then $f=Ug$ for a unit $U$. Thus, the relationship between the corresponding homogeneous polynomials is $f_m=U(0,0)g_m$, since those will be the only terms with order $m$. Therefore, $f_m$ and $g_m$ have the same linear terms in their factorisation and hence, the tangent lines do not depend on the function $f$. Furthermore, notice that, if $f$ and $g$ are functions with homogeneous polynomials $f_m$ and $g_n$ respectively, then the product $fg$ has, as homogeneous polynomial, $f_mg_n$. Thus, the tangent lines of the curves $(C,0)$ and $(D,0)$ defined by $f$ and $g$ are precisely the ones of $(C\cup D,0)$. Moreover, we will show later that branches have just one tangent line. 
\begin{example}
    Let us compute the branches and the tangent lines of the curve defined by $f(x,y)=-4x^4-x^2y^2+4x^2y^3+y^5$. First, the tangent lines of the curve will be given by the equation  $-4x^4-x^2y^2=0$. Thus, they are $x=0$, $x=2ix$ and $x=-2ix$. This assures us that there are at least 3 branches of $f$. Inspired by the solutions of the tangent lines, we can try to perform the division algorithm between $f$ and $4x^2+y^2$. This gives the satisfactory result $f(x,y)=(4x^2+y^2)(x^2-y^3)$. Therefore, the branches of $f$ are $x^2-y^3=0$ with tangent line $x=0$, and $y=\pm 2ix$, which are themselves tangent lines. \\
\end{example}
Lastly, we are going to define the concept of singular point of a curve. 
\begin{definition}
    Let $(C,0)$ be a curve defined by $f(x,y)=0$ in its reduced form, where $f:(\C^2, 0)\rightarrow (\C,0)$ is a nonzero holomorphic function germ. The origin is said to be a \textit{singular point} of $(C,0)$ if $\partial f/\partial x\,(0,0)=\partial f/\partial y\,(0,0)=0.$ Otherwise, the origin will be said to be non-singular. The \textit{multiplicity} of the singularity at the origin of the curve $(C,0)$ is defined as the multiplicity of the curve $(C,0)$, that is, as $m(C)$. 
\end{definition}
\noindent Notice that $(0,0)$ is non-singular if and only if $m(C)=1$. In fact, more can be stated:
\begin{proposition}
    The origin is a non-singular point of a curve if and only if it is a branch that can be described as the graph of a holomorphic function, that is, up to a change of coordinates of the form $(x,y)\mapsto (y,x)$, there exists a holomorphic function germ $\phi:(\C,0)\rightarrow (\C,0)$ such that $x=t, y=\phi(t)$ are parametric equations for the curve. 
\end{proposition}
\begin{proof} Let $(C,0)$ be a curve with multiplicity $1$ and take $f(x,y)$ as an equation for the curve in reduced form and order $1$. Then, one of the partial derivatives $\partial f / \partial x$ or $\partial f / \partial y$ does not vanish at the origin. Thus, the result follows by applying the Implicit Function Theorem to $f$. 
\end{proof}
\begin{example} The functions $f(x,y)=(x-y)(x+y)$ and $f(x,y)=y^2-x^3$ both define curves with a singularity of multiplicity 2. In fact, this is the general picture in this case, since if we assume that the function $f(x,y)$ is reduced (that is, it has no repeated factors), then the curves with a singularity with multiplicity 2 can be the union of two branches with non-singular points at the origin, or a single branch with multiplicity 2 (corresponding to an irreducible function of degree 2). \\
\end{example}
Lastly, an application of the Preparation Theorem yields a result concerning the curves that have singularities of order 2. These are known in the literature as the $A_k$ singularities.
\begin{theorem}
    Let $(C,0)$ be a curve with a singularity of order 2. Then, up to a holomorphic change of coordinates, $(C,0)$ is given by $v^2+u^k=0$ for some integer $k\geq 2$. 
\end{theorem}
\begin{proof}
    Let $f$ be a reduced holomorphic function defining $(C,0)$. Since $f$ has order 2, we can assume that $f$ is regular in $y$ of order 2 (otherwise, the only terms of degree 2 in $f(x,y)$ are $x^2$ and $xy$, but a linear change of coordinates would bring a function regular in $y$ of order 2). Then, by Weierstrass Preparation Theorem, $f(x,y)=U(x,y)(y^2+a(x)y+b(x))$ for some unit $U$ and some holomorphic functions $a,b$. Therefore, the curve $(C,0)$ may be given by the equation $y^2+a(x)y+b(x)=0$. A change of coordinates $v=y+a(x)/2$ transforms the equation to $v^2+c(x)=0$. If $c=0$, the curve is given by $v^2=0$, which is not in reduced form. Then, $c$ is nonzero, so it has to be a function of order $k\geq 2$. By Lemma \ref{lemma:order}, there exists a function $\phi(x)$ of order 1 such that $c=\phi^k$ (since it has order 1, it is a germ of biholomorphism). Thus, performing the change $u=\phi(x)$ yields that $(C,0)$ is given by $v^2+u^k=0$.
\end{proof}
\section{Intersection number between two curves}
In order to study how different branches interact between themselves, it will be mandatory to introduce the concept of intersection number. This section is devoted to study this notion, as well as to provide some meaningful consequences for the structure of a branch.
\begin{definition}
    Suppose that we have a branch $(C,0)$ given by the equation $f(x,y)=0$, where $f$ is supposed to be irreducible (that is, we exclude the trivialities of the kind $f^2=0$) and $(D,0)$ a branch given by parametric equations $(\phi(t), \psi(t))$ (that are by definition required to be injective). Then, the order of $f(\phi(t), \psi(t))$ is said to be the \textit{intersection number of }$(C,0)$ \textit{and} $(D,0)$, and it is denoted as $(C.D)_0$, or just as $C.D$. 
\end{definition}
\noindent Notice that, since $f(\phi(0),\psi(0))=0$, the intersection number is a positive integer. In the particular case that $f(\phi(t), \psi(t))$ vanishes identically in a neighbourhood of $0$, the branches must be equal. Therefore, the intersection number can be interpreted to be, in this case, $+\infty$. 
\begin{example}
    Let $f(x,y)=x^2-y^3$ and $g(x,y)=x^2-y^5$, and denote as $(C,0)$ and $(D,0)$ the corresponding curves associated to them, respectively. Since $\alpha(t)=(t^3, t^2)$ are parametric equations for $(C,0)$, we can write $g(\alpha(t))=t^6-t^{10}.$ Thus, $C.D=6$.
\end{example}
\noindent As in the previous example, if the curves are both given by an equation of the form $f(x,y)=0$, it will be necessary to find parametric equations for one of them. \\

The following proposition gives natural and expected properties for the intersection number, as well as a formula for computing it when the branches are both expressed as preimages of irreducible functions. 
\begin{proposition}\label{prop:intnumber}
    Let $(C,0)$ be a branch given by the irreducible equation $f(x,y)=0$ and let $(D,0)$ be a branch defined by the parametric equations $\gamma(t)=(\phi(t),\psi(t))$. Then,
    \begin{enumerate}[(1)]
        \item The intersection number $C.D$ is well-defined, that is, it does not depend on the choice of the function that defines $(C,0)$ or the choice of the parametric equations of $(D,0)$. 
        \item If $g$ is an irreducible equation for $D$, then $$C.D=\dim_{\C}\dfrac{\C\{x,y\}}{(f,g)}.$$ 
        \item The intersection number is symmetric.
    \end{enumerate}
\end{proposition}
\begin{proof} Let us first check (1). If $f,\tilde{f}$ are irreducible expressions for $C$, then they are in reduced form. Thus, as a consequence of Proposition \ref{ruckert}, $\tilde{f}=Uf$ for some unit $U$. Since units have order $0$ and the order of the product is the sum of the orders, it follows that the orders of $f(\phi(t),\psi(t))$ and $\tilde{f}(\phi(t), \psi(t))$ coincide. Now, if we consider another parametric equations for $D$, by Lemma \ref{lemma:uniqueparam} it follows that there exists a germ of biholomorphism $u\in \C\{t\}$ with $u(0)=0$ for which the parametric equations are given by $(\phi(u(t)), \psi(u(t))$. Notice that $u$ has order 1, and therefore the orders of $f(\phi(t), \psi(t))$ and $f(\phi(u(t)), \psi(u(t))$ are equal. \\

Before checking (2), notice that (3) is a trivial consequence of it. In order to check (2), write $k$ as the order of $v(t)=f(\phi(t),\psi(t))$. Notice that $v=t^k U(t)$ for a unit $U$. Thus, we may have the equality $(v)=(t^k)$ as ideals of $\C\{t\}$. Therefore, $k$ is the dimension of the quotient $\C\{t\}/(v)$ seen as a $\C$-vector space, since $\C\{t\}$ is the direct sum of the subspace formed by the polynomials of degree $<k$ and $(t^k)$. In order to prove the desired equality, we are going to apply a statement relating the topological degree of a map with the dimension of a quotient space. Namely, if $F:(\C^n,0)\rightarrow (\C^n,0)$ is a holomorphic function with $(F^{-1}(0),0)=(\{0\},0)$, then
\begin{equation*}
    \deg F=\dim_{\C}\dfrac{\C\{x_1,\ldots, x_n\}}{(F_1,\dots,F_n)}.
\end{equation*}
A proof can be found in \cite{Juanjo}, in Corollary D.7 (Recall that the degree of a holomorphic function germ $F$ is calculated fixing a representative function of the germ and counting the roots of $F(x)=t$ for small enough regular value $t$). With this, it is enough to check that $\deg F=\deg v$, where $F=(f,g):(\C^2,0)\rightarrow (\C^2,0)$. \\

Let $z\in \C\setminus\{0\}$ be small enough. Then, $z$ is a regular value of $v$. Assume that $v^{-1}(z)=\{t_1,\ldots, t_k\}$. We claim that $F^{-1}(z,0)=\{v(t_1),\dots,v(t_k)\}$ and that $(z,0)$ is a regular value of $F$. This forces $\deg F= \deg v$ as we desire to show. For the first claim, notice that the points $(x,y)\in F^{-1}(z,0)$ satisfy $z=f(x,y)$ and $g(x,y)=0$. Since $\gamma$ is the parametric equation of the curve $D$ defined by $g$, one has that $(x,y)=\gamma(t)$. Hence, $z=(f\circ \gamma) (t)=v(t)$. Thus, $t\in v^{-1}(z)=\{t_1,\dots,t_k\}$, forcing $(x,y)\in \{\gamma(t_1),\dots,\gamma(t_k)\}$. This shows the first assertion. In order to prove the second one, we may check that, for any $t\in \{t_1,\dots,t_k\}$, we have that $J_F(\gamma(t))\neq 0$. \\

First, since $t$ is a regular value of $v$, one has that $v'(t)=f_x(\gamma(t))\gamma_1'(t)+f_y(\gamma(t))\gamma_2'(t)\neq 0$. On the other hand, the equality $g(\gamma(t))=0$ yields that $g_x(\gamma(t))\gamma_1'(t)+g_y(\gamma(t))\gamma_2'(t)=0$. Thus, $(g_x,g_y)(\gamma(t))=a(-\gamma_2'(t),\gamma_1'(t))$ for some $a\in\C$. In addition, since $z\neq 0$, then $t\neq 0$ and $\gamma (t)\neq 0$ by injectivity. Thus, $\gamma(t)$ is a regular point of $g$, forcing $(g_x,g_y)(\gamma(t))\neq (0,0)$. Then, $a\neq 0$, and therefore 
\begin{equation*}
        J_F(\gamma(t))=f_x(\gamma(t))g_y(\gamma(t))-f_y(\gamma(t))g_x(\gamma(t))=a\left(f_x(\gamma(t))\gamma_1'(t)+f_y(\gamma(t))\gamma_2'(t)\right)=av'(t)\neq 0,
\end{equation*}
showing that $z$ is a regular value of $F$. 
\end{proof}
In what follows, our intention is to establish a relationship between the intersection number defined as above and the one defined topologically:
\begin{definition}
    Let $\alpha, \beta:(\C,0)\rightarrow (\C^2,0)$ be parametric equations for the branches $(C,0)$ and $(D,0)$, respectively. If we consider them as maps $\alpha,\beta:(\R^2, 0)\rightarrow (\R^4, 0)$, then its \textit{topological intersection number} is defined as 
    \begin{equation*}
        I(\alpha, \beta)=\sum_{\alpha(x)=\tilde{\beta(x)}}I_x(\alpha, \tilde{\beta}),
    \end{equation*}
    where $\tilde{\beta}$ is homotopic to $\beta$ and transversal to $\alpha$, and $I_x(\alpha, \tilde{\beta})$ is $\pm 1$ depending on whether the space $\text{Im}\,d\alpha_x \oplus \text{Im}\,d\beta_x \subset \R^4$ has the natural orientation induced by $\R^4$ or not. 
\end{definition}
\noindent It then turns out that both intersection numbers coincide.
\begin{proposition}
    If $(C,0)$ and $(D,0)$ are branches with parametric equations given by $\alpha$ and $\beta$, then $C.D=I(\alpha, \beta)$. 
\end{proposition}
\begin{proof} Since $\alpha$ and $\beta$ are holomorphic maps, its local intersection numbers in every point are $+1$, as holomorphic maps preserve the orientation. Now, if $f:(\C^2,0)\rightarrow (\C,0)$ is an irreducible holomorphic germ of mapping that defines $(C,0)$, then the positive integer $m$ of the expansion $f(\beta(t))=\sum_{n=m}^{+\infty} a_nt^n$ with $a_m\neq 0$ is the intersection number $C.D$. We claim that the number of solutions of $\sum_{n=m}^{+\infty} a_nt^n=0$ is equal to $m$ counting multiplicities. Since the local intersection number in each of them is $+1$, it follows that $m=I(\alpha, \beta)$. \\

The claim regarding the number of solutions for the given equation is a consequence of the fact that the order of a function $f(t)$ is its degree. Indeed, the order of $f$ is the dimension as a $\C$-vector space of the quotient $\C\{t\}/(f)$, and this coincides with the degree, as an application of the result mentioned in the proof of Proposition \ref{prop:intnumber} (2).
\end{proof}
Lastly, as an application of the intersection number, it can be checked whether or not a line is tangent to a given curve just by looking at the intersection number of the curve and the line.
\begin{proposition}
    Let $(B,0)$ be a branch and $L$ a line in $\C^2$. Then, $L$ is a tangent line of $B$ if and only if $L.B>m(B)$. 
\end{proposition}
\begin{proof} Let $f$ be an irreducible function defining the branch $(B,0)$, and let $m$ be the order of $f$. Assume first that the line has as equation $y=ax$ for $a\in \C$. If $y-ax$ is not a factor of the homogeneous polynomial $f_m$ obtained by removing from $f$ the terms of order greater than $m$, then the substitution $f(x,ax)$ gives a function with order $m$. If $y-ax$ is a factor of $f_m$, then $f(x,ax)$ will have order greater than $m$. Thus, the multiplicity $m$ of the curve $(B,0)$ is in general smaller than or equal to $L.B$, with equality if and only if $L$ is a tangent line. Lastly, the change of coordinates $(x,y)\mapsto (y,x)$ is valid for showing our claim if $L$ is defined by $x=ay$. Thus, both cases cover all the possibilities for $L$. 
\end{proof}
\begin{theorem}
    A branch has a unique tangent line. 
\end{theorem}
\begin{proof}
    Let $(B,0)$ be a branch and assume that the germ of $x=0$ at the origin and $(B,0)$ are not equal (otherwise, the statement would be trivial, since $x=0$ would be the only tangent line in virtue of the previous proposition). Let $x=t^m, y=\phi(t)=\sum_{r=1}^{+\infty} a_rt^r$ be parametric equations for $(B,0)$. If $x=0$ is not tangent to $(B,0)$, then substituting in $x=0$ gives a function of order $m(B)$, and hence $m=m(B)$. In this case, if $L$ is a tangent line of $(B,0)$ given by $y=ax$, then the order of $at^m-\sum_{r=1}a_rt^r$ is greater than $m(B)=m$, since it is the intersection number $L.B$. Therefore, $a_r=0$ if $r<m$ and $a=a_m$. This forces $y=a_mx$ to be the only tangent line of $(B,0)$. If $x=0$ is tangent to $(B,0)$, then substitution in $x=0$ gives a function of order greater than $m(B)$. Thus, $m>m(B)$. Hence, in this case, $m(B)=\min\{m, \ord\phi\}=\ord \phi$. Therefore, if $y=ax$ were tangent to $(B,0)$, in particular $\phi(t)-at^m$ would have order greater than $m(B)$, which is impossible whenever $\ord\phi=m(B)$. Thus, $x=0$ is the unique tangent line in this case. 
\end{proof}
\noindent In particular, in the previous proof we have seen that every branch that is not tangent to $x=0$ can be parametrised by $x=t^m, y=\sum_{r=m}^{+\infty}a_rt^r$, and the term $a_m$ gives the slope of the tangent line $y=a_mx$ of the branch. \\

The following theorem grants us that, outside the origin, the curves are one-dimensional complex manifolds, and that the tangent lines in the origin are limits of the tangent lines associated to the manifold structure.
\begin{theorem}\label{theorem:submanifold} Let $(B,0)$ be a branch. Then, $(B\setminus\{0\},0)$ is a germ of regular submanifold of $(\C^2,0)$ with complex dimension 1. Moreover, the tangent line at a point $z\neq 0$ converges to the tangent line of $(B,0)$ at $0$ as $z\to 0$, and the same holds for the line defined by the position vector $\gamma(z)$.
\end{theorem}
\begin{proof}
    Notice first that it is possible to assume that $(B,0)$ is not tangent to $x=0$. Otherwise, apply the change of coordinates $(x,y)\mapsto (y,x)$.
    Let $\gamma:(\C,0)\rightarrow (\C^2,0)$ be parametric equations for $(B,0)$ of the form $\gamma(t)=(t^m, \sum_{n\geq m} a_nt^n)$. Then, the restriction $\gamma:(\C\setminus\{0\},0)\rightarrow (\C^2,0)$ is clearly the germ of a holomorphic embedding, since it is injective, open onto its image and $\gamma'(t)=(mt^{m-1}, \sum_{n\geq m} na_n t^{n-1})\neq (0,0)$ for every $t\neq 0$. Hence, $(B\setminus\{0\},0)$ is a regular one-dimensional complex submanifold of $(\C^2,0)$, and its tangent space in a point $\gamma(t)$ is the complex vector space $T_{\gamma(t)}(B,0)=\langle \gamma'(t)\rangle$. Lastly, the convergence claim can be studied in the projective space $\CP^1$. Indeed, the vector $\gamma'(t)\in \C^2$ in the projective space has homogeneous coordinates $(mt^{m-1}: ma_mt^{m-1}+\ldots\,)=(1:a_m+\ldots\,)$. Therefore, as $t\to 0$, these tangent lines tend to $(1:a_m)$, which is the tangent line of $(B,0)$ at the origin. In addition, since the class of $\gamma(t)$ in $\CP^1$ is $(t^m, a_mt^m+\ldots)=(1:a_m+\ldots)$, then, as $t\to 0$, the complex line defined by the position vector also converges to the tangent line of $B$ at $0$. 
\end{proof}
\chapter{Arithmetic properties of branches}
This chapter is intended to introduce some arithmetic tools that will play a fundamental role in the understanding of curves and branches. The main interest of the concepts that are studied along this chapter underlies the fact that the topological nature of a curve is completely determined by them. The most important notion, presented in the first section, is the Puiseux characteristic of a branch. Since it will be shown to be an \textit{invariant} under the topological equivalence of branches, and its nature is completely arithmetic, it will be common to say that the Puiseux characteristic is an \textit{arithmetic} or a  \textit{numerical invariant}. Along the second section, the contact of branches will be introduced as an equivalent method to measure the intersection number between branches. Lastly, the semigroup of a branch will be analysed in the third section, together with some key properties that will become crucial in the last chapter of the project. 
\section{Puiseux characteristic}
In this section, a numerical invariant for branches is going to be presented. Namely, let $(B,0)$ be a branch that is not tangent to $x=0$ at the origin. The results of the previous chapter show that $(B,0)$ admits parametric equations of the form $x=t^m$, $y=\phi(t)=\sum_{r=m}^{+\infty}a_rt^r$, where the possibility of having $a_m=0$ is not excluded. Since parametric equations are injective, then $\gcd \left(\{m\}\cup \{r:a_r\neq 0\}\right)=1$ by Theorem \ref{thrm:alwaysinjective}. The terms of the Puiseux characteristic will be defined, loosely speaking, as the most relevant numbers appearing in $\{r:a_r\neq 0\}$ that are contributing to the fact that $\gcd \left(\{m\}\cup \{r:a_r\neq 0\}\right)=1$. \\

Let us formalise the previous idea. If $m\neq 1$, define $\beta_1=\min\{k\in \N: a_k\neq 0, m\nmid k\}$, \textit{i.e.}, as the first exponent in $\phi(t)$ that is not divisible by $m$, and let $e_1=\gcd\{m,\beta_1\}$. Notice that $e_1$ is a positive integer strictly smaller than $m$, and $\beta_1$ is well-defined by the injectivity clause. Inductively, if $\beta_1, \ldots, \beta_{i}$ and $e_1,\ldots, e_i$ have been defined and $e_i>1$, define  
\begin{equation*}
    \beta_{i+1}=\min \{k\in \N: a_k\neq 0, e_i\nmid k\}\text{ and } e_{i+1}=\gcd\{e_i, \beta_{i+1}\}.
\end{equation*}
Then, $\beta_{i+1}$ is the least exponent of $\phi(t)$ that is not divisible by $e_i$. Again, since $e_i\nmid \beta_{i+1}$, it follows that $e_{i+1}<e_i$. Repeat the process until $e_g=1$ for some $g\in \N$ . Once more, the injectivity clause for the parametric equations grants the existence of such a $g$. The \textit{Puiseux characteristic} of $(B,0)$ is then defined as the tuple $(m;\beta_1, \ldots, \beta_g)$, where the $\beta_i$ are known as the Puiseux terms of $(B,0)$ and the $e_i$ as the auxiliary terms. It will be convenient to define $\beta_0$ and $e_0$ as $m$. \\

Since the definition of the $\beta_i$ and $e_i$ does not depend on the value of the coefficients $a_r$ of the parametric equations as long as they are nonzero, it follows that a sufficient condition for the branches $x=t^m, y=\sum_{r=m}^{+\infty}a_rt^r$ and $x=t^m, y=\sum_{r=m}^{+\infty}b_rt^r$ to have the same Puiseux characteristic is that the exponents that appear in both equations are the same, \textit{i.e.}, $\{r:a_r\neq 0\}=\{r:b_r\neq 0\}$. As a straightforward consequence, the Puiseux characteristic of a branch does not depend on the chosen parametric equations. Indeed, if $x=t^m, y=\phi(t)$ are parametric equations for a given branch, then any other parametric equation with $x=t^m$ will have the form $x=t^m , y=\phi(\lambda t)$ for some $\lambda\in \C$ such that $\lambda^m=1$. Therefore, both equations yield the same Puiseux characteristic.
\begin{example}
    Let $(B,0)$ be the branch with equations $x=t^6$, $y=t^6+9t^{12}+2t^{27}-4t^{81}+t^{83}$. Then, $\beta_1=\min \{k\in \N: a_k\neq 0, 6\nmid k\}=27$ and $e_1=\gcd\{6,27\}=3$. With this, $\beta_2=\min \{k: a_k\neq 0, 3\nmid k\}=83$ and $e_2=\gcd\{83,3\}=1$. Thus, the Puiseux characteristic of $(B,0)$ is $(6;27,83)$. Hence, $(B,0)$ has the same Puiseux characteristic as the branch defined by $x=t^6, y=t^{27}+t^{83}$. \\
\end{example}
Notice that the Puiseux characteristic is a finite tuple of increasing numbers. In addition, its length can be arbitrarily long. For instance, let $n\in \N\setminus \{0\}$ and $$x=t^{2^n}, y=t^{2^n\cdot 3^0}+t^{2^{n-1}\cdot 3^1}+\ldots + t^{2^1\cdot 3^{n-1}}+t^{3^n}.$$
Then, its Puiseux characteristic, $(2^n; 2^{n-1}\cdot 3^1, \ldots, 2^1\cdot 3^{n-1}, 3^n)$, has length $n+1$.\\

The Puiseux characteristic of a branch which is tangent to $x=0$ is not well-defined at a first stage. The following proposition provides the intuition to do it, but its proof will be postponed until the exponent of contact is introduced (see Corollary \ref{prop:Puisuexinvariant}).
\begin{proposition}
    The Puiseux characteristic is invariant under holomorphic changes of coordinates. Formally, if $(B,0)$ and $(\tilde{B},0)$ are branches and there exists a biholomorphism $h:(\C^2,0)\rightarrow (\C^2,0)$ (\textit{i.e.} a holomorphic and bijective map with a holomorphic inverse) with $h(B,0)=(\tilde{B},0)$, then $(C,0)$ and $(\tilde{B},0)$ have the same Puiseux characteristic.
\end{proposition}
\noindent Therefore, if $(B,0)$ is a branch that is tangent to $x=0$, it is natural do define its Puiseux characteristic as the one of its image $h(B,0)$ under a biholomorphism $h$ that sends the line $x=0$ to a branch that is not tangent to $x=0$. Moreover, this characteristic is independent of the choice of $h$. In practice, the chosen change of coordinates will be the symmetry $(x,y)\mapsto (y,x)$. \\
\begin{definition}
    The branches $(B,0)$ and $(\tilde{B},0)$ are said to be \textit{equisingular} whenever they have the same Puiseux characteristic. For general curves, $(C,0)$ and $(\tilde{C},0)$ are said to be \textit{equisingular} whenever they have the same number $n$ of branches, namely $\{(B_1,0), \ldots, (B_n,0)\}$ for $(C,0)$ and $\{(\tilde{B}_1,0), \ldots, (\tilde{B}_n,0)\}$ for $(\tilde{C},0)$, and there exists a permutation $\sigma \in \Sigma_n$ such that $(B_{\sigma(i)},0)$ and $(\tilde{B}_i,0)$ are equisingular as branches, and the intersection numbers $B_{\sigma(i)}.B_{\sigma(j)}$ and $\tilde{B}_i.\tilde{B}_j$ coincide. In other words, two curves are equisingular if there exists a bijection between their branches that preserves the Puiseux characteristic and the intersection numbers. 
\end{definition}
\noindent Equisingularity gives a natural method for performing a classification of plane curves. In fact, this equivalence relation turns out to be extremely practical, since the calculation of the Puiseux characteristic of a branch given by its parametric equations is routine. However, if a branch is defined by an irreducible function $f(x,y)=0$, it is not clear how its Puiseux characteristic could be obtained. In fact, although Newton’s algorithm provides a method for calculating the parametric equations for a curve, it is not granted that the algorithm ends in a finite number of steps. Nevertheless, Newton’s algorithm yields approximations for the parametric equations, and they turn out to be enough to obtain the Puiseux characteristic of the branch, as the following example shows. 
\begin{example}
    Let $(B,0)$ be the branch defined by $f(x,y)=x^2-y^3-5x^2y$. Notice first that the homogeneous part of $f$ is given by $x$. Hence, the tangent line of the branch $(B,0)$ is $x=0$. Since our intention is to obtain the Puiseux characteristic, let us perform the holomorphic change of coordinates $(x,y)\mapsto (y,x)$. Then, $g(x,y)=f(y,x)=y^2-x^3-5xy^2$ defines a branch with the same Puiseux characteristic as $(B,0)$. An application of Newton's algorithm provides a first approximation given by the dominating terms $y^2-x^3=0$. Since $m(B)=\ord f=2$, then Newton's algorithm grants us that the parametric equation of the symmetric curve is of the form $x=t^2, y=t^3+\ldots$. Therefore, the Puiseux characteristic of $B$ is $(2;3)$ since $\gcd (2,3)=1$ and hence the terms of higher order of $y$ are not relevant for the characteristic. \\
\end{example}
In general, for a branch given by an irreducible function $f(x,y)$ not tangent to $x=0$, Newton's algorithm offers in a finite number of steps an approximation of the form $x=t^m, y=\phi(t)+\ldots $, where $m=\ord f$. Now, keep on iterating the algorithm until the terms of $\phi(t)$ are enough to obtain the Puiseux characteristic (that is, when $e_g=1$ for some $g\in \N$). \\

As it has been commented on before, equisingularity is preserved under holomorphic changes of coordinates. One could speculate on whether the converse of this statement holds, and an answer is given in the following counterexample.
\begin{exercise} (4.7.9 of \cite{Wall}) Show that the curves $y^3+x^6=0$ and $y^3+yx^4=0$ are equisingular, but cannot be reduced to each other by a holomorphic change of coordinates.
\end{exercise}
\noindent {\bf Solution} Let us check first that the given curves are indeed equisingular. In order to do so, it will be necessary to perform the decomposition of the curves in their branches. After some algebraic calculations, one obtains that
\begin{align*}
    &y^3+x^6=(y+x^2)(y^2-yx^2+x^4)=(y+x^2)(y-\omega x^2)(y-\bar{\omega}x^2),\\
    &y^3+yx^4=y(y^2+x^4)=y(y-ix^2)(y+ix^2),
\end{align*}
where $\omega=e^{2i\pi /3}$. This provides each curve its three corresponding branches. Since all six of them have parametric equations of the form $(t,\phi(t))$, it follows that the multiplicities of the branches are $m=1$, and hence all of them have Puiseux characteristic $(1)$. \\

It only remains to compute the intersection numbers between the three branches for both curves. For instance, the intersection number between $y-x^2=0$ and $y-\omega x^2=0$ is 2, since $(t,\omega t^2)$ are parametric equations for the second one, and a substitution in the first one yields the function $(\omega-1)t^2$ of order 2. A similar argument shows that any other nontrivial intersection number is $2$, as table 1 encapsulates.

\begin{table}[ht]
    \centering
    \begin{tabular}{|c|c|c|c|}
    \hline
                        & $y+x^2$  & $y-\omega x^2$ & $y-\bar{\omega}x^2$ \\ \hline
    $y+x^2$             & $\infty$ & $2$            & $2$                 \\ \hline
    $y-\omega x^2$      & $2$      & $\infty$       & $2$                 \\ \hline
    $y-\bar{\omega}x^2$ & $2$      & $2$            & $\infty$            \\ \hline
    \end{tabular}
    \caption{Intersection number between the branches of $y^3+x^6$.}
\end{table}
\noindent Since the same scenario holds for the branches of $y^3+yx^4=0$, it turns out that any permutation of the branches respects equisingularity for the branches and the intersection numbers. Hence, the given curves are equisingular. \\

Let us see that the curves $y^3+x^6=0$ and $y^3+yx^4=0$ cannot be reduced to each other by a holomorphic change of coordinates by contradiction. Assume that there exists coordinates $(x,y)\mapsto (\tilde{x}, \tilde{y})$ such that $\tilde{y}^3+\tilde{x}^6=y^3+yx^4$. Since the map has to send the origin to the origin, one can write 
\begin{equation*}
    \left\{ \begin{array}{l}
        \tilde{x}=\phi(x,y)=a_{10}x+a_{01}y+a_{20}x^2+a_{11}xy+a_{02}y^2+\ldots \\
        \tilde{y}=\psi(x,y)=b_{10}x+b_{01}y+b_{20}x^2+b_{11}xy+b_{02}y^2+\ldots \\
        \end{array}
\right.
\end{equation*}
Comparing the expressions, it follows that $y^3+yx^4$ has the same coefficients with terms of degrees 3, 4 and 5 that appear in the expansion 
\begin{equation*}
    (b_{10}x+b_{01}y+b_{20}x^2+b_{11}xy+b_{02}y^2+b_{30}x^3+b_{21}x^2y+b_{12}xy^2+b_{03}y^3)^3.
\end{equation*}
First, since the coefficients in the last expansion of $x^3$ and $y^3$ are, respectively, $b_{10}^3$ and $b_{01}^3$, then $b_{10}=0$ and $b_{01}=1$. Taking this into account, it follows that the coefficient of $x^2y^2$ is given by $3b_{20}$. Thus, $b_{20}=0$. However, in this case the term $x^4y$ cannot appear in the previous expansion. Therefore, these curves are not equivalent under a holomorphic change of coordinates. \\

The previous example shows that the classification of plane curves through holomorphic changes of coordinates is in fact a much finer question than the one posed by the equisingularity relation. The difficulty that lies in the analytic classification forces it to be out of the reach of the dissertation. This classification was incomplete until March 2021, when Marcelo Escudeiro Hernandes and Maria Elenice Rodrigues Hernandes achieved the complete analytic classification of plane curves in \cite{analyticclassification}. \\

When it comes to the classification of curves through the topological equivalence (which is, loosely speaking, a classification via changes of coordinates that are only required to be homeomorphisms), the situation is completely different. Not only is the Puiseux characteristic an invariant under the topological equivalence of plane branches, but it also turns out to be a complete invariant for the classification. The proof of these assertions are the object of study of the third chapter. 
\section{Exponent of contact}
This section provides a way to measure how branches interact between themselves. This new tool, called the exponent of contact between branches,  is deeply related with the intersection number. As we know, if $(B,0)$ is a branch that is not exactly the germ of vertical line $x=0$, the results of the first chapter of the project show that it can be expressed through equations of the form $x=t^m, y=\sum_{r=1}^{+\infty} a_r t^r.$ Our first goal is to to provide a precise meaning for the expression of $y$ in terms of $x$ given by $y=\sum_{r=1}^{+\infty} a_rx^{r/m},$ further from considering it just as a formal series in the ring $\C\{x^{1/m}\}$. In order to do so, it would be mandatory to define properly what $x^{1/m}$ means, in terms of $m$-th roots. As we know, every nonzero complex number has exactly $m$ different $m$-th roots in $\C$. Hence, the correspondence $\C\rightarrow \C$ that assigns $m$-th roots to every complex number is clearly not a function. \\

Let us restrict our attention to the domain $\S=\C\backslash (-\infty,0]=\{re^{i\theta}: r>0, -\pi <\theta<\pi\}$, and define the \textit{canonical }$m$\textit{-th root function in }$\S$, denoted as $x^{1/m}$, as the function $\S\rightarrow \C$ given by $re^{i\theta}\mapsto r^{1/m} e^{i\theta/m}$, where $r>0$ and $\theta \in (-\pi,\pi)$. This clearly gives a well-defined and continuous function in $\S$, and satisfies that $(x^{1/m})^m=x$ for every $x\in \S$. Furthermore, the other $m$-th roots for a fixed complex number are obtained through multiplication by an $m$-th root of the unity. Hence, any other continuous function $f:S\rightarrow \C$ that satisfies $f(x)^m=x$ is of the form $re^{i\theta} \mapsto \omega r^{1/m} e^{i\theta/m}$, for $r>0$, $\theta\in (-\pi,\pi)$ and $\omega\in \C$ with $\omega^m=1$. Thus, there are exactly $m$ distinct $m$-th root functions defined in $\S$, and are of the form $\omega x^{1/m}$ for $\omega\in \C$ an $m$-th root of the unity. \\

Returning to the branch $(B,0)$ with equations $x=t^m, y=\sum_{r=1}^{+\infty} a_r t^r$, a \textit{pro-branch} of $(B,0)$ in $\S$, or just a \textit{pro-branch} of $(B,0)$, is a function of $y$ in terms of $x$ extracted from the parametric equations, where a choice of an $m$-th root in $\S$ has been made. The set of such functions is denoted as $\pro(B)$. In other terms, $\gamma\in \pro(B)$ means that $\gamma(x)=\sum_{r=1}^{+\infty} a_r \omega^r x^{r/m}$, for some $\omega\in \C$ with $\omega^m=1$, where $x^{r/m}$ is just $(x^{1/m})^r$. Notice that pro-branches can be similarly defined over any sector of $\C$, just by performing a choice of an $m$-th root function along the sector. However, it will be enough for our purposes to restrict our attention to $\S$. \\

Recall that every branch that is not tangent to $x=0$ admits a parametric equation of the form $x=t^m, y=\sum_{r=m}^{+\infty}a_rt^r$, where $m$ is precisely the multiplicity of the curve. From now on, and unless otherwise is stated, it would be assumed that every branch is not tangent to $x=0$. \\

Let $(B,0)$ and $(B',0)$ be branches. Then, for every $\gamma\in \pro (B)$ and $\gamma'\in \pro (B')$, one may write them as $\gamma(x)=\sum_s c_s x^s$ and $\gamma'(x)=\sum_s c'_s x^s$, where $s\geq 1$ is not necessarily an integer. define the \textit{exponent of contact of }$\gamma$ and $\gamma'$ as
\begin{equation*}
    \cont (\gamma, \gamma')=\min \,\{s\geq 1: c_s\neq c'_s\},
\end{equation*}
and $\cont (\gamma, \gamma')=+\infty$ just in case that $\gamma=\gamma'$. A straightforward property that follows from the definition is that, if $\gamma,\gamma',\gamma''$ are pro-branches, then 
\begin{equation}\label{eq:contact_pro}
    \cont (\gamma, \gamma'')\geq \min\,\{\cont(\gamma, \gamma'), \cont(\gamma',\gamma'') \}.
\end{equation}
Indeed, let $\kappa=\cont (\gamma, \gamma'')$, and write $\gamma(x)=\sum c_s x^s$, and similarly for $\gamma'$ and $\gamma''$ with coefficients $c_s'$ and $c_s''$. If $\kappa$ is a finite value, it has to be the smallest value satisfying $c_\kappa=c_\kappa''$. Let us distinguish from two possibilities. If $c_\kappa\neq c_\kappa'$, then in particular $\kappa \geq \cont (\gamma, \gamma')$ by the definition of exponent of contact. Now, if $c_\kappa=c_\kappa'$, then $c_\kappa'\neq c_\kappa''$, yielding that $\kappa \geq \cont(\gamma', \gamma'')$. Hence, $\kappa\geq \min \{\cont(\gamma, \gamma'), \cont(\gamma',\gamma'')$ and the result follows. In particular, equation (\ref{eq:contact_pro}) states that the two smaller values of $\cont (\gamma, \gamma'), \cont (\gamma, \gamma''), \cont (\gamma', \gamma'')$ coincide. \\

At first sight, one could think that the exponent of contact concept depends on the choice of coordinates of $(\C^2,0)$, but far from it, it has a clear geometrical expression.
\begin{lemma}\label{lemma:contactinvariant} Let $B$ and $B'$ be representative sets for branch germs, and let $\gamma\in \textup{\pro} (B)$ and $\gamma'\in \textup{\pro}(B)$ with $\kappa= \cont (\gamma, \gamma')$. For $P\in B$, write $d(P,B)$ for the distance from $P$ to $B$, and $d(P,0)$ to the distance from $P$ to the origin. Then, $\lim_{P\to 0} d(P,B)/d(P,0)^\kappa$ is a finite nonzero value.
\end{lemma}
\noindent For a proof, see \cite{Wall}, Lemma 4.1.1. In general, one can define the \textit{exponent of contact} for the branches $(B,0)$ and $(B',0)$ to be the value $$\cont (B,B')=\min \left\{\cont (\gamma, \gamma'): \gamma\in \pro (B), \gamma'\in \pro (B') \right\}.$$ 
It turns out that, in the previous definition, one of the pro-branches may be fixed in virtue of the following result:
\begin{proposition}\label{prop:contact}
    Let $\gamma\in \textup{\pro} (B)$. Then, $\{\cont (\gamma, \gamma') : \gamma'\in \textup{\pro} (B')\}$ does not depend on $\gamma$. 
\end{proposition}
\begin{proof} Notice that, if $\gamma(x)=\sum_{r=1}^{+\infty}a_rx^{r/m}$, the functions $\gamma_k(x)=\sum_{r=1}^{+\infty}a_r e^{2\pi i rk/m} x^{r/m}$ defined for $x\in S$ and $k\in \{0,\dots,m-1\}$ satisfy that $\pro(B)=\{\gamma_0, \gamma_1, \ldots, \gamma_{m-1}\}$. If this set is considered to be ordered and possible repeated elements are allowed, it follows that any other choice of the pro-branch yields $\pro(B)$ to be a cyclic permutation of the previous case. Hence, a change in $\gamma$ would permute the set $\{\cont (\gamma, \gamma') : \gamma'\in \textup{\pro} (B')\}$ cyclically. Thus, the set does not depend on the choice of $\gamma\in \pro(B)$.     
\end{proof}
\noindent Thus, in the calculation of $\cont (B,B')$ one may fix one of the pro-branches (but not both). It would be convenient from now on to consider the set of the previous proposition with possibly repeated elements. \\ 

A straightforward consequence of equation (\ref{eq:contact_pro}) is that 
\begin{equation}\label{eq:contact_branches}
    \cont (B,B')\geq \min \,\{\cont (B,B''), \cont (B',B'')\},
\end{equation}
where $(B,0), (B',0)$ and $(B'',0)$ are generic branches. As before, this claim forces the two smaller of the three involved exponents to be equal. 
\begin{exercise}(4.7.1. of \cite{Wall}) Consider the branches $(B,0)$ and $(B',0)$ defined by the equations $x=t^2, y=t^3$ and $x=t^4, y=t^6+t^7$, respectively. For $(B,0)$, one has that $x^{1/2}$ and $-x^{1/2}$ are the only possible square root functions defined along $\S$. Hence, $\pro (B)=\{\gamma_1, \gamma_2\}$, where $\gamma_1(x)=x^{3/2}$ and $\gamma_2(x)=-x^{3/2}$. For the second one, notice that the fourth root functions are given by $x^{1/4}, ix^{1/4}, -x^{1/4}$ and $-ix^{1/4}$. Hence, $\pro(B')=\{\phi_1, \phi_2, \phi_3, \phi_4\}$ is given by
\begin{align*}
    \phi_1(x)&=x^{3/2}+x^{7/4}, &\phi_2(x)=-x^{3/2}-ix^{7/4}, \\
    \phi_3(x)&=x^{3/2}-x^{7/4}, &\phi_4(x)=-x^{3/2}+ix^{7/4}.
\end{align*}
Hence, $\cont (\gamma_1, \phi_1)=\cont (\gamma_1, \phi_3)=7/4$, $\cont (\gamma_1, \phi_2)=\cont (\gamma_1, \phi_4)=3/2$ and $\cont (\gamma_2, \phi_1)=\cont (\gamma_2, \phi_3)=3/2$, $\cont (\gamma_2, \phi_2)=\cont (\gamma_2, \phi_4)=7/4$. This yields that $\cont (B,B')=\min \{3/2, 7/4\}=3/2$ (however, it would have been enough to find the minimum fixing the pro-branch $\phi_1$ of $B'$, and reducing to 2 the number of exponents of contact to be calculated). \\
\end{exercise}
Since the most important coefficients in a parametric equation for a branch are the $\beta_q$, appearing as $t\,^{\beta_q}$, it follows that the coefficients $\alpha_q=\beta_q/m$ play a key role in terms of pro-branches, for $q\in \{1,\dots g\}$. These terms are known as the \textit{Puiseux exponents} of a branch. Notice that, if two given pro-branches have an exponent of contact $\kappa$ greater than $\alpha_q$, then they share the same first $q$ Puiseux exponents $\alpha_1, \ldots, \alpha_q$. 
\begin{proposition}\label{prop:cont} Let $(B,0)$ and $(B',0)$ be branches with exponent of contact $\kappa =\cont (B,B')$. Write $\alpha_1, \ldots, \alpha_g$ for the Puiseux exponents of $(B,0)$, and $e_1,\dots, e_g$ for the auxiliary terms of the Puiseux characteristic of $(B,0)$. 
    \begin{enumerate}[(i)]
        \item Let $\gamma'\in \textup{\pro} (B')$ and let $q\in \N$ be such that $\alpha_q<\kappa \leq \alpha_{q+1}$. Then, $\{\cont (\gamma, \gamma'): \gamma\in \textup{\pro} (B)\}$ consists of:
        \begin{itemize}
            \item $\alpha_r$ occurring $e_{r-1}-e_r$ times, for $r\in \{1,\dots,q\}$, and 
            \item $\kappa$ occurring $e_q$ times.
        \end{itemize}
        \item If $\gamma'\in \textup{\pro} (B)$, then $\{\cont (\gamma, \gamma'):\gamma\in \textup{\pro}(B), \gamma\neq \gamma'\}$ consists of $\alpha_i$ occurring $e_{i-1}-e_i$ times, for $i\in \{1,\dots,g\}$. 
    \end{enumerate}
\end{proposition}
\begin{proof}
    Let $\gamma\in \pro (B)$ be the pro-branch that satisfies $\cont (\gamma,\gamma')=\kappa$, and write $\gamma(x)=\sum_{j=1}^{+\infty}a_j x^{j/m}$. If $\gamma_k(x)=\sum_{j=1}^{+\infty} a_j e^{2\pi i kj/m} x^{j/m}$ for $k\in \{0,\dots, m-1\}$, then $\pro(B)=\{\gamma_0, \dots, \gamma_{m-1}\}$. Notice that the coefficient associated to $x^{j/m}$ of $\gamma$ and $\gamma_k$ is unchanged if and only if $e^{2\pi i kj/m}=1$, which occurs when $kj/m\in \Z$. If $j=\beta_1$, it follows that the coefficients of $x^{\alpha_1}$ in $\gamma$ and $\gamma_k$ is different if and only if $k\alpha_1 \notin \Z$, and in this case $\cont (\gamma_k, \gamma')=\alpha_1$. Since $e_1=\gcd \{m,\beta_1\}$, then there exists exactly $m-e_1=e_0-e_1$ values of $k$ for which this holds. If $j=\beta_r$ with $1\leq r\leq q$, it follows that the coefficients of $x^{\alpha_r}$ in $\gamma$ and $\gamma_k$ is different but the previous are equal if and only if $k\alpha_r \notin \Z$, but $k\alpha_{r-1}\in \Z$. Therefore, $\cont (\gamma_k, \gamma')=\alpha_r$ if and only if $k$ is divisible by $m/e_{r-1}$ but not divisible by $m/e_r$, which occurs for $(m-e_{r})-(m-e_{r-1})=e_{r-1}-e_r$ values of $k$. This proves the first claim of $(i)$. For the second one, notice that the remaining possible changes in the coefficients of $\gamma_k$ with respect to $\gamma$ are in terms of order greater than $\alpha_{q+1}$. Since $\gamma$ is chosen so that $\cont (\gamma, \gamma')=\kappa=\min \{\cont (\gamma_k, \gamma'):0\leq k \leq m-1\}$, it follows that, for the values of $k$ that are divisible by $m/e_q$, one has that $\cont (\gamma_k, \gamma')=\kappa$, and this holds for $e_q$ values of $k$. The same arguments given for $(i)$ shows $(ii)$ since, in this case, $\kappa=+\infty$. 
\end{proof}
\begin{corollary}\label{prop:Puisuexinvariant}
    The Puiseux characteristic of a branch $(B,0)$ is invariant under holomorphic changes of coordinates. 
\end{corollary}
\begin{proof}
    For $\gamma'\in \pro (B)$, one has that $\{\cont (\gamma, \gamma'): \gamma\in \pro (B)\}$ is invariant by Lemma \ref{lemma:contactinvariant}. Now, the number of entries in this set (counting repetitions), is $e_0-e_g=m-1$. Hence, $m$ is invariant. The different entries in the set are the $\alpha_i$ for $i\in \{1,\dots,g\}$, and hence the $\beta_i=\alpha_i m$ is invariant under holomorphic changes of coordinates.
\end{proof}
It turns out that there is a strong connection between the exponent of contact and the intersection number of branches. Namely:
\begin{proposition}\label{prop:int_cont}
    Let $(B,0)$ and $(B',0)$ be branches. Then, 
    \begin{equation*}
        B.B'=\sum\left\{\cont (\gamma, \gamma'):\gamma\in \textup{\pro} (B), \gamma'\in \textup{\pro} (B')\right\}.
    \end{equation*}
\end{proposition}
\begin{proof}
    Choose coordinates in $(\C^2,0)$ such that $x=0$ is not tangent to either branch. By the results of the preceding chapter, one has that one can find parametric equations for $(B,0)$ of the form $x=t^m, y=\phi(t)$, and $x=u^m, y=\xi(u)$ for $(B',0)$. It follows that an equation for $(B,0)$ is given by $f(x,y)=\prod_{k=0}^{m-1} \big(y-\phi(e^{2\pi i k/m}x^{1/m})\big)$. Writing $\psi_k(x)=\phi(e^{2\pi i k/m}x^{1/m})$ for $k\in \{0, \dots, m-1\}$ yields $f(x,y)=\prod_{k=0}^{m-1} \left(y-\psi_k(x)\right)$ and $\pro (B)=\{\psi_0, \ldots, \psi_{m-1}\}$. It then follows by the definition of the intersection number that $B.B'=\ord\!_u f(u^n, \xi(u))$. Moreover, if we write $\eta(x)=\xi(x^{1/n})$, it follows that $\eta \in \pro(B')$. Since $x=u^n$, it follows that $\ord\!_u$ equals $n\,\ord\!_x$. Hence,
    \begin{equation*}
        B.B'=n\,\ord\!_x \prod_{k=0}^{m-1} \left(\eta(x)- \psi_k(x)\right)=n\sum_{k=0}^{m-1}\ord\!_x \big( \eta(x)-\psi_k(x)\big)=n\sum\big\{ \cont (\eta, \psi_k): 0\leq k \leq m-1\big\},
    \end{equation*}
    and the result follows by Proposition \ref{prop:contact}.
\end{proof}
A straightforward consequence of the previous propositions is an equation relating $B.B'$ and $\cont (B,B')$ in terms of the Puiseux characteristic.
\begin{corollary}\label{cor:contact}
    Let $(B,0)$ and $(B',0)$ be branches and let $\kappa=\cont (B,B')$. If $q\in \N$ is the number for which $\alpha_q<\kappa\leq \alpha_{q+1}$. Then, 
    \begin{equation*}
        B.B'=\dfrac{m(B')}{m(B)}\Big( \beta_1(m-e_1)+ \beta_2(e_1-e_2)+ \ldots + \beta_q (e_{q-1}-e_q)+\kappa e_q m(B)\Big) .
    \end{equation*}
\end{corollary}
\begin{proof} By propositions \ref{prop:cont} and \ref{prop:int_cont}, it follows that
\begin{equation*}
    B.B'=n\Big(\alpha_1(m-e_1)+ \alpha_2(e_1-e_2)+\dots+\alpha_q(e_{q-1}-e_q)+\kappa e_q\Big),
\end{equation*}
and the result follows since $\beta_i=\alpha_i/m$ and $n=m(B'), m=m(B)$.
\end{proof}
\noindent This corollary gives an equivalent formulation of the equisingularity notion for general curves:
\begin{theorem} $(C,0)$ and $(\tilde{C},0)$ are equisingular if and only if they have the same number $n$ of branches, namely $\{(B_1,0), \ldots, (B_n,0)\}$ for $(C,0)$ and $\{(\tilde{B}_1,0), \ldots, (\tilde{B}_n,0)\}$ for $(\tilde{C},0)$, and there exists a permutation $\sigma \in \Sigma_n$ such that $(B_{\sigma(i)},0)$ and $(\tilde{B}_i,0)$ are equisingular as branches, and the exponents of contact $\cont (B_{\sigma(i)},B_{\sigma(j)})$ and $\cont (\tilde{B}_i,\tilde{B}_j)$ coincide.
\end{theorem}
\begin{proof} This follows immediately from the definition of equisingularity of curves, and from Corollary \ref{cor:contact}, it follows that two pairs of branches have the same intersection numbers if and only if they have the same exponents of contact. 
\end{proof}
Let $(B,0)$ be a branch with Puiseux characteristic $(m;\beta_1, \ldots, \beta_g)$. define, for $i\in \{1,\dots, g\}$, the number 
\begin{equation*}
    \barbeta_i=\dfrac{1}{e_{i-1}}\Bigl(m\beta_1+e_1(\beta_2-\beta_1)+ \ldots + e_{i-1}(\beta_i-\beta_{i-1})\Bigr),
\end{equation*}
and let $\barbeta_0=m$. It then follows that $e_i\,\barbeta_{i+1}-e_{i-1}\barbeta_i = e_i(\beta_{i+1}-\beta_i)$ for $i\in \{1,\dots,g-1\}$, and hence 
\begin{equation}\label{eq:betabarra}
    \barbeta_{i+1}-\dfrac{e_{i-1}}{e_i} \barbeta_i=\beta_{i+1}-\beta_i.
\end{equation}
Since $\barbeta_0=m$ and $e_{i-1}/e_i$ are integers, it follows by induction that $\barbeta_i$ is an integer for every $0\leq i \leq g$. This new coefficient gives a more compact form for the equation that relates $B.B'$ and $\cont (B,B')$. Indeed, it is straightforward to check that the formula of Corollary \ref{cor:contact} can be rewritten as 
\begin{equation}\label{eq:intnumber}
    B.B'=e_q\dfrac{m(B')}{m(B)}\Bigl(\barbeta_{q+1}-\beta_{q+1}+m(B)\kappa\Bigr),
\end{equation}
where $\kappa=\cont (B,B')$, $q$ is the natural number that satisfies $\alpha_q<\kappa \leq \alpha_{q+1}$ and the involved Puiseux coefficients are the ones of $(B,0)$. \\

Although it could seem at first glance that the $\barbeta_i$ are an artificial invention, it turns out that not only do they completely encode the equisingularity class of a branch, but they also have a natural way of looking at them. This will be clear after having the semigroup of a branch introduced.
\section{The semigroup of a branch}
This last section is devoted to introduce a complete algebraic invariant for the equisingularity relation, namely the semigroup of a branch. Let $(C,0)$ be a curve defined by a reduced equation of the form $f(x,y)=0$. The quotient ring $\cont_C=\C\{x,y\} / (f)$ is defined to be the \textit{local ring of} $C$. For a single branch, one may use parametric equations $\gamma:(\C,0)\rightarrow (\C^2,0)$ to induce a ring homomorphism $\gamma^*:\C\{x,y\} \rightarrow \C\{t\}$ given by the composition $\gamma^*(g)=g\circ \gamma$. Thus, $\cont_B$ can be identified with the image of $\gamma^*$, a subring of $\C\{t\}$. \\

define the \textit{semigroup of a branch} $(B,0)$ as the set $\S(B)=\{\ord \phi: \phi \in \cont_B\}$, equipped with the sum of natural numbers. Since $\ord \phi+ \ord \psi=\ord (\phi\psi)$ and $1\in \S(B)$, it follows that $\S(B)$ is indeed a semigroup. define $\delta(B)=|\{r\geq 0 : r\notin \S(B)\}|$ as the number of gaps in $\S(B)$. Since the order is invariant under holomorphic changes of coordinates, it turns out that both $\S(B)$ and $\delta(B)$ do not depend on the coordinates of $(\C^2,0)$. \\
\begin{remark} Although $\S(B)$ is called \textit{semigroup}, some Algebra references would call it  \textit{monoid}, since it is a set equipped with an associative binary operation that has an identity element. However, in a singularity theoretic context, it is common to use the former term to refer these algebraic objects.\\
\end{remark}
Some immediate consequences of the definition are the following:
\begin{proposition}\label{prop:semigroup} Let $(B,0)$ be a branch. Then,
\begin{enumerate}[(i)]
    \item $\S(B)$ consists of the intersection numbers $B.\tilde{B}$ for branches $(\tilde{B},0)$ not equal to $(B,0)$. 
    \item For every $i\in \{0,\dots,g\}$, one has that $\barbeta_i\in \S(B)$. 
\end{enumerate}
\end{proposition}
\begin{proof} Let $x=t^m, y=\gamma(t)=\sum_{r=1}^{+\infty}a_rt^r$ be a parametric equation for $(B,0)$. The first claim follows since intersection numbers are the orders in $t$ of compositions $f\circ \gamma$, where $f$ is the defining irreducible function for a generic branch $(\tilde{B},0)$, and this coincides with the image of $f$ by $\gamma^*$. For the second one, let $(B_q^-,0)$ be the branch defined by $x=t^m, y=\sum_{1\leq r <\beta_q} a_rt^r$ (which is in general not injective). By the definition of $\beta_q$, it follows that the powers of $t$ that appear in $y$ are also powers of $t^{e_{q-1}}$. Hence, $B_q^-$ has multiplicity $m/e_{q-1}$ and the expression $x=t^{m/e_{q-1}}, y=\sum_{1\leq r <\beta_q} a_rt^{r/e_{q-1}}$ is a parametric equation for $(B_q^-,0)$ (\textit{i.e.}, it is injective). Therefore, $\kappa=\cont (B,B_q^-)=\beta_q/m=\alpha_q$, and by the equation (\ref{eq:intnumber}), since the exponent of contact lies in the interval $\alpha_{q-1}<\kappa\leq \alpha_q$, it follows that 
\begin{equation*}
    B.B_q^-= e_{q-1}\dfrac{m(B_q^-)}{m(B)}\Bigl(\barbeta_q-\beta_q+m(B)\kappa\Bigr)=e_{q-1}\dfrac{m/e_{q-1}}{m}\Bigl(\barbeta_q-\beta_q+m\alpha_q\Bigr)=\barbeta_q,
\end{equation*}
and the result follows.
\end{proof}
It turns out that much more can be stated:
\begin{theorem} Let $(B,0)$ be a branch. Then,
\begin{enumerate}[(i)]
    \item $\delta(B)$ is finite. 
    \item The powers $t^r$ with $r\notin \S(B)$ form a basis of $\C\{t\}/\cont_B$. Hence, $\dim_\C \C\{t\}/\cont_B=\delta(B)$. 
    \item If $N(\S(B))$ is the greatest integer not lying in $\S(B)$, one has that any element of $\C\{t\}$ with order greater than $N(\S(B))$ is in $\cont_B$. 
    \item $N(S(B))=\sum_{q=1}^g (e_{q-1}-e_q)(\beta_q-1)-1$.
\end{enumerate}
\end{theorem}
\noindent For a proof, see 4.3.1, 4.3.2 and 4.3.3 of \cite{Wall}. Moreover, it can be shown that:
\begin{theorem}\label{theorem:semigroup} $\barbeta_0, \dots, \barbeta_g$ is a minimal set of generators of $\S(B)$. Furthermore, $\barbeta_q$ is the least element of $\S(B)$ not divisible by $e_{q-1}$.  
\end{theorem}
\noindent The proof can be found in 4.3.5 of \cite{Wall}. \\
\begin{example}\label{example} Consider the branch given by $x=t^4, y=t^6+t^7$. Then, its Puiseux characteristic is $(4;6,7)$, and $e_0=4,e_1=2,e_2=1$. An immediate application of equation (\ref{eq:betabarra}) yields $\barbeta_0=4, \barbeta_1=6, \barbeta_2=13$. Hence, $\S(B)$ is the semigroup generated by $4,6$ and $13$. Thus,
$$\S(B)=\{0,4,6,8,10,12,13,14,16,17,18,\dots\}$$
and $N(\S(B))=(4-2)(6-1)+(2-1)(7-1)-1=15$. It can be directly checked that $4,6$ and $13$ lie in $\S(B)$. Indeed, $4$ is the order of the image of $x$ by $\gamma^*$, where $\gamma(t)=(t^4,t^6+t^7)$. Moreover, $6$ is the order of the image of $y$, and a straightforward computation shows that $x^3-y^2$ has an image of order $13$.\\
\end{example}
The last aim of the chapter is to show that both the $\barbeta_i$ and the semigroup of a branch completely determine the equisingularity class of the branch. 
\begin{theorem}\label{thrm:semigroup_determines_puiseux}
    Let $(B,0)$ be a branch not tangent to $x=0$. Then, any of the following data determines the others:
    \begin{enumerate}[(i)]
        \item the Puiseux characteristic,
        \item the terms $\barbeta_0, \ldots, \barbeta_g$, and 
        \item the semigroup $\S(B)$. 
    \end{enumerate} 
\end{theorem}
\begin{proof}
    It is clear by definition that the Puiseux characteristic determines the $\barbeta_i$. As a consequence of Theorem \ref{theorem:semigroup}, one has that $\barbeta_0, \dots, \barbeta_g$ determines $\S(B)$. Conversely, it can be shown by induction that $\gcd\{\barbeta_0, \ldots, \barbeta_q\}=e_q$, for $q\in \{0,\dots, g\}$. Since $\barbeta_0=m$, it follows by (\ref{eq:betabarra}) that the $\beta_q$ can be written in terms of the $\barbeta_0, \dots, \barbeta_q$. This shows that the $\barbeta_i$ determine the Puiseux characteristic. Lastly, the $\barbeta_i$ can be deduced from $\S(B)$. Indeed, $\barbeta_0$ is the least nonzero integer of $\S(B)$, and inductively $\barbeta_q$ is the least element in $\S(B)$ not divisible by $e_{q-1}$, and $e_q=\gcd\{e_{q-1}, \barbeta_q\}$.  
\end{proof}
Therefore, the equisingularity class of a branch is completely determined by any of the three invariants listed above. In the last part of the third chapter, the $\barbeta_i$ will be shown to play a fundamental role when it comes to encode the topological information of a branch. At that moment, the results given in this chapter will flourish and provide the tools to complete the desired topological classification. 
\chapter{The topology of a curve}
Equisingularity provides a merely arithmetic classification for curves, due to the fact that it just involves some calculations made from the exponents appearing in the parametric equations of a branch. In a more topological mindset, it would be of great interest to perform a classification of branches in terms of their topological structure. Notice first that any branch is homeomorphic to the set-germ $(\C,0)$, due to the fact that parametric equations are homeomorphisms. Hence, the intrinsic topological structure of the curves gives no information of them. This makes clear that the classification of plane curves has to be performed by looking at how they are embedded in $(\C^2,0)$. This motivates the following definition:
\begin{definition}
    The representative sets of curve-germs $C$ and $\tilde{C}$ are said to be \textit{topologically equivalent} whenever there exists an $\epsilon_0>0$ such that, for every $0<\epsilon<\epsilon_0$, there is an orientation-preserving homeomorphism of the pair $h:(D_\epsilon,C\cap D_\epsilon)\rightarrow (D_\epsilon, \tilde{C}\cap D_\epsilon)$, where $D_\epsilon=\{(z,w)\in \C^2 : |z|^2+|w|^2\leq \epsilon^2\}$. The previous condition regarding $h$ refers to the existence of a homeomorphism $h:D_\epsilon \rightarrow D_\epsilon$ that preserves the orientation of $D_\epsilon$, and such that $h(C\cap D_\epsilon)=\tilde{C}\cap D_\epsilon$.
\end{definition}
\noindent The previous definition makes clear that the topological class of a curve does not depend on the chosen representative set. In particular, if $C$ and $\tilde{C}$ are topologically equivalent, it follows just by taking germs, that $(C,0)$ and $(\tilde{C}, 0)$ are equivalent in the sense that there exists an orientation-preserving homeomorphism $h: (\C^2,0)\rightarrow (\C^2,0)$ such that $h(C,0)=(\tilde{C},0)$. It is a nontrivial fact that both notions of equivalence coincide, \textit{i.e.} two curves are topologically equivalent if and only if there exists a germ of orientation-preserving homeomorphism $(\C^2,0)\rightarrow (\C^2,0)$ carrying one curve to the other. We will not prove here such a stronger fact (see \cite{le}). \\

The main objective of the chapter is to show the surprising fact that curves are equisingular if and only if they are topologically equivalent. Therefore, the topological classification of plane curves, which has a clear underlying geometrical intuition, can be performed just by looking at some arithmetic data from the curve, namely the Puiseux characteristic of their branches, and the intersection numbers between them. \\

The proof of this main result will be constructed along the whole chapter. In the first section, the link of a curve will be presented as the link arising from the intersection of the curve with a small enough $3$-dimensional sphere. It will be shown that both the topological class of a curve and the isotopy class of its link determine each other. The second section will be devoted to show that the isotopy class of the knot associated to a branch depends only on the Puiseux characteristic of the branch. Hence, the equisingularity class of a branch will be shown to determine the isotopy class of the link. This result will be generalised to general curves with more than one branch, as an application of the exponent of contact's notion. The remaining part of the chapter will be devoted to show the converse statement, namely, that the isotopy type of the knot determines the equisingularity class of the branch. In the third section, the nature of the knots that arise as intersections of branches with $3$-dimensional spheres will be analysed: they will be shown to be iterated torus knots. In the last section, this geometrical nature of the knot will be harnessed to extract the Puiseux characteristic of its corresponding branch. More precisely, the Alexander polynomial will be shown to encode completely the equisingularity class of the branch that is associated to the studied knot. Therefore, the main result of the project will follow in its general version, namely, that \textit{curves are equisingular if and only if they are topologically equivalent.} \\

The origin of these ideas dates back to 1905, when Writinger started thinking about the possible groups that could arise from algebraic knots. He suggested this problem to his PhD student Brauner, who gave the first steps towards the plane curve classification. These correspond to the ideas developed in the first two sections. Brauner described the links via stereographic projections in his paper \cite{Brauner}, in 1928. However, a clearer viewpoint was given by Kähler in \cite{kahler}. The presented ideas in the last section, regarding the usage of the Alexander polynomial to extract the information from the equisingularity class of the curve were noticed independently by Burau and Zariski in 1932 (see articles \cite{burau} and \cite{zariski}).
\vfill
\section{The link of a curve}
Along this first section, our focus is to analyse the intersection of a representative set for a plane curve $C$ with the sphere $\S_\epsilon=\{(z,w)\in \C^2 : |z|^2+|w|^2=\epsilon^2\}$. The topological intuition behind it comes from transversality theory. Indeed, whenever $C$ and $\S_\epsilon$ meet transversely, the real codimension of $C\cap \S_\epsilon$ in $\C^2$ will be $1+2=3$. Thus, the intersection $C\cap \S_\epsilon$ should have the structure of a $1$-dimensional real manifold in $\S_\epsilon$, which is in fact a link immersed the $3$-dimensional sphere. The nature of curves will be shown to be conic, meaning that the topology of $C\cap D_\epsilon$ is the same as the one given by the cone on $C\cap \S_\epsilon$. Therefore, the topological class of a curve will be determined by the kind of knot or link arising in $C\cap \S_\epsilon$. In what follows, our intention will be to bring mathematical consistence to the previous assertions. \\

In order to do so, it will be necessary to make use of vector fields. Recall that, if $\Omega \subset \C^2$ is an open set, a \textit{vector field} is a real-differentiable map $\xi:\Omega \rightarrow \C^2$ (that is, holomorphy is not expected in vector fields, and hence they have to be understood as smooth vector fields in $\R^4$). An \textit{integral curve} of a vector field $\xi$ with basepoint $x\in \Omega$ is a real-differentiable map $\gamma: I\subset \R\rightarrow \Omega$ such that $\gamma'(t)=(\xi\circ \gamma)(t)$ and $\gamma(0)=x$. Moreover, a \textit{flow} of the vector field $\xi$ is a map $\varphi: W\subset \Omega \times \R \rightarrow \Omega$ such that, for a fixed $x\in \Omega$, the map $\varphi(x, \cdot)$ is an integral curve with basepoint $x$. \\

Since vector fields are deeply studied in an undergraduate level, the main results about them are not stated, and the main theorems regarding the existence, uniqueness and maximality of integral curves will be applied without further mention. In spite of this, the following lemma concerning vector fields will be stated and proved, due to the fact that it is not a common result that could be found in an undergraduate degree course.
\begin{lemma}\label{lemma:convexity}
    Let $\Omega \subset \C^2$ be an open set. Assume that, for each $x\in \Omega$, we are given a convex subset of vectors $A_x\subset \C^2$. If there exists an open cover $\{U_k\}_k$ of $\Omega$ and vector fields $\xi_k$ defined in $U_k$ such that for every $x\in U_k$ one has that $\xi_k(x)\in A_x$, then there exists a vector field $\xi$ in $\Omega$ such that $\xi(x)\in A_x$ for every $x\in \Omega$.  
\end{lemma}
\begin{proof}
    Let $\{\varphi_k\}_k$ be a differentiable partition of the unity such that $\phi_k$ is zero outside $U_k$, but also in a neighbourhood of its boundary. Then, the product $\varphi_k\xi_k$ is a vector field on $U_k$ vanishing in a neighbourhood of the boundary of $U_k$. Thus, it can be extended to a vector field in $\Omega$ defining it as $0$ outside $U_k$. Hence, $\xi=\sum \varphi_k \xi_k$ is a vector field in $\Omega$, since for every $x\in \Omega$, there exists a neighbourhood where the previous sum is finite, and hence well-defined. Now, if $\varphi_k(x)\neq 0$, in particular $x\in U_k$ and therefore $\xi_k(x)\in A_x$. Thus, $\xi(x)$ is a finite sum of elements in $A_x$, whose coefficients sum up $1$. The convexity of $A_x$ implies that $\xi(x)\in A_x$ for every $x\in \Omega$. 
\end{proof}
As it will be clarified later on, it will be convenient to consider the \textit{cone} of a topological space $X$ as the quotient topological space given by 
\begin{equation*}
    \mathcal{C}(X)=\dfrac{X\times [0,\epsilon]}{X\times \{\epsilon\}}.
\end{equation*}
Hence, $\mathcal{C}(X)$ is the quotient space $X\times[0,\epsilon]$, where the copy of $X$ given by $X\times\{\epsilon\}$ is collapsed to a point. \\

In any research area, or any research where a major result is achieved, there are crucial moments that could be called turning points, where the understanding of the treated objects is deepened. Theorem \ref{thrm:vector} is one of these turning points. Before this result, it seemed impossible to establish a relationship between the topological nature of a curve with its equisingularity class, which is completely described by arithmetic invariants. After this theorem, the knot theory reinvigorates the understanding that we have of plane curves, providing an unexpected way to look at them. Not only does Theorem \ref{thrm:vector} state that $C\cap \S_\epsilon$ is a union of embedded copies of the $\S^1$, which is a link in the knot-theoretical language, but it also enlightens us on how to perform the classification of plane curves: just by studying these arising links.
\begin{theorem}\label{thrm:vector}
    Let $C\subset \C^2$ be a representative set of a curve-germ. Then, there exists an $\epsilon_0>0$ such that, for every $0<\epsilon<\epsilon_0$, the intersection $C\cap \S_\epsilon$ is a 1-dimensional real manifold smoothly embedded in $\S_\epsilon$, and there is a homeomorphism of the pair $(D_\epsilon, C\cap D_\epsilon)$ to the cone on $(S_\epsilon, C\cap \S_\epsilon)$ that preserves orientations. 
\end{theorem}
\begin{remark} 
        Theorem \ref{thrm:vector} states that a curve $C$ can be visualised as the cone of $C\cap \S_\epsilon$. Figure 3 shows the case of having $C\cap \S_\epsilon$ as a knotted copy of $\S^1$. Outside the origin, a curve is a $1$-dimensional complex manifold, and hence a real surface. However, one could think at first glance that the previous picture does in fact have self-intersections. This would contradict the manifold structure that a curve has outside the origin. \\
        Indeed, this illustration in $\R^3$ does intersect itself, although these intersections are 1-dimensional real manifolds. Hence, if the extra dimension of $\C^2=\R^3\times \R$ is thought as a colouring, the curve can be painted in such a way that, along the ``intersections'', both regions have different colours. This makes clear that the picture does not have self-intersections in $\C^2$ whenever we remove the origin from it.
    \end{remark}
    \begin{figure}
    \centering
    \includegraphics[width=0.4\textwidth]{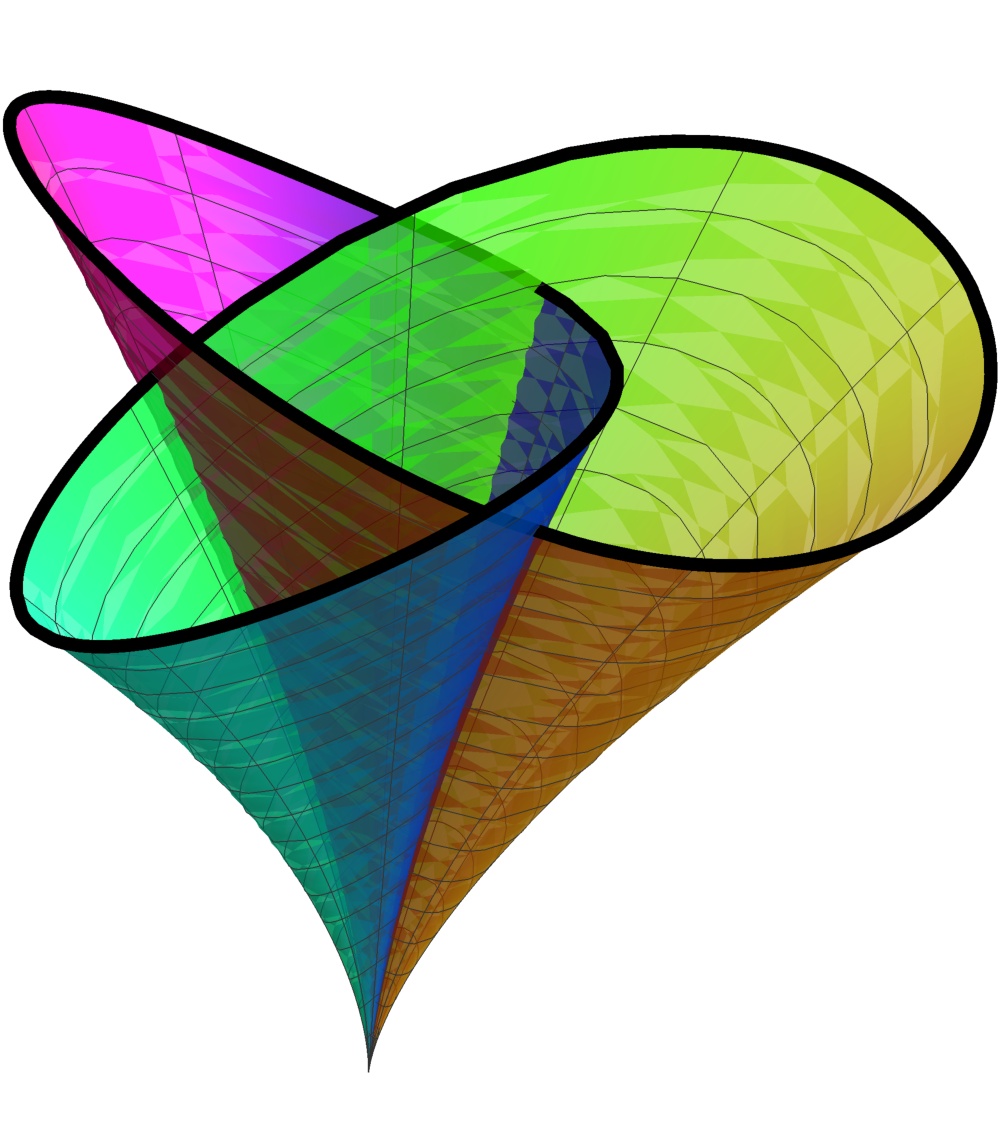}
    \caption{the curve $C$ given by $x^2-y^3=0$, and the intersection $C\cap \S_\epsilon$ shown in black.}
    \end{figure}
\begin{proof}[Proof of Theorem \ref{thrm:vector}] Let $\epsilon'>0$ be a radius that is small enough to have the origin as the unique singularity in $C\cap B_{\epsilon'}$. Let $r:B^*_{\epsilon'}\rightarrow (0,\epsilon)$ be the radius function $r(x)=|x|$, and $\rho(x)=x/|x|$ its gradient vector field. \\

The procedure will be to construct a vector field $\xi$ in $B^*_{\epsilon'}$ with $\langle \xi, \rho\rangle =-1$ and such that, for every $x\in C$, $\xi(x)$ is tangent to $C$ at $x$. Once again, it is important to keep in mind that vector fields have to be understood in $\R^4$. Hence, the inner product is the one in $\R^4$, and the tangent space of $C$ has to be seen as a 2-dimensional real vector space. \\

Assume first that such a vector field $\xi$ exists. Let us show that the maximal flow of $\xi$ is enough to construct the desired homeomorphism. Let $x\in B^*_{\epsilon'}$, and consider the maximal integral curve $\gamma:I\rightarrow B^*_{\epsilon'}$ of $\xi$ that passes through $x$ if $t=0$. We claim that
\begin{enumerate}[(1)]
    \item $(r\circ \gamma)(t)=-t+r(x)$ for every $t\in I$,
    \item $I=(r(x)-\epsilon', r(x))$,
    \item $\lim_{t\to r(x)^-}\gamma(t)=0$,
    \item If $x\in C$, then $\gamma(t)\in C$ for every $t\in I$. 
\end{enumerate}
    Indeed, for every $t\in I$, one has that 
    \begin{equation*}
        (r\circ \gamma)'(t)=\langle (\nabla r\circ \gamma) (t), \gamma'(t)\rangle =\langle (\rho\circ\gamma)(t), (\xi\circ\gamma)(t)\rangle =-1.
    \end{equation*}
    Hence, since $(r\circ\gamma)(0)=r(x)$, it follows that $(r\circ \gamma)(t)=-t+r(x)$, proving (1). In order to check (2), notice that $\gamma$ is well-defined whenever $\gamma(t)\in B_{\epsilon'}^*$, which means that $0<(r\circ \gamma)(t)<\epsilon'$. As a consequence of (1), this is in turn equivalent to have that $r(x)-\epsilon'<t<r(x)$. This proves (2). Condition (3) is a straightforward consequence of both (1) and (2), due to the fact that $\lim_{t\to r(x)^-}(r\circ \gamma)(t)=-r(x)+r(x)=0$. Hence, $\lim_{t\to r(x)^-}\gamma(t)=0$. Lastly, (4) is followed by the fact that the restriction of $\xi$ to $C$ is a well-defined vector space, since $\xi$ is tangent to $C$. Therefore, if $x\in C$, the maximal integral curve $\gamma$ is also an integral curve of the restricted field. \\

    It follows from the previous assertions that, if $\varphi:W\rightarrow B^*_{\epsilon'}$ is the maximal flow of $\xi$, then 
    \begin{equation*}
        W=\{(x,t)\in B_{\epsilon'}^*\times \R: r(x)-\epsilon'<t<r(x)\}.
    \end{equation*}
    Let $0<\epsilon<\epsilon'$, and let $D_{\epsilon}=\{(z,w)\in \C^2: |z|^2+|w|^2\leq \epsilon^2\}$ be the closed ball in $\C^2$. Notice that $S_\epsilon\times [0,\epsilon)\subset W$ trivially, and that $\varphi(S_\epsilon\times [0,\epsilon))\subset D^*_\epsilon$ by condition (1). Hence, the restriction of $\varphi$ given by $f:S_\epsilon \times [0,\epsilon)\rightarrow D_{\epsilon}^*$ is a well-defined and smooth map. In addition, $f$ is a diffeomorphism, since the map $g:D^*_\epsilon\rightarrow S_\epsilon \times [0,\epsilon)$ given by $g(x)=(\varphi(x,r(x)-\epsilon), \epsilon-r(x))$ is the inverse of $f$. Indeed, $g$ is a smooth map, and it satisfies that 
    \begin{equation*}
        (g\circ f)(x,t)=(\varphi(\varphi(x,t), -t+r(x)-\epsilon), \epsilon+t-r(x))= (\varphi(\varphi(x,t),-t),t)=(x,t)
    \end{equation*}
    as a consequence of (1) and the properties of the flow. Moreover,  
    \begin{equation*}
        (f\circ g)(x)=\varphi(\varphi(x,r(x)-\epsilon), \epsilon-r(x))=x.
    \end{equation*}
    Thus, $f$ is a diffeomorphism with inverse $g$. In addition, $\lim_{t\to \epsilon^-}f(x,t)=0$ for every $x\in \S_\epsilon$ by the third condition stated before. Moreover, $f((C\cap \S_\epsilon)\times [0,\epsilon)) = C\cap D_\epsilon^*$ as a consequence of condition (4). \\

    Thanks to the map $f$, it is now possible to define, for $\mathcal{C}(\S_\epsilon)=(\S_\epsilon\times [0,\epsilon])/(\S_\epsilon\times \{\epsilon\})$, the map $\tilde{f}:\mathcal{C}(S_\epsilon)\rightarrow D_\epsilon$ by $\tilde{f}([x,t])=f(x,t)$ whenever $0\leq t <\epsilon$, and $\tilde{f}([x,\epsilon])=0$. This clearly yields a well-defined map, and since $\lim_{t\to \epsilon^-}f(x,t)=0$, then $\tilde{f}$ is continuous. Furthermore, since $f$ is bijective,  the same holds for $\tilde{f}$. This is enough to grant that $\tilde{f}$ is a homeomorphism, due to the fact that it is a continuous bijection, $\mathcal{C}(\S_\epsilon)$ is a compact space and $D_\epsilon$ is Hausdorff. In addition, $\tilde{f}(\mathcal{C}(C\cap \S_\epsilon))=D_\epsilon\cap C$. Lastly, the orientation-preserving clause for $\tilde{f}$ follows since the homeomorphism comes from a diffeomorphism $f$ arising as the flow of a vector field. Hence, the restriction of $\tilde{f}$ to $\S_\epsilon \rightarrow \S_\epsilon \times\{0\}\subseteq \mathcal{C}(\S_\epsilon)$ is the identity map $x\rightarrow [(x,0)]$. This proves the theorem. \\

    It only remains to check the existence of such a vector field $\xi$ in $B^*_{\epsilon'}$, namely, with $\langle \xi, \rho\rangle =-1$ and such that $\xi(x)$ is tangent to $C$ if $x\in C$. Notice first that the impositions for $\xi(x)$ define a convex subset of possible vectors. Indeed, if $v,w\in \C^2$ satisfy that $\langle v,\rho\rangle=\langle w,\rho\rangle=-1$ and $\alpha, \beta$ are positive real values with $\alpha+\beta=1$, then $\langle \alpha v+\beta w, \rho\rangle =-\alpha-\beta=-1$. Moreover, if $v$ and $w$ are tangent vectors to $C$ at $x\in C$, then the same holds for $\alpha v+\beta w$. Hence, by Lemma \ref{lemma:convexity}, it is enough to construct the vector field $\xi$ in a neighbourhood of every $x\in B_\epsilon'$. Assume first that $x\notin C$. Since $C$ is closed, there exists a neighbourhood of $x$ that does not contain points of $C$. Hence, for every point $y$ near $x$ define $\xi(y)=-y/|y|=-\rho(y)$. This clearly satisfies that $\langle \xi(y), \rho(y)\rangle=-1$. If $x\in C$ and $\rho(x)$ is not perpendicular to $T_xC$, then there exists a neighbourhood $U$ of $x$ such that, for every $y\in U$, the radius $\rho(y)$ is not perpendicular to $T_yC$. Consider the orthogonal projection of $y$ in $T_yC$, which is clearly a nonzero vector. Then, multiply this vector by a real number, in such a way that the resulting vector $\xi(y)$ satisfies $\langle \xi(y), \rho(y)\rangle =-1$. This defines $\xi$ in the open neighbourhood $U$.  Lastly, it was discussed in \ref{theorem:submanifold} that the radius vector $\rho(x)$ defines a complex line convergent to the complex tangent line at the origin. In addition, the complex tangent lines of nonzero points also converge to the tangent line at 0. Therefore, it is possible to shrink the ball $B_{\epsilon'}$ if necessary to have that, for every $x\in B_{\epsilon'}^*$, the radius vector $\rho(x)$ is not perpendicular to $T_xC$. \\

    The clause regarding that $C\cap \S_\epsilon$ is a regular 1-dimensional real submanifolds follows from the fact that both $C\setminus \{0\}$ and $\S_\epsilon$ are regular submanifold of $\C^2$, with real codimensions 2 and 1, respectively. In addition, for every $x\in C\cap \S_\epsilon$, $\xi(x)$ is a tangent vector of $C$ at $x$ and $\langle \xi(x), \rho(x)\rangle=-1\neq 0$. Since $T_x\S_\epsilon$ consists of the perpendicular vectors to $\rho(x)$, it follows that $C$ and $\S_\epsilon$ meet transversely. Therefore, $C\cap \S_\epsilon$ is a submanifold of $\C^2$ with real codimension $\codim_\R\, C + \codim_\R\, \S_\epsilon = 2+1=3$, and hence real dimension one.  
\end{proof}
\begin{remark}
    Notice that, for the map $\tilde{f}$ of the previous lemma, the points on the cone $\mathcal{C}(\S_\epsilon)$ with height $0<\delta<\epsilon$ are mapped to points in $\S_\delta$. Therefore, $(\S_\epsilon, C\cap \S_\epsilon)$ and $(\S_\delta, C\cap \S_\delta)$ are homeomorphic. This can be informally stated just by saying that the pair $(\S_\epsilon, C\cap \S_\epsilon)$ is essentially independent of $\epsilon$, provided $\epsilon$ is small enough. \\
\end{remark}
The most significant consequence of the previous result is that it is enough to compare the intersection of a curve with an small enough sphere in order to determine its topological class. 
\begin{corollary}\label{corollary:conic}
    Let $C$ and $\tilde{C}$ be representatives of curve-germs. Then, they are topologically equivalent if and only if there exists an $\epsilon_0>0$ such that for every $0<\epsilon<\epsilon_0$ the pairs $(\S_\epsilon, C\cap \S_\epsilon)$ and $(\S_\epsilon, \tilde{C}\cap \S_\epsilon)$ are homeomorphic through an orientation-preserving homeomorphism. 
\end{corollary} 
\begin{proof}
    The direct implication is straightforward by the definition of the equivalence for curves, since a homeomorphism $D_\epsilon\rightarrow D_\epsilon$ restricts to a homeomorphism $\S_\epsilon\rightarrow \S_\epsilon$. This follows from homology, since if there were a point $x\in \S_\epsilon$ with $h(x)\in B_\epsilon$ or \textit{vice versa}, the restricted homeomorphism $h:D_\epsilon\setminus \{x\}\rightarrow D_\epsilon \setminus \{h(x)\}$ would induce an isomorphism of the homology groups $0\rightarrow \Z$ or $\Z\rightarrow 0$, but this it not possible.\\ 
    For the converse implication, notice that the previous theorem gives a homeomorphism between the pair $(D_\epsilon,C\cap D_\epsilon)$ and the cone of $(\S_\epsilon, C\cap \S_\epsilon)$, and the same for $D$. The result follows since a homeomorphism of the pairs $(\S_\epsilon, C\cap \S_\epsilon)\rightarrow(\S_\epsilon, D\cap \S_\epsilon)$ extends to a homeomorphism of the cones. The orientation claims are clear, since every homeomorphism that appears preserves the orientations. 
\end{proof}
In order to properly understand $C\cap \S_\epsilon$, consider first the intersection of $\S_\epsilon$ with a branch. 
\begin{proposition}
    If $B$ is a representative of a branch-germ, then there exists an $\epsilon_0>0$ such that, for every $0<\epsilon<\epsilon_0$, the intersection $B\cap \S_\epsilon$ is homeomorphic to the 1-sphere $\S^1$. 
\end{proposition}
\begin{proof} 
    Let $\epsilon>0$ be small enough to have that $B$ and $\S_\epsilon$ meet transversely, and that $B\cap D_\epsilon$ is homeomorphic to the cone on $K_\epsilon=B\cap \S_\epsilon$. Let $\varphi:U\subset \C \rightarrow \C^2$ be a parametric equation for $B\cap D_\epsilon$, where $U$ is a small enough open neighbourhood of the origin. \\
    Since $\varphi$ is immersive in $U\setminus \{0\}$, in particular it is transversal to $\S_\epsilon$. Therefore, if $g:U\rightarrow \R$ is defined by $g(t)=|\varphi(t)|^2$, one has that $\epsilon^2$ is a regular value of $g$. By the regular value theorem, it follows that $g^{-1}((-\infty,\epsilon^2])=\varphi^{-1}(B\cap D_\epsilon)$ is a real submanifold of $\C$ with boundary $g^{-1}(\epsilon^2)=\varphi^{-1}(B\cap \S_\epsilon)$. Now, notice that $\varphi$ is a homeomorphism, and that $B\cap D_\epsilon$ is a contractible space (since it is homeomorphic to the cone on $K_\epsilon$). Hence, $\varphi^{-1}(B\cap D_\epsilon)$ is a compact contractible surface with boundary $\varphi^{-1}(B\cap \S_\epsilon)$. By the topological classification of compact surfaces with boundary, it follows that $\varphi^{-1}(B\cap D_\epsilon)$ is homeomorphic to the closed 2-dimensional disk and its boundary is homeomorphic to $\S^1$. Since $\varphi$ is a homeomorphism, it follows that $B\cap D_\epsilon$ is homeomorphic to $\S^1$. 
\end{proof}
\noindent Therefore, $B\cap \S_\epsilon$ is a compact manifold homeomorphic to $\S^1$. In general, $C\cap \S_\epsilon$ consists of a finite disjoint union of embedded copies of $\S^1$, one for each branch. Therefore, the intersection of a branch or a curve with a small enough $3$-dimensional sphere yields what we call a knot or a link in the knot-theoretical framework. \\

Recall that an \textit{embedding} between the (real) manifolds $X$ and $Y$ is a smooth injective map $f:X\rightarrow Y$ that is immersive at every point $x\in X$ (\textit{i.e.}, $df_x:T_xX \rightarrow T_{f(x)}Y$ is injective), and whose restriction onto its image $f:X\rightarrow f(X)$ is a homeomorphism. Provided $X$ is compact, it turns out that the third condition is superfluous, since, in this case, every continuous injective map is a homeomorphism onto its image.  
\begin{definition}
    A \textit{knot} in the 3-dimensional sphere $\S^3$ is a subset $K\subset \S^3$ that can be described as the image of an embedding $f:\S^1\rightarrow \S^3$. Similarly, a \textit{link} in $\S^3$ is a subset $L\subset \S^3$ that can be described as the image of an embedding $f$ from a finite disjoint union of copies of $\S^1$, to $\S^3$. Then, a link is a disjoint union of a finite number of knots.
\end{definition}
\noindent With the previous concepts it is now clear that, under the identification of $\S_\epsilon$ with $\S^3$ through a homothety, $C\cap \S_\epsilon$ is a link with as many connected components as there are branches in $C$. The particular knots and links that arise from plane curves are known as \textit{algebraic links}. \\


The equivalence of knots and links is performed from an extrinsic viewpoint: is is studied how are they embedded in $\S^3$. Since a knot is homeomorphic to $\S^1$, there is no relevant information in its intrinsic topological structure. The following  concept is fundamental for the definition of link equivalence. 
\begin{definition}
    Let $I=[0,1]$ and let $X$ and $Y$ be real manifolds. An \textit{isotopy} is a smooth map $F:X\times I \rightarrow Y\times I$ of the form $F(x,t)=(f_t(x),t)$, where for each $t\in [0,1]$ the maps $f_t:X\rightarrow Y$ are embeddings.
\end{definition}
\noindent Therefore, the equivalence of knots and links will be performed through the isotopy concept.
\begin{definition}
    The links $L_0$ and $L_1$ are said to be isotopic whenever they have the same number of connected components, say $n$, and if $X$ denotes a disjoint union of $n$ copies of $\S^1$, then there exists an isotopy $F:X\times I\rightarrow \S^3\times I$ such that $f_0(X)=L_0$ and $f_1(X)=L_1$. \\
\end{definition}
Intuitively, two links $L_0$ and $L_1$ are isotopic if there is a one-parametric family of links $L_t$ between them, for $t\in [0,1]$, that ``depends smoothly'' on $t$. The following result states that, in this case, there is a family of diffeomorphisms $h_t:\S^3\rightarrow \S^3$ carrying $L_0$ to $L_t$ for $t\in [0,1]$. In other words:
\begin{lemma}[Isotopy extension theorem] If $F:X\times I\rightarrow Y\times I$ is an isotopy and $F(x,t)=(f_t(x),t)$, then there exists a diffeomorphism $H:Y\times I\rightarrow Y\times I$ of the form $H(y,t)=(h_t(y),t)$ such that $f_t=h_t\circ f_0$. 
\end{lemma}
\noindent The proof can be found in \cite{Wall}, lemma 5.1.6. It consists of an application of the vector field integration technique. \\ 

The following proposition offers a useful result regarding different alternative definitions one could give for the equivalence of knots and links. Namely, 
\begin{proposition}\label{prop:equivalences} Let $L_0$ and $L_1$ be links. The following are equivalent:
\begin{enumerate}[(1)]
    \item $L_0$ and $L_1$ are isotopic. 
    \item There exists a diffeomorphism $H:\S^3\times I\rightarrow \S^3\times I$ of the form $H(y,t)=(h_t(y),t)$ such that $h_0=id_{\S^3}$ and $h_1(L_0)=L_1$. 
    \item There exists an orientation-preserving diffeomorphism $h:\S^3\rightarrow \S^3$ such that $h(L_0)=L_1$. 
\end{enumerate}
\end{proposition}
\begin{proof}
    Notice first that (1) implies (2) as an immediate consequence of the Isotopy extension theorem. Moreover, (2) implies (3) trivially, since the diffeomorphism $h_1$ preserves orientations as it is isotopic to the identity. On the other hand, (3) implies (2) as a consequence of the Smale conjecture (shown by Hatcher in \cite{hatchersmale} in 1983), which states that $O(4)$ is a deformation retract of the space of diffeomorphisms $\S^3\rightarrow \S^3$. Hence, any orientation-preserving diffeomorphism $h:\S^3\rightarrow \S^3$ is smoothly isotopic to the identity map. It follows that there exists a smooth isotopy $H:\S^3\times I\rightarrow \S^3\times I$ of the form $H(y,t)=(h_t(y),t)$, and such that $h_0=id_{\S^3}$ and $h_1=h$. It only remains to check that (2) implies (1). If $f:X\rightarrow \S^3$ is a smooth embedding (where $X$ is a disjoint union of $\S^1$) such that $f(X)=L_0$, it follows that the map $F=H\circ (f\times id_I):X\times I\rightarrow \S^3\times I$ given by $f_t=h_t\circ f$ is a smooth isotopy satisfying $f_0(X)=(h_0\circ f)(X)=f(X)=L_0$ and $f_1(X)=(h_1\circ f)(X)=h_1(L_0)=L_1$. 
\end{proof}

In a knot-theoretical framework, it is quite common to introduce knots as images of embeddings that are not required to be differentiable, \textit{i.e.}, as topological spaces that are homeomorphic to $\S^1$. In this case, knots or links are said to be isotopic if there exists a homeomorphism $H:\S^3\times I\rightarrow \S^3\times I$ of the form $H(y,t)=(h_t(y),t)$ such that $h_0=id_{\S^3}$ and that $h_1$ carries one link to the other (this is commonly known as an \textit{ambient isotopy}). However, some pathological kind of knots arise in this general setting, having an infinite number of crossings with themselves, namely, the \textit{wild} knots. In order to avoid them, \textit{tame} or \textit{polygonal} knots are the kind of knots that are isotopic to a closed polygonal curve in $\S^3$. It turns out that a similar result to the given in the previous proposition holds for polygonal knots, regarding that polygonal knots or links are isotopic if and only if there exists a homeomorphism $\S^3\rightarrow \S^3$ preserving the orientations and sending one link to the other (see \cite{burde}, Proposition 1.10). \\ 

It turns out that this connection between both definitions for knots does not end up here. In fact, it can be shown that any link in the differentiable sense is ambient isotopic to a polygonal link, and the reciprocal statement does also hold. Then, the kind of knots arising in both parallel theories are completely equivalent. For a reference, see \cite{cromwell}, Theorems 1.11.6 and 1.11.7. The same arguments that shed light to the above are also valid to show that links in the differentiable sense are (smoothly) isotopic if and only if its corresponding polygonal versions are (continuously) ambient isotopic. \\ 

In terms of the above, the following result offers a stronger result compared to the one given in \ref{prop:equivalences}, and becomes crucial in the connection of curves and their associated links. 
\begin{proposition}
    The links $L_0, L_1\subset \S^3$ are isotopic if and only if the pairs $(\S^3, L_0)$ and $(\S^3,L_1)$ are homeomorphic through an orientation-preserving homeomorphism. In other words, $L_0$ and $L_1$ are isotopic if and only if there exists an orientation-preserving homeomorphism $h:\S^3\rightarrow \S^3$ such that $h(L_0)=L_1$. 
\end{proposition}
\begin{proof} 
The direct implication is trivial by Proposition \ref{prop:equivalences}. For the converse, assume the existence of an orientation-preserving homeomorphism $h:(\S^3,L_0)\rightarrow (\S^3,L_1)$. Then, if $\bar{L}_0$ and $\bar{L}_1$ are polygonal links that are ambient isotopic to $L_0$ and $L_1$, respectively, in particular there exists homeomorphisms $(\S^3,L_i)\rightarrow (\S^3, \bar{L}_i)$ that preserve orientations. Hence, $\bar{L}_0$ and $\bar{L}_1$ are ambient isotopic as polygonal links. This forces $L_0$ and $L_1$ to be smoothly isotopic.      
\end{proof}
These knot-theoretical notions let us reinterpret the remark given after the proof of Theorem \ref{thrm:vector} to state that there exists an $\epsilon_0>0$ such that, for every $0<\epsilon<\epsilon_0$, the isotopy type of the link $C\cap \S_\epsilon$ is independent of $\epsilon$. In this case, any representative link in this isotopy class will be called the \textit{link of the curve} $(C,0)$. Moreover, by Corollary \ref{corollary:conic}, this isotopy class completely determines the topological class of the curve. In other words, two curves are topologically equivalent if and only if their associated links are isotopic. \\

Since our intention is to analyse the intersections of curves with a small enough sphere, it will be mandatory to give a parametric equation for a knot $B\cap \S_\epsilon$ by means of an equation for a branch $B$. The following construction was done first by Kähler, in article \cite{kahler}. Let $B$ be a representative set of a branch, and assume that its tangent line is given by $y=0$ (which imposes no restriction, since this can be achieved through a linear change of coordinates). Then, $B$ can be parametrised by $x=t^m, y=\sum_{r>m}a_rt^r$. Since $t^{-m}y$ converges to $0$ as $t$ does, one can take an $\epsilon>0$ such that, if $|x|<\epsilon$, then $|t^{-m}y|<1$, and hence $|y|<|x|<\epsilon$. Therefore, it is possible to replace $S_\epsilon$ by $\S'_\epsilon=\{(x,y)\in \C^2: |x|=\epsilon,\, |y|<\epsilon\}$, which yields an immediate parametric equation for $B\cap \S'_\epsilon$ just by letting $t=\epsilon^{1/m}e^{i\theta}$. Indeed, in this case $|x|=|t|^m=\epsilon$ and therefore the substitution of $t$ in terms of $\theta$ in the parametric equations for $B$ provides parametric equations for $B\cap \S'_\epsilon$ given by
\begin{equation*}
    x=\epsilon e^{im\theta}, y=\sum_{r>m}a_r\epsilon^{r/m}e^{ir\theta}.
\end{equation*}
Furthermore, it can be shown that the replacement of $\S_\epsilon$ with $\S'_\epsilon$ does not affect the isotopy class of the knot in virtue of the following proposition:
\begin{proposition}
    There exists an $\epsilon_0>0$ such that, for every $0<\epsilon<\epsilon_0$, $(\S_\epsilon, B\cap \S_\epsilon)$ is homeomorphic to $(\S'_\epsilon, B\cap \S'_\epsilon)$ through an orientation-preserving homeomorphism.
\end{proposition}
\begin{proof} The same proof of Theorem \ref{thrm:vector} works in this case just by replacing the radius vector $(x,y)$ in the construction of the vector field with the vector $(x,0)$. 
\end{proof}
The last aim of the section will be to introduce the linking number between knots. Turning back to knot theory, the following result shows the existence of the so-called Seifert surfaces in the particular case of algebraic knots.
\begin{theorem}\label{theorem:seifert} Let $K$ be an algebraic knot. Then, there exists an oriented $2$-dimensional real manifold $X$ with boundary and smoothly embedded in $\S^3$, whose boundary is precisely $\partial X=K$.
\end{theorem}
\begin{proof}
    Let $B$ be a representative set of a branch defining $K$, and let $f$ be the irreducible function associated with $(B,0)$. In particular, there exists an $\epsilon>0$ such that $K=B\cap \S_\epsilon$ and that $B$ and $\S_\epsilon$ meet transversely. define $X_u=\{(x,y)\in \C^2: f(x,y)=u\}$ and $K_u=X_u\cap \S_\epsilon$, for $u\in \C$ small enough. Since $\S_\epsilon$ is compact and $f$ has continuous derivatives, one has that there exists $\delta>0$ such that, for every $u\in \C$ with $|u|<\delta$, $X_u$ meets $\S_\epsilon$ transversely. Hence, $K_u$ is a knot isotopic to $K$. Now, by Sard's Theorem, the set of regular values of $f$ has a null complement. Therefore, there exists a regular value $u\in \C$ with $|u|<\delta$. As a result, $X_u \cap D_\epsilon$ is a one-dimensional complex manifold smoothly embedded $\C^2$ with boundary $K_u$. Lastly, an isotopy between $K_u$ and $K$ carries $X_u$ to a surface $X$ with boundary $K$. 
\end{proof}
Such a surface $X$ as in the previous theorem is called a \textit{Seifert surface} for $K$. The reader with some knot-theoretical background may feel surprised with the proof given above. Indeed, for a general knot (\textit{i.e.}, not necessarily algebraic), the Seifert surface is commonly constructed through Seifert's algorithm. In the particular case of algebraic knots, however, there is a special election for such a surface, which is known as the \textit{Milnor fibration} of the branch defined by $f$, and corresponds to the regular surface described by the equation $f(x,y)=u$, where $u$ is a small enough regular value of $f$. This slight perturbation of the branch generates a regular real surface, and a great amount of information is encoded in its homology. However, for the sake of brevity, this is the only mention that is done concerning Milnor fibrations. With the Seifert surfaces being presented, it is now possible to introduce the linking number concept.
\begin{definition}
    Let $K_1, K_2\subset \S^3$ be algebraic knots, and let $X$ be a Seifert surface of $K_1$ that meets $K_2$ transversely. define the \textit{linking number} of $K_1$ and $K_2$ as the topological intersection number between $X$ and $K_2$, and denote it as $Lk(K_1,K_2)=I(X,K_2)$. 
\end{definition}
\noindent Notice that there is no restriction to demand transversality between $X$ and $K_2$, since a slight deformation of any Seifert surface of $K_1$ could be performed to achieve the transversality condition. Moreover, the previous definition does not depend on the chosen Seifert surface $X$. An intuitive sketch for this is that, if $\tilde{X}$ is any other oriented Seifert surface and transversal to $K_2$, then $\tilde{X}-X$ defines, up to a smooth deformation near $K_1$, a compact manifold. This forces $I(\tilde{X}-X, K_2)=0$, and thus $I(\tilde{X},K_2)=I(X,K_2)$. For a more formal justification, we refer the reader to \cite{rolfsen}, page 132. In the reference, eight equivalent definitions of linking number are given. Moreover, it is shown that the linking number is invariant under the action of ambient isotopies in $\S^3$. \\

In what follows, our intention will be to give a connection of the intersection number of branches and the linking number of its corresponding knots. As it is suggested in the given reference, the linking number can be obtained as the following lemma claims. 
\begin{lemma}
    Let $K_1,K_2\subset \S^3$ be algebraic knots. Since $\S^3$ is the boundary of $D^4=\{(x,y)\in \C^2: |x|^2+|y|^2\leq \epsilon^2\}$, consider Seifert surfaces $X_1$ and $X_2$ of $K_1$ and $K_2$, respectively, embedded in $D^4$. If $X_1$ and $X_2$ meet transversely, then $Lk(K_1, K_2)=I(X_1,X_2)$.
\end{lemma}
\noindent In particular, it follows that the linking number is indeed symmetric. Furthermore, it can be shown that the intersection number of branches coincides with the linking number of its corresponding knots. 
\begin{corollary}\label{cor:linking}
    Let $(B_1,0)$ and $(B_2,0)$ be branches, and consider $K_1$ and $K_2$ their associated knots, respectively. Then, $B_1.B_2=Lk(K_1,K_2)$. 
\end{corollary}
\begin{proof}
    Let $f_1,f_2$ be irreducible branches defining $(B_1,0)$ and $(B_2,0)$, respectively. Let $U\subset \C$ be a neighbourhood of the origin, and take representative sets $B_1, B_2$ for the branches in $U$. Let $\epsilon>0$ be small enough that $K_1=B_1\cap \S_\epsilon$ and $K_2=B_2\cap \S_\epsilon$ are the corresponding knots of $(B_1,0)$ and $(B_2,0)$, respectively. Consider the mapping $F:U\rightarrow \C^2$ given by $F=(f_1,f_2)$, and let $X_1=\{(x,y)\in U : f_1(x,y)=u\}$ and $X_2=\{(x,y)\in U : f_2(x,y)=v\}$ be defined for $u,v\in \C$. Since $B_1$ and $B_2$ meet $\S_\epsilon$ transversely, it follows that there exists $\delta>0$ such that, whenever $|(u,v)|<\delta$, one has that $X_1$ and $X_2$ meet $\S_\epsilon$ transversely, just as in Theorem \ref{theorem:seifert}. We then have a smooth ambient isotopy that carries $K_i$ to $K_i'=X_i\cap \S_\epsilon$, for $i\in \{1,2\}$. Since the linking number is an invariant under ambient isotopies, it follows that $Lk(K_1,K_2)=Lk(K_1',K_2')$. By Sard's Theorem, one can take $(u,v)$ to be a regular value of $F$. Hence, $X_1\cap D_\epsilon$ and $X_2\cap D_\epsilon$ are transverse Seifert surfaces for $K_1'$ and $K_2'$. Since they are in fact one-dimensional complex manifolds, the local intersection numbers are $+1$, and hence $I(X_1,X_2)=|X_1\cap X_2|$. Hence, $Lk(K_1',K_2')=|X_1\cap X_2|$. Lastly, $|X_1\cap X_2|$ is exactly the number of preimages of the regular value $(u,v)$ by $F$, and hence it equals $\deg F$. The proof is now complete, since $C.D=\dim_\C \C\{x,y\}/(f,g)=\deg F$ by Proposition \ref{prop:intnumber} and Corollary D.7 of \cite{Juanjo}. 
\end{proof}

\section{Equisingularity determines the isotopy type of the link}
In the previous section, it has been shown that the understanding of the topological class of a curve and a branch can be performed just by looking at its associated link or knot, respectively; and that a parametric equation for a branch $B$ induces a parametric equation for its associated knot $B\cap \S'_\epsilon$. As we have commented on then, the knots defined by $B\cap \S'_\epsilon$ and $B\cap \S_\epsilon$ are isotopic provided $\epsilon>0$ is small enough. Hence, the achieved equations for the knot $B\cap \S_\epsilon'$ could be seen as equations for the knot given by $B\cap \S_\epsilon$ without further mentions. \\

In what follows, a geometrical description of the resulting knot is about to be performed through the analysis of the provided equations for it. This analysis culminates in the fact that the isotopy type of the knot depends only on the Puiseux characteristic of the branch, and hence on its equisingularity class. Equisingular branches are, therefore, topologically equivalent. \\ 

Lastly, the exponent of contact will flourish in the last theorem of the section, where the result is generalised for curves with several branches. Hence, equisingular curves will be shown to define isotopic links, and thus, to be in the same topological class.\\

With the intention of supporting the underlying intuition, let us consider the particular example given by the branch $x=t^4, y=t^6+t^7$. Focus first on the branch obtained by considering the first term of $y$, namely, $x=t^4, y=t^6$. Its associated knot in $\S_\epsilon'$ can be parametrised by the equations $x=\epsilon e^{2i\theta}, y=\epsilon^{3/2}e^{3i\theta}$. Therefore, the knot lies in the subset of $\C^2$ given by $\{(x,y)\in \C^2 : |x|=\epsilon, |y|=\epsilon^{3/2}\}$, which corresponds to a flat torus. Moreover, for a fixed value of $x$, two values of $y$ are associated with $x$ and, if $\theta\in [0,\pi]$, the values of $x$ rounds the circle $|x|=\epsilon$ once, while $y$ rounds its corresponding circle $3/2$ times. Therefore, the knot can be visualised as in Figure 4: in a vertical circumference (painted red at the first image) draw two antipodal points, which correspond to the values of $y$ that are determined for a given $x$. As $x$ rounds the circle once, the meridian starts moving around the torus and turns back to its starting point. At the same time, the points of $y$ complete a whole lap and a half, which means that not only is the meridian being translated around the torus, but it is also rotating over itself. \\

This rotation movement forces both drawn points to complete $3/2$ laps, and hence when the meridian returns to its original starting point after a lap, the points interchange its positions. The trajectory of these points defines a closed connected curve in the torus, being described in two halves attached due to the spinning movement of the points. \\
\begin{figure}
\centering
\includegraphics[width=0.3\textwidth]{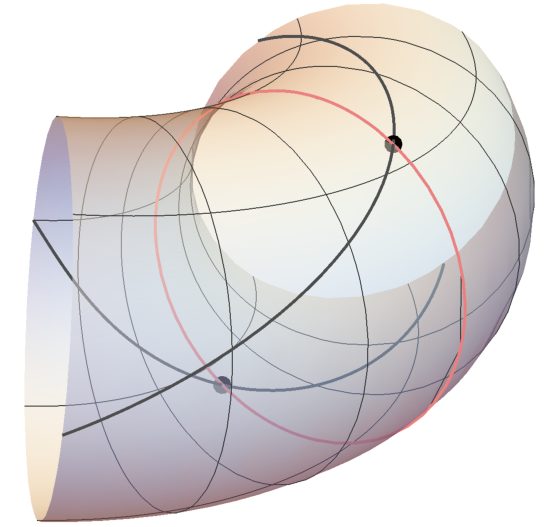}
\includegraphics[width=0.4\textwidth]{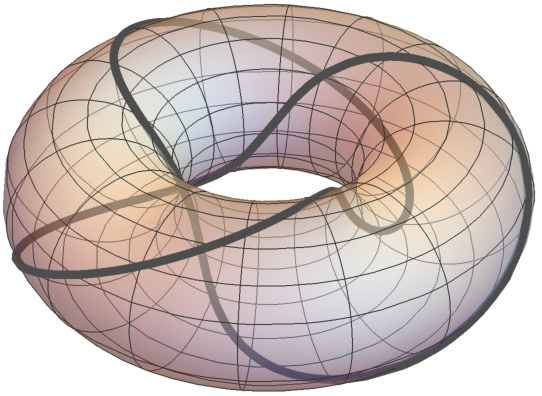}
\caption{Representation of the knot $x=\epsilon e^{2i\theta}, y=\epsilon^{3/2}e^{3i\theta}$ in the torus.}
\end{figure}
In general, the flat torus of radii $\epsilon, \delta>0$ is the set $\{(x,y)\in \C^2 : |x|=\epsilon, |y|=\delta\}=\S^1_{\epsilon}\times \S^1_{\delta}$. A $(p,q)$\textit{-torus knot} is a knot that is isotopic to the one described in the flat torus by the equations $x=\epsilon e^{ip\theta}, y=\delta e^{iq\theta}$, where $p$ and $q$ are nonzero coprime integers. It then follows that the knot described before, given by the branch $x=t^2, y=t^3$, is a $(2,3)$-torus knot. In general, the knot defined by the branch $x=t^p, y=t^q$ for $p<q$ being coprime integers, is a $(p,q)$-torus knot. \\

Turning back to the original chosen branch, given by $x=t^4, y=t^6+t^7$, one has that its knot in $\S_\epsilon'$ is described by the equations
\begin{equation}\label{eq:knot}
    x=\epsilon e^{4i\theta}, y=\epsilon^{3/2}e^{6i\theta}+\epsilon^{7/4} e^{7i\theta}.
\end{equation}
Notice that a fixed value of $x$ determines 4 values of $y$, although the resulting knot is no longer in the flat torus as before. However, since $\epsilon$ can be chosen to be arbitrarily small, it follows that $\epsilon^{7/4}$ tends to 0 faster than $\epsilon^{3/2}$ does. Therefore, one can thicken the knot given by $x=\epsilon e^{2i\theta}, y=\epsilon^{3/2}e^{4i\theta}$ to obtain a surface homeomorphic to a torus encircling it, with radius $\epsilon^{7/4}$. It follows that the knot given by (\ref{eq:knot}) lies in this knotted torus. \\

Indeed, restricting our attention to the meridian of our original flat torus, one can draw both circumferences of radius $\epsilon^{7/4}$ centred at the previous points. Plot antipodal points in both circumferences at the same positions, as shown in Figure 5. As $x$ rounds the circle once (\textit{i.e.}, as the red circle completes a whole translation around the main flat torus), the red circle spins $3/2$ times, and the blue smaller circles spin $7/4$ times. Not only are the blue circles swapped in the process, but their associated blue points are also interchanged. Therefore, the union of the four trajectories described by the blue points form a connected real curve given by (\ref{eq:knot}).
\begin{figure}
    \centering
    \includegraphics[width=0.32\textwidth]{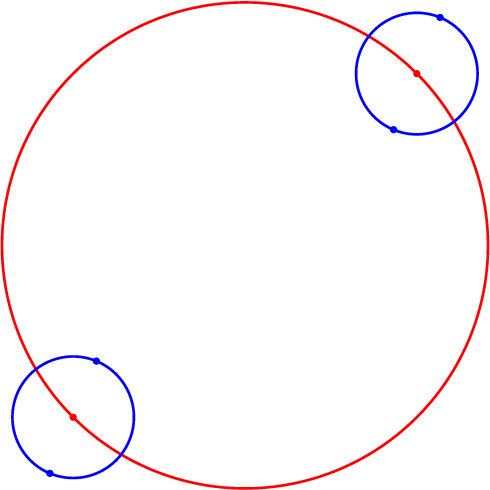}

    \caption{the carousel of the knot given by the branch $x=t^4, y=t^6+t^7$.}
\end{figure}
This kind of construction of nested circles with centers in bigger circles, and all of them spinning at a certain rate is usually referred as a \textit{carousel}. Moreover, the knots that are constructed from a flat torus by iterating a finite number of times the process of taking a torus knot and thickening it to obtain a knotted torus are known as \textit{iterated torus knots}. In the following chapter, this kind of knot construction will be analysed in a more formal treatment. \\

In order to provide a deeper intuition of the concept, let us discuss the first steps that one should follow to construct the knot of a generic branch. Take a branch $(B,0)$ tangent to $y=0$ and with parametric equations $x=t^m, y=\sum_{r=n}^{+\infty}a_rt^r$, where $n>m$. As it is usual, denote its Puiseux characteristic as $(m;\beta_1, \ldots, \beta_g)$, and the auxiliary terms as $e_i=\gcd \{e_{i-1}, \beta_i\}$ for $i\in \{1,\ldots, g\}$, and $e_0=m$. For each $k\geq n$, define the branch $(B_k,0)$ as the one described  by $x=t^m, y_k=\sum_{r=n}^k a_r t^r$ (which is not necessarily an injective function, but it still defines a branch, as it was discussed after Theorem \ref{thrm:alwaysinjective}). Therefore, the knot $K_k=B_k\cap \S'_\epsilon$ can be described by $x=\epsilon e^{im\theta}, y=\sum_{r=n}^k a_r \epsilon^{r/m}e^{ir\theta}$. \\

If $k<\beta_1$, it follows that all the exponents of $t$ that appear in $y_k$ divide $m$, due to the fact that $\beta_1$ is the smallest integer not divisible by $m$. Hence, $y_k$ is a polynomial in $x$. This forces that every value of $x$ corresponds to a unique value of $y$. Therefore, $K_k$ is the knot associated to the carousel formed just by one circumference at each level. Thus, there exists an isotopy carrying $K_k$ to the longitude of the flat torus, forcing $K_k$ to be the unknotted. \\
If $k=\beta_1$, one can write $x=t^m, y_k=y_{k-1}+a_{\beta_1}t\,^{\beta_1}$. Since $\beta_1\nmid m$, then $t\,^{\beta_1}$ cannot be expressed as a natural power of $x$. Then, since $e_1=\gcd \{m,\beta_1\}$, it follows that, for each value of $x$, there corresponds $m/e_1$ equispaced values for $y$. Hence, the knot $K_{\beta_1}$ lies in the torus that is centred at the knot $K_{k-1}$ with radius $|a_{\beta_1}|\epsilon\,^{\beta_1/m}$ and, since $K_{k-1}$ is isotopic to the unknotted, it follows that $K_{\beta_1}$ is isotopic to the torus knot with parameters $(m/e_1, \beta_1/e_1)$. \\

If $\beta_1<k<\beta_2$, one has that the powers in $t$ of order greater than $\beta_1$ are a natural power of $t^{e_1}=x\,^{\beta_1/m}$, as a consequence of being $\beta_2$ the least integer not divisible by $e_1$. Hence, for each value of $x$, there corresponds $m/e_1$ values of $y$ as before. Therefore, the picture of the carousel for $K_k$ is exactly the same as of $K_{\beta_1}$, with the difference that there are now extra levels, but with just one circumference attached at each level. Thus, as in the case $k<\beta_1$, there is an isotopy carrying $K_k$ to $K_{\beta_1}$. Formally, the isotopy is given by 
\begin{equation*}
    x=\epsilon e^{ie_1\theta}, y=\sum_{r=n}^{\beta_1} a_r\epsilon^{r/m}e^{ire_1\theta} + \sum_{r=\beta_1+1}^{k} sa_r \epsilon^{r/m}e^{ire_1\theta/m}
\end{equation*}
taking $0\leq s \leq 1$. \\

The general picture is that, for $\beta_i\leq k < \beta_{i+1}$,  the isotopy type of $K_k$ is the same as $K_{\beta_i}$. In the following theorem, not only are the previous statements formally shown, but it is proven that the Puiseux characteristic determines the isotopy type of the knot. 
\begin{theorem}\label{equisingular_implies_isotopic} The isotopy type of the knot associated to a branch depends only on its Puiseux characteristic. In other words, equisingular branches define isotopic knots. 
\end{theorem}
\begin{proof}
    Let $B$ a representative set-germ of a branch tangent to $y=0$ and parametrised by $x=t^m, y=\sum_{r=n}^{+\infty}a_rt^r$ for $n>m$, and let $K$ be its associated knot in $\S'_\epsilon$, provided $\epsilon>0$ is small enough. Then, $K$ is described by the equations $x=t^m, y=\sum_{r=n}^{+\infty} a_r\epsilon^{r/m}e^{ir\theta}$. We claim that, for $0\leq s \leq 1$, the equations 
    \begin{equation*}
        x=\epsilon e^{im\theta}, y=\sum_{r=n}^{+\infty} s_r a_r \epsilon^{r/m} e^{ir\theta},
    \end{equation*}
    where $s_r=1$ if $r=\beta_q$ and $s_r=s$ otherwise, define an isotopy from $K$ to the knot associated to the branch $x=t^m, y=\sum_{q=1}^g c_q t\,^{\beta_q}$, where $c_q=a_{\beta_q}$. \\
    Notice first that the previous equations are naturally smooth in the variables $(\theta, s)$, and that for a fixed $0\leq s \leq 1$, the equations define an immersion, due to the fact that $\partial x / \partial \theta = im\epsilon e^{im\theta}\neq 0$. It just remains to check the injectivity clause to show that the previous equations define embeddings for a fixed $s$. We should then verify that $\theta$ and $\theta+2\pi k/m$ lead to different values of $y$ unless $m$ divides $k$. As it was commented on before, it is enough to verify that, for each $q\in \N$, the sum of the contributions of the terms between $\beta_q$ and $\beta_{q+1}$ is small compared to the distance of the points in the smallest circle in the carousel of the knot $K_{\beta_q}$. In this case, the carousel construction will define a proper knot without self-intersections.  \\

    Let us restrict first to the case $q=g$, where the contributing terms are infinite. In this case, let $R>0$ be a radius for which $|a_r|R^r\leq C$ for every $r\in \N$ and a certain bound $C>0$. Therefore, 
    \begin{equation*}
        \sum_{r=\beta_g+1}^{+\infty} |a_r|\epsilon^{r/m} \leq \sum_{r=\beta_g+1}^{+\infty} CR^{-r} \epsilon^{r/m}= C \sum_{r=\beta_g+1}^{+\infty} \left(\dfrac{\epsilon^{1/m}}{R}\right)^r = \dfrac{C(\epsilon^{1/m}/R)\,^{\beta_g+1}}{1-\epsilon^{1/m}/R}=\dfrac{C\epsilon^{(\beta_g+1)/m}}{R\,^{\beta_g}(R-\epsilon^{1/m})}
    \end{equation*}
    and, as $\epsilon\to 0$, this tends to $0$ faster than $\epsilon\,^{\beta_g/m}$. If $q<g$, one can bound $\sum_{r=\beta_g+1}^{\beta_{g+1}-1} |a_r| \epsilon^{r/m} \leq \sum_{r=\beta_q+1}^{+\infty} |a_r| \epsilon^{r/m}$, and apply the same procedure of the previous case to grant that this sum has a contribution of order greater than the one of $\epsilon\,^{\beta_q/m}$. Therefore, there exists an $\epsilon_0>0$ such that, for every $0<\epsilon<\epsilon_0$, the sum of the contributions of the terms of order higher than $\beta_q$ is strictly smaller than the distance between two consecutive points in the smallest circle associated to the carousel of $K_{\beta_q}$. Injectivity then follows, proving the desired claim. \\
    
    Therefore, the knot $K$ is isotopic to the one described by the equations 
    \begin{equation*}
        x=\epsilon e^{im\theta}, y=\sum_{q=1}^g a_{\beta_q} \epsilon\,^{\beta_q/m} e^{i\beta_q\theta}.
    \end{equation*}
    Writing $a_{\beta_q}=e^{l_q}$ for some $l_q\in \C$ yields that a further isotopy transforms $c_q$ to $1$ through 
    \begin{equation*}
        x=\epsilon e^{im\theta}, y=\sum_{q=1}^g e^{sl_q} \epsilon\,^{\beta_q/m} e^{i\beta_q\theta}
    \end{equation*}
    for $0\leq s \leq 1$. The injectivity clause follows from the fact that the corresponding branch $x=t^m, y=\sum_{q=1}^g e^{sl_q}t\,^{\beta_q}$ is injective by Theorem \ref{thrm:alwaysinjective} and the fact that $e^{sl_q}\neq 0$ for every $0\leq s \leq 1$. Hence, $K$ is isotopic to the knot given by the equations $x=\epsilon e^{im\theta}, y=\sum_{q=1}^g \epsilon\,^{\beta_q/m}e^{i\beta_q\theta}$, which is precisely the knot associated to the branch $x=t^m, y=\sum_{q=1}^g t\,^{\beta_q}$. Since this branch depends only on the Puiseux characteristic of $(B,0)$, it follows that the isotopy type of the knot is determined by the equisingularity class of $(B,0)$. 
\end{proof}
Therefore, equisingular branches determine isotopic knots. As we have shown in Corollary \ref{corollary:conic}, two branches define isotopic knots if and only if they are topologically equivalent. Hence, equisingular branches are topologically equivalent. The last aim of the section is to extend the previous result to the case of general curves. In order to do so,  the exponent of contact will become crucial.
\begin{theorem}\label{thrm:implication} Let $(C,0)$ and $(\tilde{C},0)$ be equisingular curves, and let $L$ and $\tilde{L}$ be their associated links, respectively. Then, $L$ and $\tilde{L}$ are isotopic. 
\end{theorem}
\begin{proof} The proof is just sketched, since a complete and formal proof could be tedious and technical, and our intention is to provide a dynamic exposition of the topic. For a further reference, see 5.3.2 of \cite{Wall}. If the curves $(C,0)$ and $(\tilde{C},0)$ are equisingular, in particular they have the same number of branches, say $n$. Write $(B_1,0), \ldots, (B_n,0)$ for the branches of $(C,0)$, and $(\tilde{B}_1,0), \ldots, (\tilde{B}_n,0)$ for the ones of $(\tilde{C},0)$. From the equisingularity condition follows that, up to a relabelling of the branches of $(\tilde{C},0)$, the branches $(B_i,0)$ and $(\tilde{B}_i,0)$ share the same Puiseux characteristic for $i\in \{1,\ldots, n\}$, and that $\cont (B_i,B_j)=\cont (\tilde{B}_i, \tilde{B}_j)$ for $1\leq i\neq j \leq n$. \\
Choose coordinates in $(\C^2,0)$ in such a way that $x=0$ is not the tangent line of the branches of $(C,0)$ and $(\tilde{C},0)$. It turns out that one can choose pro-branches $\gamma_i$ for $(B_i,0)$ simultaneously for $i\in \{1,\dots, n\}$, in such a way that $\cont (\gamma_i, \gamma_j)=\cont (B_i,B_j)$, for $1\leq i \neq j \leq n$, and the same for the branches of $(\tilde{C},0)$. This condition is vacuous if $n=1$, and trivial for $n=2$. The verification of the claim for a generic natural number $n$ can be tackled through an induction reasoning, and the details can be found in the provided reference. Write $\gamma_i(x)=\sum_{r\geq 1}a_{r,i} x^{r}$ for $1\leq i \leq n$, and similarly for $\tilde{\gamma_i}$, with coefficients $\tilde{a}_{r,i}$. We claim that deformations $a_{r,i}(s)$ with $s\in [0,1]$ can be constructed with the property that $a_{r,i}(0)=a_{r,i}$, $a_{r,i}(1)=\tilde{a}_{r,i}$, and where:
\begin{enumerate}[(1)]
    \item The Puiseux characteristic of the branch defined by $(B_i(s),0)$ with the coefficients $a_{r,i}(s)$ is independent of $s\in [0,1]$.
    \item The exponents of contact $\cont (B_i(s), B_j(s))$ do not depend on $s\in [0,1]$. 
\end{enumerate}
Let us assume first that such deformations exist and see how could one construct the desired isotopy between the links $L$ and $\tilde{L}$. For convenience, let us define the curve $(C(s),0)$ as the one described by the branches $(B_i(s),0)$, for $i\in \{1,\dots, n\}$. These coefficients $a_{r,i}(s)$ provide a deformation from $(C,0)$ to $(\tilde{C},0)$ through the curves $(C(s),0)$, where the Puiseux characteristic of the branches and the exponents of contact are preserved. Intersecting with $\S_\epsilon$ yields a series of parametric equations for the arising links of $(C(s),0)$. More precisely, one obtains a map $f:X\times I\rightarrow \S^3$, where $X$ is a disjoint union of $n$ copies of $\S^1$, and for $s\in I=[0,1]$, the associated functions $f_s$ satisfy that $f_s(X)$ is the link associated to $(C(s),0)$. In particular, for a fixed $s\in [0,1]$, the map $f_s$ is immersive by the same justification provided in Theorem \ref{equisingular_implies_isotopic}, since $f_s$ is immersive in each connected component of $X$. Moreover, the maps $f_s$ are injective as a consequence of (1) and (2). This mainly follows from the provided description of the knots in the previous theorem. Indeed, if $(x,y)\in \C^2$ are the chosen coordinates, then, for a fixed value of $x$, there are finitely many values of $y$ with the property that $(x,y)$ lies in the link generated by $(C(s),0)$. Now, as it was discussed in the previous theorem, the values of $y$ that are associated to the same branch are all different. Now, since the exponent of contact $\cont (B_i(s), B_j(s))$ does not depend on $s$, it follows from Lemma \ref{lemma:contactinvariant} that the orders of the distances between the knots in the link of $(C(s),0)$ are constant in $s$. Therefore, the values of $y$ that correspond to different branches are also distinct. This shows that $f_s$ is indeed an embedding for every $s\in [0,1]$, forcing $f$ to be an isotopy between $L$ and $\tilde{L}$. \\

It just remains to construct the deformations $a_{r,i}(s)$ satisfying the above-mentioned conditions. In order to provide a precise and clear definition for the $a_{r,i}(s)$, a generic expression is omitted. The deformations will then be  defined through an induction process. Notice that, for $n=1$, Theorem \ref{equisingular_implies_isotopic} provides a method to satisfy (1), since the given equations keep the Puiseux characteristic constant in $s\in [0,1]$, and (2) is vacuous for a single branch. Assume that $a_{r,i}(s)$ has been defined for $1\leq i \leq n-1$, and let us construct the remaining $a_{r,n}(s)$. Write
\begin{equation*}
    \kappa_n=\max \, \{\cont (B_i, B_n): 1\leq i \leq n-1\},
\end{equation*}
and choose $j\in \{1,\dots, n-1\}$ satisfying $\cont (B_j, B_n)=\kappa_n$. define then $a_{r,n}(s)=a_{r,j}(s)$ for $r<\kappa_n$. \\ 
Let us establish some notation for a more concise explanation. If $(B,0)$ is a branch with Puiseux characteristic $(m;\beta_1, \dots, \beta_g)$, then a \textit{characteristic exponent} is a value $r=\alpha_q=\beta_q/m$ for some $q\in {1,\dots, g}$. The values $\alpha_{q-1}<r<\alpha_q$ that satisfy $mr/e_{q-1}\in \Z$ are precisely the ones whose coefficients could be altered in the pro-branch without changing the Puiseux characteristic. This justifies to call them \textit{free exponents}. The remaining ones will be referred as the \textit{forbidden exponents}. It then follows that a generic pro-branch has the Puiseux characteristic of $(B,0)$ if and only if the coefficients associated with powers of $x$ with forbidden exponents are equal to 0, and the ones associated with the characteristic exponents are nonzero. This notation will help to construct the deformations for $r\geq \kappa_n$. \\

If $r=\kappa_n$, there are two main possibilities. If $r$ were a forbidden exponent of $(B_n,0)$, it would then be a characteristic one for $(B_j,0)$. Therefore, we have the chance to let $a_{r,n}(s)=0$, since $a_{j,n}(s)\neq 0$ for every $s\in [0,1]$. If $r=\kappa_n$ is not forbidden, it can be a free or a characteristic exponent. In this case, one has to make sure that the values of $a_{r,n}(s)$ and $a_{j,n}(s)$ do not coincide for $s\in [0,1]$. Since this holds for $s=0,1$, we are able to take $\epsilon>0$ small enough to have a distance $d>0$ satisfying that $a_{r,j}(s)\in B(a_{r,j},d)$ for $s\in [0,\epsilon]$, $a_{r,j}(s)\in B(\tilde{a}_{r,j}, d)$ for $s\in [1-\epsilon, 1]$, and that $a_{r,n}, \tilde{a}_{r,n}$ do not lie in the previous balls. Let $c>0$ be a value greater than $\max_{s\in [0,1]}|a_{r,j}(s)|$. Since, for $s\in [0,\epsilon]$, the values $a_{r,j}(s)$ are contained in the small ball mentioned before, define $a_{r,n}(s)$ for $s\in [0,\epsilon]$ to be a path joining $a_{r,j}$ with $c$, and avoiding the ball $B(a_{j,n},d)$. For $s\in [\epsilon, 1-\epsilon]$, define it to be constant at $c$. Lastly, for $s\in [1-\epsilon, 1]$, choose a path from $c$ to $\tilde{a}_{r,n}$ avoiding the ball $B(\tilde{a}_{r,j},d)$. \\

\begin{figure}[h]
    \centering
    \includegraphics[scale=0.35]{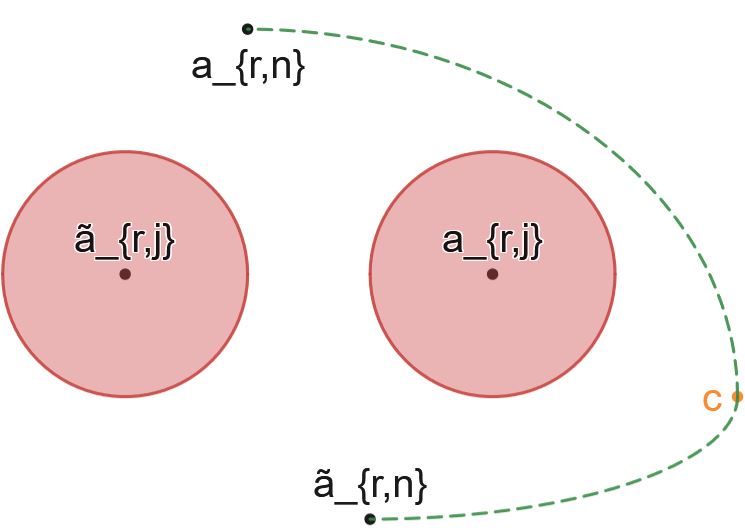}

    \caption{Representation of $a_{r,n}(s)$ if 
    $r=\kappa_n$ is a forbidden exponent}
\end{figure}

This gives a well-defined path connecting $a_{r,n}$ and $\tilde{a}_{r,n}$ that does not coincide with $a_{r,j}(s)$ for any value of $s\in [0,1]$. If $r$ is a characteristic exponent, the above-chosen paths should also avoid the origin. Lastly, if $r>\kappa_n$, the exponents of contact with all the other branches are already fixed, and we are able to join $a_{r,n}$ with $\tilde{a}_{r,n}$ freely through a path (that cannot pass through the origin whenever $r$ is characteristic). This completes the construction of the desired deformations. 
\end{proof}
Theorem \ref{thrm:implication} shows that equisingular curves are topologically equivalent. The remaining part of the chapter will be devoted to prove the converse statement. Namely, it will be shown that the isotopy type of a link associated to a curve $(C,0)$ completely determines the equisingularity class of $(C,0)$. It will then follow that topological equivalence and equisingularity define the same equivalence among plane curves.
\section{Knots of branches as iterated torus knots}
This section is devoted to analyse the knots that appear as the intersection of a branch with a small enough $3$-dimensional sphere, with the intention of extracting the Puiseux characteristic by means of isotopy invariants of the knot. Consider a generic branch that is not tangent to $x=0$ and that generates the knot $K$. Let $(m; \beta_1, \dots, \beta_g)$ be the Puiseux characteristic of the branch, and let $e_0, \dots, e_g$ be its auxiliary terms. The previous section shows that $K$ is isotopic to the knot generated by the branch $(B,0)$ given by  $x=t^m, y=\sum_{q=1}^g t\,^{\beta_q}$. Hence, $K$ can be identified with the knot associated to $(B,0)$. \\ 

Let $(B_q,0)$ be the branch defined by $x=t^m, y=\sum_{i=1}^q t\,^{\beta_i}$, and let $K_q$ be its associated knot. Then, as we have commented on in the previous chapter, the knot $K_q$ is constructed from the knot $K_{q-1}$ by what we will call a \textit{cabling process}: as a torus knot in the knotted torus arising from thickening $K_{q-1}$. \\

In order to formalise this concept, let $K$ be a knot in $\S^3$ given by an embedding $f:\S^1\rightarrow \S^3$, and consider a small enough compact tubular neighbourhood of $K$ given by the image of an embedding $F:\S^1\times D^2\rightarrow \S^3$ with the property that $F(x,0)=f(x)$, for $x\in \S^1$. It is a well-known fact that the image of $F$ is a deformation retract of $K$ and that it is unique up to an isotopy in $\S^3$. Let us focus our attention on the image $T=F(\S^1\times \S^1)$, which is a manifold intrinsically diffeomorphic to a torus, but extrinsically knotted as $K$ is. Now, any torus knot is sent via $F$ to a knot $K'$ in the knotted torus. In this case, we will say that $K'$ is a \textit{cable knot} about $K$. Notice that, although the tubular neighbourhood is unique up to an isotopy, one may have different options for the embedding $F$. Hence, the image of the standard longitude and meridian of the flat torus depend on $F$. The meridian of $T$ will then be chosen as a curve in $T$ with trivial homology in $F(\S^1\times D^2)$. This determines the meridian up to the orientation and up to isotopies in $T$. The choice of a longitude in $T$ will be performed through a curve that has zero linking number with the original knot $K$. With this, a study of the homology class of the knot $K'$ constructed in $T$ expressed as an integer linear combination of the homology classes of the longitude and meridian show that the knot $K'$ is determined up to isotopy by some tuple of natural numbers $(p,q)$, where $p$ is the number of crossings of $K'$ with the meridian, and $q$ is the linking number of $K$ and $K'$. In this case, we say that $K'$ is a cable knot about $K$ with \textit{cabling invariant} $(p,q)$. For more details on this construction, see \cite{Wall}, chapter 5.4. \\

In terms of the above-mentioned concepts, it follows that $K_q$ is a cable knot about $K_{q-1}$. In what follows, the cabling invariant is determined. 
\begin{theorem} The first cabling invariant number is given by $e_{q-1}/e_q$.
\end{theorem}
\begin{proof} Let us first obtain an explicit parametrisation of the knotted torus around $K_{q-1}$. Since $K_{q-1}$ is given by the equations $x=\epsilon e^{im\theta}, y=\sum_{r=1}^{q-1} \epsilon\,^{\beta_r/m} e^{i\beta_r\theta}$, then $$x=\epsilon e^{im\theta}, \left|y-\sum_{r=1}^{q-1} \epsilon\,^{\beta_r/m}  e^{i\beta_r\theta}\right|=\epsilon\,^{\beta_q/m}$$ defines a knotted torus of radius $\epsilon\,^{\beta_q/m}$ around $K_{q-1}$, for $\epsilon>0$ small enough. Since the parametric equations of the knot are not injective in general, it follows that the equation
\begin{equation*}
    x=\epsilon e^{im \theta/e_{q-1}}, y=\sum_{r=1}^{q-1}  \epsilon\,^{\beta_r/m} e^{i\beta_r\theta/e_{q-1}}+\epsilon\,^{\beta_q/m} e^{i\phi}
\end{equation*}
provides an injective parametrisation of the torus. Now, notice that the knot $K_q$ can be easily parametrised just by taking $(\theta, \phi)=(e_{q-1}/e_q \psi, \beta_q/e_q \psi)$. Since $e_{q-1}/e_q$ is a natural number, it is then clear that the knot $K_q$ intersects $e_{q-1}/e_q$ times the meridian of the torus. This is precisely the definition of the first cabling invariant. 
\end{proof}
\noindent In the following result, the second cabling invariant is determined. 
\begin{theorem} The second cabling invariant number is given by $Lk(K_q,K_{q-1})=\barbeta_q/e_q\,$.
\end{theorem}
\begin{proof} Notice first that  $Lk(K_q,K_{q-1})=B_q.B_{q-1}$ as an application of Corollary \ref{cor:linking}. Notice that the defining equations for $B_q$ are not injective in general, but the equations $x=t^{m/e_q}, y=\sum_{i=1}^q t\,^{\beta_i/e_q}$ indeed define injective parametric equations for the branch $(B_q,0)$ due to the fact that $e_q=\gcd \{m,\beta_1, \dots, \beta_q\}$. Hence, $m(B_q)=m/e_q$ and $m(B_{q-1})=m/e_{q-1}$. Moreover, $\kappa=\cont (B_{q}, B_{q-1})=\beta_q/m$, since pro-branches for $(B_{q-1},0)$ and $(B_q,0)$ are, respectively, $y=\sum_{i=1}^{q-1} x\,^{\beta_i/m} \text{ and } y=\sum_{i=1}^{q} x\,^{\beta_i/m}$. Now, the Puiseux characteristic of $(B_q,0)$ is $(m/e_q; \beta_1/e_q, \ldots, \beta_q/e_q)$, and hence if $e_i'$ and $\barbeta_i'$ are the Puiseux terms associated to the branch $(B_q,0)$, it follows from its definitions and by a straightforward induction that $e_i'=e_i/e_q$ and $\barbeta_i'=\barbeta_i/e_q$, for $i\in \{0,\dots,q\}$. Since $\kappa$ lies in the interval $\beta_{q-1}/m<\kappa\leq \beta_q/m$, it follows by equation (\ref{eq:intnumber}) that
    \begin{equation*}
        B_{q}.B_{q-1}=e_{q-1}'\dfrac{m(B_{q-1})}{m(B_q)}\Bigl(\barbeta_q'-\beta_q'+m(B_{q})\kappa\Bigr)=\dfrac{e_{q-1}}{e_q}\dfrac{m/e_{q-1}}{m/e_{q}}\left( \dfrac{\barbeta_q}{e_q} - \dfrac{\beta_q}{e_q} + \dfrac{m}{e_q} \dfrac{\beta_q}{m}\right)=\dfrac{\barbeta_q}{e_q},
    \end{equation*} 
    and the claim follows.
\end{proof}
In conclusion, $K_{q}$ is a cable about $K_{q-1}$, with cabling invariant $(e_{q-1}/e_q, \barbeta_q/e_q)$. Hence, the construction of the knot $K$ arising from the intersection of a branch with a small enough 3-sphere can be described through an iterated cabling process, starting from the unknotted $K_0$. Therefore, there is a finite sequence of knots $K_0, K_1, \ldots, K_g=K$, starting from the unknotted, and where each $K_q$ is a cable knot about $K_{q-1}$. Hence, $K$ is what is known as an \textit{iterated torus knot}, determined by its cabling invariants $\{ (e_{q-1}/e_q, \barbeta_q/e_q)\}_{q=1}^g$.
\section{The isotopy type of the link determines the equisingularity class of the curve}
The aim of this last section is to show the converse statement of Theorem \ref{thrm:implication}, namely that, if two curves have isotopic associated links, then the curves are equisingular. As it is routine, the result will be shown first for branches, and then it will be extended to general curves. This, in particular, shows that the arithmetic tools presented in the second chapter, namely the Puiseux characteristic and the semigroup of a branch, are topological invariants of the curve. The main tool that will be used along this section to distinguish from algebraic knots is the Alexander polynomial, a well-known isotopy invariant for knots that assigns to each knot a polynomial of a single variable, which is well-defined up to a multiplication with a factor of the form $\pm t^k$. \\

The construction of this polynomial could be performed through different angles, some being more computational or more theoretical procedures. The aim of this project does not lie in the study of knot-theoretical invariants, but it is to apply them to achieve the classification of plane curve singularities. Thus, the definition of the Alexander polynomial will be merely sketched, leaving the details to further references. \\

For a full discussion on the definition of the Alexander polynomial of a knot, we refer the reader to \cite{burde} for a theoretical approach, and \cite{murasugi} for a more computational and practical development. The main idea is to start from a knot $K$ and a Seifert surface $\S$ for it, and to create a certain amount of closed curves in $\S$ that arise naturally in the Seifert algorithm for constructing the surface $\S$. If $\alpha_1, \ldots, \alpha_n$ are the above-mentioned closed curves on $\S$, one constructs the \textit{Seifert matrix of }$K$ with respect to $\S$ as the matrix $M$ whose $(i,j)$ entry is the linking number of $\alpha_i$ and a slight deformation of $\alpha_j$, in a precise way explained in \cite{murasugi}. This matrix depends on the chosen surface $\S$, and hence it is far from being an isotopy invariant for $K$. Nevertheless, it can be shown that the polynomial 
\begin{equation*}
    \Delta_K(t)=\det (M-tM^T),
\end{equation*} 
defined up to multiplication by a factor $\pm t^k$, is indeed independent on $S$ and on the isotopy type of $K$. This polynomial, defined as the \textit{Alexander polynomial} of $K$, is therefore a knot invariant. A similar process could be performed to construct the Alexander polynomial for a link. However, for the sake of brevity, the details are omitted. \\

For our purposes, some relevant facts concerning the Alexander polynomial of a knot are going to be stated. The first of them is that the Alexander polynomial of the unknotted is given by the constant polynomial $1$. On the other hand, it is not difficult to check that the Alexander polynomial of the $(p,q)$-torus knot is
\begin{equation*}
    \Delta_{(p,q)}(t)= \dfrac{(t^{pq}-1)(t-1)}{(t^p-1)(t^q-1)},
\end{equation*}
and that this is, in fact, a polynomial (\textit{i.e.}, the denominator divides the numerator in the ring of polynomials $\C[t]$). For a reference, see 9.15 of \cite{burde}. Moreover, it turns out that the Alexander polynomial of a knot arising from a cabling process can be easily obtained, as the following proposition holds. 
\begin{proposition}\label{prop:alexander} If $K'$ is a cable knot about $K$ with cabling invariant $(p,q)$, then 
\begin{equation*}
    \Delta_{K'}(t)=\Delta_K(t^p)\Delta_{(p,q)}(t).
\end{equation*}
\end{proposition}
\noindent We refer the reader for a proof to \cite{eisenbud}, in Theorem 5.3. The deduction of Proposition \ref{prop:alexander} from the referred theorem is given below Theorem 9.7.2 of \cite{Wall}.  \\ 

Turning back to the case of the knot $K$ associated to a branch $(B,0)$, it has been shown in the previous chapter the existence of a finite sequence of knots $K^0, K^1, \dots, K^g$, where $K^0$ is the unknotted, $K^g=K$ and, for every $q\in \{1,\dots,g\}$, $K^q$ is a cable knot about $K^{q-1}$ with cabling invariant $(e_{q-1}/e_q, \barbeta_q/e_q)$. This gives us the opportunity to obtain the Alexander polynomial of $K$ inductively. \\ 

In order to give a clear and concise exposition, it is mandatory to introduce a special notation. Let us define $\S(n)$ to be the \textit{symbol} of the polynomial $t^n-1$, just as a method to shorten notation. For a product of polynomials of the previous kind, an additive notation in the symbols will be used. More precisely, define the \textit{symbol} of the expression $\prod_{n=1}^N (t^n-1)^{a_n}$ as the formal sum $\sum_{n=1}^N a_n \S(n)$, where $a_n\in \Z$ for every $n\in \{1,\dots, N\}$. As a special case, consider the empty sum $0$ (with $N=0$) to be the symbol of the constant polynomial $1$. This provides a more compact way to refer some polynomials, as for instance the symbol of $\Delta_{(p,q)}(t)$ is given by $\S(pq)+\S(1)-\S(p)-\S(q)$. Notice that not every polynomial admits a symbol of this kind, but the ones we will work with will be shown to admit it. Notice that the symbol completely determines its associated polynomial:
\begin{proposition} The symbol completely determines the expressions of the form $\prod_{n=1}^N (t^n-1)^{a_n}$, for $a_n\in \Z$ and $N\in \N$. 
\end{proposition}
\begin{proof} One has to show that the equality $\prod_{n=1}^N (t^n-1)^{a_n}=\prod_{n=1}^N (t^n-1)^{b_n}$ for every $t\in \C$ forces $a_n=b_n$ for every $n\in \{1,\dots,N\}$. The verification of this statement will be performed through an induction on $N$. For $N=0$ the result trivially holds, since both products are empty and corresponding to the constant polynomial $1$. Assume that the result holds for $N-1$. Notice that, after cancelling the powers in $t^N-1$ in both sides, and possibly interchanging the term from the left to the right hand side, one can assume that $a_N\geq b_N=0$. Assume that $a_N>0$. In this case, the substitution of $t=e^{2\pi i /N}$ in the left hand side is $0$, but it is nonzero in the right one since $b_N=0$. Then, a contradiction has reached, forcing $a_N=0$. The result follows by the inductive hypothesis.
\end{proof}
\noindent In the following theorem, the Alexander polynomial of $K$ is obtained through its symbol. 
\begin{theorem} Let $q\in \{0,\dots,g\}$. Then, the Alexander polynomial of $K^q$ admits a symbol, and it is equal to
\begin{equation*}
    \sum_{i=1}^q \S\left(\dfrac{e_{i-1}\barbeta_i}{e_i e_q}\right) - \sum_{i=1}^q \S\left( \dfrac{\barbeta_i}{e_q}\right) -\S\left(\dfrac{m}{e_q}\right) + \S(1).
\end{equation*}
Moreover, the Alexander polynomial of $K$ also admits a symbol, and it is given by 
\begin{equation*}
    \sum_{i=1}^g \S\left(\dfrac{e_{i-1}\barbeta_i}{e_i}\right) - \sum_{i=1}^g \S\left( \barbeta_i\right) -\S\left(m\right) + \S(1).
\end{equation*}
\end{theorem}
\begin{proof} It is clear that the second statement follows from the first one, due to the fact that $K^g=K$ and $e_g=1$. For the first claim, notice that $K^0$ has as Alexander polynomial the constant $1$ since it is the unknotted, and the formula for $q=0$ yields the symbol $-\S(m/e_0)+\S(1)=-\S(1)+\S(1)=0$, corresponding to the polynomial $1$. Hence, the result is true for $q=0$. Assume that the result holds for $q$, and let us show that this implies the case associated to $q+1$. Notice that the Alexander polynomial of $K^{q+1}$ is given by Proposition \ref{prop:alexander} through the relation 
\begin{equation*}
    \Delta_{K^{q+1}}(t)=\Delta_{K^q}(t^{e_{q}/e_{q+1}})\Delta_{(e_{q}/e_{q+1}, \,\barbeta_{q+1}/e_{q+1})}(t).
\end{equation*}
Now, if a polynomial admits a symbol $\S(n)$ and one replaces $t$ by $t^k$, the symbol of the composition is just $\S(kn)$. Moreover, the symbol of a product of polynomials is the formal sum of the symbols. These claims, together with the inductive hypothesis, yield that the symbol of $\Delta_{K^{q+1}}(t)$ is given by 
\begin{equation*}
    \sum_{i=1}^q \S\left(\dfrac{e_{i-1}\barbeta_i}{e_i e_{q+1}}\right) - \sum_{i=1}^q \S\left( \dfrac{\barbeta_i}{e_{q+1}}\right) -\S\left(\dfrac{m}{e_{q+1}}\right) + \S\left(\dfrac{e_q}{e_{q+1}}\right) + \S\left(\dfrac{e_q}{e_{q+1}}\dfrac{\barbeta_{q+1}}{e_{q+1}}\right)- \S\left(\dfrac{e_q}{e_{q+1}}\right)-\S\left(\dfrac{\barbeta_{q+1}}{e_{q+1}}\right)+\S(1).
\end{equation*}
As one may notice, there are two equal terms that cancel out, namely, the given by $\S(e_q/e_{q+1})$. Moreover, the first and third symbols that come from the torus knot are the cases $i=q+1$ of the first and second sums, respectively. This provides exactly the expected formula of the symbol of $\Delta_{K^{q+1}}(t)$ claimed in the statement. The result follows by induction. 
\end{proof}
The previous theorem shows that the Alexander polynomial of $K$ is given by
\begin{equation*}
    \Delta_K(t)=\dfrac{(t-1)\prod_{i=1}^g (t^{e_{i-1}\barbeta_i/e_i}-1)}{(t^m-1)\prod_{i=1}^g (t\,^{\barbeta_i}-1)}.
\end{equation*}
As a straightforward consequence, the desired converse result follows:
\begin{corollary}\label{isotopy_implies_equisingular} The isotopy type of the knot $K$ determines the Puiseux characteristic of the branch $(B,0)$ that defines it. 
\end{corollary}
\begin{proof} The previous result shows that the symbol of the Alexander polynomial of $K$ is an expression regarding some of the arithmetic invariants of $(B,0)$ that have been studied. Moreover, the given formula admits no further cancellations (\textit{i.e.}, there are no positive terms that can be cancelled with the negative ones). Indeed, if $e_{i-1}\barbeta_i / e_i=\barbeta_j$ for some $i,j\in \{1,\dots,g\}$, it would then follow that $\barbeta_i$ properly divides $\barbeta_j$, since $e_{i-1}/e_i$ is a nonzero natural number. However, this contradicts \ref{theorem:semigroup}, where it is stated that the $\barbeta_i$ are a minimal set of generators of $\S(B)$. Hence, for a knot $K$, the negative terms appearing in the symbol of its Alexander polynomial are precisely the increasing sequence given by $\barbeta_0, \dots, \barbeta_g$ associated to $(B,0)$. By Theorem \ref{thrm:semigroup_determines_puiseux}, it follows that these coefficients determine the Puiseux characteristic of the branch. 
\end{proof}
\begin{remark} The polynomials of the form $t^n-1$ factor through a special kind of polynomials, which are cyclotomic polynomials. For a fixed $d\in \N$, the \textit{cyclotomic polynomial} $\Phi_d$ is the monic polynomial whose roots are the primitive $d$-th roots of the unity. It is a well-known fact that $t^n-1=\prod_{d|n} \Phi_d(t)$. This gives us the opportunity to write the above-mentioned fractional expression for the Alexander polynomial as a product of cyclotomic polynomials through cancellations in the denominator. 
\end{remark}
\noindent In the following example, this simplification procedure will be performed in a particular case.
\begin{example} Let us consider the branch of example \ref{example}, given by $x=t^4, y=t^6+t^7$, and let $K$ be its associated knot. As we have checked, $e_0=4, e_1=2, e_2=1$, and $\barbeta_0=4, \barbeta_1=6, \barbeta_2=13$. Hence, $\Delta_K(t)$ has the symbol $\S(26)+ \S(12)+\S(1)- \S(13)-\S(6)-\S(4)$, and thus
\begin{equation*}
    \Delta_K(t)=\dfrac{(t^{26}-1)(t^{12}-1)(t-1)}{(t^{13}-1)(t^6-1)(t^4-1)}.
\end{equation*}
After the factorisation of the involved terms through cyclotomic polynomials is performed, all the terms in the denominator cancel out (as it was expected to happen), leaving that 
\begin{equation*}
    \Delta_K(t)=\Phi_{26}(t)\Phi_{12}(t)=(t^{12}-t^{11}+t^{10}-t^9+t^8-t^7+t^6-t^5+t^4-t^3+t^2-t+1)(t^4-t^2+1)
\end{equation*}
is the expression of $\Delta_K(t)$ as a product of cyclotomic polynomials. \\
\end{example}
Corollary \ref{isotopy_implies_equisingular} completes the desired classification for plane branches, as the following result encapsulates:
\begin{theorem} Let $(B,0)$ and $(\tilde{B},0)$ be branches, and let $K$ and $\tilde{K}$ be their associated knots. Then, the following are equivalent:
\begin{enumerate}[(i)]
    \item $(B,0)$ and $(\tilde{B},0)$ are equisingular.
    \item $K$ and $\tilde{K}$ are isotopic.
    \item The Alexander polynomials of $K$ and $\tilde{K}$ are equal.
    \item $(B,0)$ and $(\tilde{B},0)$ are topologically equivalent. 
\end{enumerate}
\end{theorem}
\begin{proof} $(i)$ implies $(ii)$ is given by Theorem \ref{equisingular_implies_isotopic}, and $(ii)$ implies $(iii)$ is a consequence of being the Alexander polynomial an isotopy invariant of knots. Now, Theorem \ref{isotopy_implies_equisingular} shows that $(iii)$ implies $(i)$. Lastly, $(ii)$ if and only if $(iv)$ is shown in Corollary \ref{corollary:conic}. 
\end{proof}
This concludes the desired classification of branches, and the following result generalises it to curves with several branches.  
\begin{theorem} Let $(C,0)$ and $(\tilde{C},0)$ be branches, and let $L$ and $\tilde{L}$ be their associated links. Then, the following are equivalent:
    \begin{enumerate}[(i)]
        \item $(C,0)$ and $(\tilde{C},0)$ are equisingular.
        \item $L$ and $\tilde{L}$ are isotopic.
        \item $(C,0)$ and $(\tilde{C},0)$ are topologically equivalent. 
    \end{enumerate}
\end{theorem}
\begin{proof} Notice first that Corollary \ref{corollary:conic} shows that $(ii)$ holds if and only if $(iii)$ does, as it happened in the case of branches. Now, $(i)$ implies $(ii)$ is shown in Theorem \ref{thrm:implication}. Let us check that $(ii)$ implies $(i)$. If $L$ and $\tilde{L}$ are isotopic, there exists a diffeomorphism $H:\S^3\times I\rightarrow \S^3\times I$ of the form $H(y,t)=(h_t(y),t)$, satisfying that $h_0=id_{\S^3}$ and that $h_1(L)=\tilde{L}$. Denote the connected components of $L$ as $K_1, \dots, K_n$, and the ones of $\tilde{L}$ as $\tilde{K}_1, \dots, \tilde{K}_n$. Therefore, the diffeomorphism $H$ preserves these connected components, and yields a bijection between the components of $L$ and the ones of $\tilde{L}$. In particular, up to a relabelling, one has that $K_i$ and $\tilde{K}_i$ are isotopic knots through $H$. Thus, the associated branches that generate these knots are equisingular. Furthermore, the action of such a diffeomorphism $H$ preserves the linking numbers, and hence $Lk(K_i, K_j)=Lk(\tilde{K}_i, \tilde{K}_j)$ for every $i\neq j\in \{1,\dots, n\}$. Since the linking number of knots coincides with the intersection number between their associated branches, it follows that the previous induces a bijective correspondence on the branches of $(C,0)$ and $(\tilde{C},0)$ where the equisingularity class of the branches is preserved, as well as the intersection numbers between them. This forces $(C,0)$ and $(\tilde{C},0)$ to be equisingular. 
\end{proof}
It can be shown that the Alexander polynomial of a link is a complete invariant for the classification of plane curves under the topological equivalence. However, the proof of this is far from being straightforward. In fact, it was proved 50 years later by Yamamoto in \cite{yamamoto}. \\

This completes the topological classification of plane curve singularities, and yields that the topological nature of a curve-germ is completely determined by the equisingularity class of the curve, which can, in turn, be described by some arithmetic data: the Puiseux characteristic of its branches, and the intersection numbers between them. Hence, the topological information that a curve embodies can be completely described by means of some numerical parameters that are extracted from the defining equation of the curve. This concludes the exposition. \\

The given definition of equivalence is topological. Hence, it overlooks some geometrical aspects of curves regarding their analytic structures. As it was commented on before, one could consider curves to be analytically equivalent if one curve can be reduced to the other through a holomorphic change of coordinates in $\C^2$. This equivalence gives a much finer relation, and a complete analytic classification deals with the differences that two equisingular curves may have in an analytic sense. The complete analytic classification was achieved in 2021, in \cite{analyticclassification}. \\

The possible generalizations for the curve concept are diverse. One could study, for instance, spaces in $(\C^{n+1},0)$ defined by a nontrivial equation $f(x_1,\ldots, x_{n+1})=0$, yielding the so-called \textit{complex hypersurfaces}. In this context, it is possible to generalise the equisingularity relation to hypersurfaces through \textit{Whitney stratifications}. However, the classifications through the topological and equisingular equivalences yield no longer the same classes for $n>1$. This makes clear that the case $n=1$, \textit{i.e.}, the one studied along this project, is a privileged setting. Indeed, although some results that have appeared in the project admit a general version in this wider context, it turns out that branches of complex hypersurfaces do not admit, in general, parametric equations. This is the key difference between curves in $\C^2$ and hypersurfaces in $\C^{n+1}$, and also the main ingredient that enriches the theory of complex plane curve singularities. 

\thispagestyle{plain}
\printbibliography
\addcontentsline{toc}{chapter}{References}

\end{document}